\documentclass[a4paper,11pt,reqno]{smfart}
\usepackage{amssymb,amsmath,mathrsfs,graphics,graphicx,mathptm,float}
\usepackage[T1]{fontenc}
\usepackage[english, francais]{babel}
\usepackage[all]{xy}
\theoremstyle{plain}
\newtheorem{thm}{Th\'eor\`eme}[section]
\newtheorem{pro}[thm]{Proposition}
\newtheorem{lem}[thm]{Lemme}
\newtheorem{cor}[thm]{Corollaire}
\newtheorem{proposition}{Proposition}

\theoremstyle{definition}
\newtheorem*{defi}{D\'efinition}
\newtheorem*{defis}{D\'efinitions}
\newtheorem{eg}[thm]{Exemple}
\newtheorem*{ap}{Applications de la Proposition \ref{calculformel}}
\newtheorem{prob}{Probl\`eme}
\newtheorem{probs}[prob]{Probl\`emes}

\newtheorem{rem}[thm]{Remarque}
\newtheorem{rems}[thm]{Remarques}

\usepackage[usenames,dvipsnames]{color}
\usepackage[colorlinks, linktocpage, citecolor = blue, linkcolor = blue]{hyperref}

\def\og{\leavevmode\raise.3ex\hbox{$\scriptscriptstyle\langle\!\langle$~}}
\def\fg{\leavevmode\raise.3ex\hbox{~$\!\scriptscriptstyle\,\rangle\!\rangle$}}

\def\transp #1{\vphantom{#1}^{\mathrm t}\! {#1}}

\addtocounter{section}{0}             
\numberwithin{equation}{section}       

\begin{document}
\selectlanguage{french}

\title[It\'eration d'applications rationnelles dans les espaces de matrices]{It\'eration d'applications rationnelles \\dans les espaces de matrices}

\author{Dominique \textsc{Cerveau}}

\address{Membre de l'Institut Universitaire de France.
IRMAR, UMR 6625 du CNRS, Universit\'e de Rennes $1$, $35042$ Rennes, France.}
\email{dominique.cerveau@univ-rennes1.fr}

\author{Julie \textsc{D\'eserti}}

\address{Institut de Math\'ematiques de Jussieu, Universit\'e Paris $7$, Projet G\'eom\'etrie et Dynamique, Site Chevaleret, Case $7012$, $75205$ Paris Cedex 13, France.}
\email{deserti@math.jussieu.fr}

\maketitle{}

\selectlanguage{english}
\begin{abstract}
The iteration of rational maps is well-understood in dimension $1$ but less so in higher dimensions. We study some maps on spaces of matrices which present a weak complexity with respect to the ring structure. First we give some properties of certain rational maps; the simplest example is the rational map which sends the matrix $\mathrm{M}$ onto $\mathrm{M}^2$ for which we exhibit some dynamical properties. Finally we deal with some small perturbations of this map.

\noindent\emph{2010 Mathematics Subject Classification. --- $14$E$05$, $32$H$50$, $37$B$05$.}
\end{abstract}

\selectlanguage{french}

\section*{Introduction}

\noindent L'it\'eration des applications rationnelles est tr\`es bien comprise en dimension un, un peu moins en dimension deux et encore moins en dimension plus grande. Nous nous proposons d'\'etudier ici des applications sp\'eciales sur les espaces de matrices qui pr\'esentent une \og faible complexit\'e\fg\, par rapport \`a la structure d'anneau. Par souci de simplicit\'e nous travaillons sur les matrices $2\times 2$, bien que la plupart du discours se laisse g\'en\'eraliser sans probl\`eme.

\smallskip

\noindent Dans un premier temps nous nous int\'eressons aux transformations compatibles \`a la conjugaison, \emph{i.e.} aux applications rationnelles $\Phi\colon\mathcal{M}(2;\mathbb{C}) \dashrightarrow\mathcal{M}(2;\mathbb{C})$ telles que $\mathrm{A}\Phi(\mathrm{M})\mathrm{A}^{-1}= \Phi(\mathrm{A}\mathrm{M}\mathrm{A}^{-1})$ pour tout $\mathrm{A}$ dans $\mathrm{GL}(2;\mathbb{C})$ et tout $\mathrm{M}$ l\`a o\`u cela a un sens. Un exemple de ce type d'applications est donn\'e par les polyn\^omes de matrices. Nous commen\c{c}ons par pr\'esenter des propri\'et\'es satisfaites par ces transformations comme par exemple: l'invariance du groupe diagonal $\mathcal{D}$ ou encore le fait qu'une telle transformation est birationnelle si et seulement si sa restriction \`a~$\mathcal{D}$ l'est. Au centre $\mathcal{C}=\big\{\lambda \mathrm{Id}\,\big\vert\lambda\in\mathbb{C}\big\}$ on peut associer la fibration $\mathcal{P}$ en $2$-plans d\'efinie comme suit: si $\mathrm{M}$ d\'esigne un \'el\'ement de~$\mathcal{M}(2;\mathbb{C}) \setminus\mathcal{C},$ on d\'efinit $\mathcal{P}(\mathrm{M})$ comme l'unique plan contenant~$\mathrm{M}$ et $\mathcal{C}$: le $2$-plan $\mathcal{P}(\mathrm{M})$ n'est rien d'autre que l'ensemble des matrices qui commutent \`a~$\mathrm{M}$. Cette fibration est invariante fibre \`a fibre par toute application rationnelle compatible \`a la conjugaison. Bien s\^ur on peut consid\'erer l'application $\mathrm{Inv}$ de $\mathcal{M}(2;\mathbb{C})$ dans~$\mathbb{C}^2$ qui \`a une matrice $\mathrm{M}$ associe ses invariants de similitude $(\mathrm{tr}\,\mathrm{M},\det\mathrm{M})$. Par d\'efinition une application rationnelle~$\Phi$ compatible \`a la conjugaison laisse invariant le feuilletage associ\'e \`a la fibration $\mathrm{Inv}$; plus pr\'ecis\'ement il existe une application rationnelle $\mathrm{Sq}\,\Phi\colon\mathbb{C}^2\to\mathbb{C}^2$ telle que $\mathrm{Inv}\,\circ \Phi= \mathrm{Sq}\,\Phi\circ\mathrm{Inv}$. R\'eciproquement on peut se demander \`a quelle condition une application rationnelle de $\mathbb{C}^2$ se rel\`eve \`a $\mathcal{M}(2; \mathbb{C})$, question \`a laquelle nous r\'epondrons. \'Evidemment la transforma\-tion~$\mathrm{Sq}\,\Phi$ contient une grande partie de la dynamique de l'application initiale $\Phi$.

\noindent Nous nous int\'eressons ensuite tout particuli\`erement \`a l'application $\Phi_\mathrm{Id}\colon\mathcal{M}(2;\mathbb{C})\to\mathcal{M}(2;\mathbb{C})$, $\mathrm{M} \mapsto\mathrm{M}^2$, exemple typique d'application compatible \`a la conjugaison. Apr\`es avoir donn\'e quelques propri\'et\'es satisfaites par cette application, nous d\'ecrivons l'adh\'erence de ses points p\'eriodiques ainsi que le bord du bassin d'attraction de la matrice nulle qui donne naissance \`a une int\'eressante hypersurface \textsc{Levi}-plate. Une fa\c{c}on de mesurer la complexit\'e d'une transformation est d'examiner son centralisateur. C'est dans cette optique que nous d\'eterminons le groupe $\mathrm{Aut}(\mathcal{M}(2;\mathbb{C});\Phi_\mathrm{Id})$ des automorphismes holomorphes de $\mathcal{M}(2;\mathbb{C})$ qui commutent \`a $\Phi_\mathrm{Id}$ ainsi que $\mathrm{Bir}(\mathcal{M}(2;\mathbb{C});\Phi_\mathrm{Id})$ le groupe des transformations birationnelles de $\mathcal{M}(2;\mathbb{C})$ qui commutent \`a~$\Phi_\mathrm{Id}$; le premier est engendr\'e par les applications de conjugaison $\sigma_{\mathrm{P}}\colon\mathrm{M}\mapsto \mathrm{P} \mathrm{M}\mathrm{P}^{-1}$ o\`u $\mathrm{P}$ appartient \`a $\mathrm{GL}(2;\mathbb{C})$, et la transposition $\mathrm{M}\mapsto\transp\, \mathrm{M}$. Pour obtenir le second il faut ajouter l'application \og inverse\fg: $\mathrm{M}\mapsto\mathrm{M}^{-1}$.

\begin{proposition}
Le groupe $\mathrm{Aut}(\mathcal{M}(2;\mathbb{C});\Phi_\mathrm{Id})$ est isomorphe \`a $\mathrm{PGL}(2;\mathbb{C})\rtimes\mathbb{Z} /2\mathbb{Z}$.

\medskip

\noindent Le groupe $\mathrm{Bir}(\mathcal{M}(2;\mathbb{C});\Phi_\mathrm{Id})$ est engendr\'e par $\mathrm{Aut} (\mathcal{M}(2;\mathbb{C});\Phi_\mathrm{Id})$ et par l'involution $\iota\colon\mathrm{M}\mapsto\mathrm{M}^{-1}$.
\end{proposition}

\noindent Enfin nous consid\'erons des d\'eformations sp\'eciales de $\Phi_\mathrm{Id}$, $\Phi_\mathrm{A}\colon\mathcal{M}(2;\mathbb{C})\to\mathcal{M}(2;\mathbb{C})$, $\mathrm{M}\mapsto\mathrm{A}\mathrm{M}^2$ avec $\mathrm{A}$ dans~$\mathrm{GL}(2;\mathbb{C})$. Ce sont les applications \og monomiales\fg\, les plus simples et elles sont en g\'en\'eral non compatibles \`a la conjugaison. Alors qu'\`a une variable les transformations $z\mapsto z^2$ et $z\mapsto az^2$ sont lin\'eairement conjugu\'ees, la situation ici est plus complexe. Nous pr\'ecisons tout du moins pour $\mathrm{A}$ g\'en\'erique, {\it i.e.} pour $\mathrm{A}$ de la forme $\left[\begin{array}{cc}\lambda & 0\\0&\frac{1}{\lambda}\end{array}\right]$, l'ensemble des orbites p\'eriodiques et l'adh\'erence de cet ensemble:

\begin{proposition}\label{add}
Pour $\lambda$ g\'en\'erique l'adh\'erence des points p\'eriodiques de $\Phi_\mathrm{A}$ est constitu\'ee  
\begin{itemize}
\item d'un tore $\mathbb{S}^1\times\mathbb{S}^1$ contenu dans $y=z=0;$

\item de deux $\mathbb{S}^1\times\mathbb{C},$ pr\'ecis\'ement $\frac{1}{\lambda}\mathbb{S}^1\times\mathbb{C}\times\big\{0\big\}\times\big\{0\big\}$ et $\big\{0\big\}\times\big\{0\big\}\times\mathbb{C}\times\lambda\mathbb{S}^1$ $($en identifiant $\mathcal{M}(2;\mathbb{C})=\left\{\left[\begin{array}{cc}x & y \\ z & t\end{array}\right]\right\}$ \`a $\mathbb{C}^4=\big\{(x,y,z,t)\,\vert\, x,\,y,\,z,\, t\in\mathbb{C}\big\})$; 

\item de la matrice nulle.
\end{itemize}
\end{proposition}

\noindent Alors que pour $z\mapsto z^2$ l'adh\'erence des points p\'eriodiques (diff\'erent de $0$) est le bord du bassin d'attraction de l'origine, nous d\'eduisons de la Proposition \ref{add} que ce n'est pas le cas ici. Nous donnons aussi des exemples d'orbites born\'ees, non p\'eriodiques et non contenues dans le bassin d'attraction de l'origine. Le cas o\`u $\mathrm{A}$ est une matrice quaternionique est particuli\`erement riche. Nous effectuons quelques exp\'eriences num\'eriques afin de mieux comprendre le bassin d'attraction de la matrice nulle, son bord ainsi que la dynamique de ce type de transformations. En particulier nous nous int\'eressons \`a la restriction de $\Phi_\mathrm{A}$ \`a l'ensemble $\Bbb H_1$ des quaternions de module $1$ qui est d\'ecrite par $f_\vartheta\colon(x,y)\mapsto\big(\mathrm{e}^{\mathrm{i}\vartheta}(x^2+\vert x\vert^2-1),\mathrm{e}^{\mathrm{i}\vartheta}y(x+\overline{x})\big)$; remarquons qu'elle laisse la famille de cercles param\'etr\'ee par $\eta\mapsto(x,y\mathrm{e}^{\mathrm{i}\eta})$ globalement invariante. La premi\`ere composante $f_\vartheta^1$ de $f_\vartheta$ indique comment passer d'un cercle \`a l'autre; elle s'identifie \`a une application de $\mathbb{R}^2$ dans lui-m\^eme ce qui nous permet de d\'eterminer quelques-unes de ses propri\'et\'es. Les applications $f_\vartheta^1$ pr\'eservent le disque unit\'e de~$\mathbb{R}^2$ et c'est la dynamique dans ce disque qui nous int\'eresse. Nous avons par exemple r\'ealis\'e des exp\'eriences num\'eriques visant \`a mesurer la fa\c{c}on dont les orbites s'approchent du bord du disque unit\'e et \`a quelle vitesse:

\medskip
\begin{center}
\begin{tabular}{ccc}
\begin{picture}(0,0)%
\includegraphics{theta0b.pstex}%
\end{picture}%
\setlength{\unitlength}{3947sp}%
\begingroup\makeatletter\ifx\SetFigFont\undefined%
\gdef\SetFigFont#1#2#3#4#5{%
  \reset@font\fontsize{#1}{#2pt}%
  \fontfamily{#3}\fontseries{#4}\fontshape{#5}%
  \selectfont}%
\fi\endgroup%
\begin{picture}(1200,1200)(1201,-5161)
\end{picture}%
 \hspace*{10mm}& \hspace*{10mm}\begin{picture}(0,0)%
\includegraphics{theta1b.pstex}%
\end{picture}%
\setlength{\unitlength}{3947sp}%
\begingroup\makeatletter\ifx\SetFigFont\undefined%
\gdef\SetFigFont#1#2#3#4#5{%
  \reset@font\fontsize{#1}{#2pt}%
  \fontfamily{#3}\fontseries{#4}\fontshape{#5}%
  \selectfont}%
\fi\endgroup%
\begin{picture}(1200,1200)(1201,-5161)
\end{picture}%
\hspace*{10mm}& \hspace*{10mm}\begin{picture}(0,0)%
\includegraphics{thetapi2b.pstex}%
\end{picture}%
\setlength{\unitlength}{3947sp}%
\begingroup\makeatletter\ifx\SetFigFont\undefined%
\gdef\SetFigFont#1#2#3#4#5{%
  \reset@font\fontsize{#1}{#2pt}%
  \fontfamily{#3}\fontseries{#4}\fontshape{#5}%
  \selectfont}%
\fi\endgroup%
\begin{picture}(1200,1200)(1201,-3961)
\end{picture}%
 \\
\hspace*{-12mm}$\vartheta=0$ &$\vartheta=1$&\hspace*{10mm}$\vartheta=\frac{\pi}{2}$
\end{tabular}
\end{center}
\medskip

\noindent Comme dans le cas sp\'ecial o\`u $\mathrm{A}=\mathrm{Id}$ nous d\'eterminons le sous-groupe des transformations biholomorphes de $\mathcal{M}(2;\mathbb{C})$ qui commutent \`a $\Phi_\mathrm{A}$. 

\begin{proposition}
Si $\mathrm{A}$ est un multiple de l'identit\'e alors $\Phi_\mathrm{A}$ est conjugu\'e \`a $\Phi_\mathrm{Id}$, sinon $\mathrm{Aut}(\mathcal{M}(2;\mathbb{C});\Phi_\mathrm{A})$ est form\'e des~$\sigma_\mathrm{P}$ o\`u  $\sigma_\mathrm{P}$ d\'esigne l'application $\mathrm{M}\mapsto \mathrm{P}\mathrm{M}\mathrm{P}^{-1}$ et $\mathrm{P}$ une matrice qui commute \`a $\mathrm{A}$. En fait $\mathrm{Aut}(\mathcal{M}(2;\mathbb{C});\Phi_\mathrm{A})$ s'identifie \`a $\mathbb{C}^*$ agissant sur $\mathcal{M}(2;\mathbb{C})$ de la fa\c{c}on suivante: $(x,y,z,t,\alpha)\mapsto\left(x,\alpha y,\frac{z}{\alpha},t\right)$. Les orbites de cette action sont aussi celles du champ de vecteurs invariant $y\frac{\partial}{\partial y}-z\frac{\partial}{\partial z}$.
\end{proposition}

\noindent Remarquons que le bord du bassin d'attraction de l'origine suivant $\Phi_\mathrm{A}$ est invariant par l'action de ce groupe. Beaucoup de questions concernant ces transformations $\Phi_\mathrm{A}$ restent ouvertes; nous listons celles qui nous paraissent  les plus pertinentes, tant sur le plan th\'eorique que d'un point de vue num\'erique.

\subsection*{Remerciements} Nous remercions M. \textsc{Baker}, D. \textsc{Boucher}, S. \textsc{Cantat}, G. \textsc{Casale}, S. \textsc{Gou\"{e}zel} et F. \textsc{Loray} pour leur disponibilit\'e. 

\section{G\'en\'eralit\'es}

\noindent Rappelons quelques d\'efinitions. 

\begin{defis}
Une application m\'eromorphe $f\colon X\dashrightarrow Z$ entre deux vari\'et\'es complexes compactes est d\'efinie par son graphe $\Gamma(f)\subset X\times Z$; ce graphe est une sous-vari\'et\'e irr\'eductible pour laquelle la projection $\pi_1\colon\Gamma(f)\to~X$ sur le premier facteur est une application holomorphe surjective propre dont la fibre g\'en\'erique est un point (\cite{Fi}). Le \emph{lieu d'ind\'etermination} de~$f$ est l'ensemble des points o\`u $\pi_1$ n'admet pas d'inverse local, on le note $\mathrm{Ind}\,f$. L'application~$f$ est \emph{dominante} si la seconde projection $\pi_2\colon\Gamma(f)\to Z$ est surjective. Notons $\mathrm{Exc}\,\pi_2$ l'ensemble des points o\`u~$\pi_2$ n'est pas une application finie; on d\'efinit l'\emph{ensemble exceptionnel} de~$f$ par~$\mathrm{Exc}\,f=\pi_1(\mathrm{Exc}\,\pi_2)$.
 
\noindent Dans le cas particulier o\`u $f\colon\mathbb{C}^{n+1}\to\mathbb{C}^{n+1}$ est une application polynomiale homog\`ene ad-hoc repr\'esentant la transformation rationnelle $$\mathbb{P}(f) \colon\mathbb{P}^n(\mathbb{C})\dashrightarrow\mathbb{P}^n(\mathbb{C}),$$ un point $[m]=(m_0: \ldots :m_n)$ est un point d'ind\'etermination de $\mathbb{P}(f)$ si $f(m)=0$. Une sous-vari\'et\'e irr\'educti\-ble~$V\subset\mathbb{P}^n(\mathbb{C})$ est \emph{contract\'ee} par $\mathbb{P}(f)$ si la dimension de $\overline{\mathbb{P}(f)(V\setminus\mathrm{Ind}\,\mathbb{P}(f))}$ est strictement inf\'erieure \`a celle de $V$. Une telle sous-vari\'et\'e est contenue dans l'ensemble $\mathrm{Exc}\,\mathbb{P}(f)$.

\noindent Une application homog\`ene $f\colon\mathbb{C}^{n+1}\to\mathbb{C}^{n+1}$ respecte la fibration de \textsc{Hopf} par les droites passant par l'origine: une droite $d$ est envoy\'ee par $f$ sur une autre droite $f(d)$ \`a moins que $f$ s'annule sur $d$. L'application~$\mathbb{P}(f)$ explique comment sont \'echang\'ees ces droites par $f$.

\noindent Nous travaillerons aussi avec des applications rationnelles $f\colon\mathbb{C}^n\to\mathbb{C}^n$ qui sont simplement les applications dont les composantes $f_i$ sont rationnelles, \emph{i.e.} des quotients de polyn\^omes. Pour ces applications on a aussi la notion de sous-vari\'et\'es contract\'ees et d'ensemble d'ind\'etermination qui co\"{i}ncide avec l'union des p\^oles des composantes.
\end{defis}

\noindent Soit $\mathcal{M}(2;\mathbb{C})=\left\{
\left[\begin{array}{cc}a&b\\c&d\end{array} \right]\,\Big\vert\,a,\,b,\,c,\,d\in\mathbb{C}\right\}\simeq \mathbb{C}^4$ l'ensemble des matrices $2\times 2$ \`a coefficients complexes. Notons $\mathrm{Inv}$ l'application de $\mathcal{M}(2;\mathbb{C})$ \`a valeurs dans $\mathbb{C}^2$ qui \`a une matrice associe ses invariants de similitude
\begin{align*}
\mathrm{Inv}\colon\mathcal{M}(2;\mathbb{C})\to\mathbb{C}^2,&& \mathrm{M}\mapsto( \mathrm{tr}\,\mathrm{M},\det\,\mathrm{M}) 
\end{align*}
\noindent avec les notations habituelles. La sous-alg\`ebre de $\mathcal{M}(2;\mathbb{C})$ des matrices diagonales sera not\'ee $\mathcal{D}$ $$\mathcal{D}=\left\{\left[\begin{array}{cc}\lambda_1&0\\0&\lambda_2\end{array} \right]\,\Big\vert\,\lambda_i\in\mathbb{C}\right\}.$$ D\'esignons par $\mathrm{diag}(\lambda_1,\lambda_2)$ l'\'el\'ement $\left[\begin{array}{cc}\lambda_1 & 0\\0&\lambda_2\end{array}\right]$ de $\mathcal{D}$ et par ${\bf 0}$ la matrice nulle $\mathrm{diag}(0,0)$.

\noindent Soit $\Phi\colon\mathcal{M}(2;\mathbb{C})\dashrightarrow\mathcal{M}(2;\mathbb{C})$ une application rationnelle dominante. Si l'on identifie $\mathcal{M}(2;\mathbb{C})$ \`a $\mathbb{C}^4$ que l'on voit comme carte affine de $\mathbb{P}^4(\mathbb{C})$, alors $\Phi$ induit une application rationnelle not\'ee $\widetilde{\Phi}$ de $\mathbb{P}^4(\mathbb{C})$ dans lui-m\^eme. Dans le cas sp\'ecial o\`u $\Phi$ est homog\`ene ($\Phi(s\mathrm{M})=s^d\Phi(\mathrm{M})$), $\Phi$ induit une application rationnelle not\'ee $\mathbb{P}(\Phi)\colon\mathbb{P}^3(\mathbb{C}) \dashrightarrow\mathbb{P}^3( \mathbb{C})$. 

\begin{defi}
Une application $\Phi\colon\mathcal{M}(2;\mathbb{C})\to \mathcal{M}(2;\mathbb{C})$ sera dite \emph{monomiale} (relativement \`a la multiplication) si elle est du type $$\Phi(\mathrm{M})=\mathrm{A}_1\mathrm{M}\mathrm{A}_2\mathrm{M}\ldots \mathrm{A}_p\mathrm{M}\mathrm{A}_{p+1}$$ o\`u les $\mathrm{A}_i$ sont des \'el\'ements de $\mathrm{GL}(2;\mathbb{C})$ fix\'es.
\end{defi}

\begin{eg}
L'application $\Phi_\mathrm{Id}\colon\mathcal{M}(2;\mathbb{C})\to\mathcal{M}(2;\mathbb{C})$ d\'efinie par $\Phi_\mathrm{Id}(\mathrm{M})=\mathrm{M}^2$ est monomiale. Un calcul \'el\'ementaire montre que si $\mathrm{M}=\left[\begin{array}{cc} x & y\\ z & t \end{array}\right]$, alors $\Phi_\mathrm{Id}(\mathrm{M})=\left[\begin{array}{cc} x^2+yz& y(x+t)\\ z(x+t)& t^2+yz\end{array}\right]$. L'application $\Phi_\mathrm{Id}$ est g\'en\'eriquement finie au sens suivant: pour $\mathrm{M}$ g\'en\'erique, $\#\,\Phi_\mathrm{Id}^{-1} (\mathrm{M})=~4$.  Si $\mathrm{M}$ est une matrice nilpotente, \emph{i.e.} $\mathrm{Inv}(\mathrm{M})=(0,0)$, alors $\Phi_\mathrm{Id}(\mathrm{M})={\bf 0}$ de sorte que $[\mathrm{M}]$ est un point d'ind\'etermination de $\mathbb{P}(\Phi_\mathrm{Id})$. En fait l'ensemble $\mathrm{Ind}( \mathbb{P}(\Phi_\mathrm{Id}))$ des points d'ind\'etermination de $\mathbb{P}(\Phi_\mathrm{Id})$ est exactement le projectivis\'e de l'ensemble $\mathcal{N}(2;\mathbb{C})=\big\{\mathrm{M} \in \mathcal{M}(2;\mathbb{C})\,\big\vert\,\mathrm{Inv}(\mathrm{M})=(0,0)\big\}$ des matrices nilpotentes; c'est donc une conique plane lisse. Remarquons que l'image de $\Phi_\mathrm{Id}$ est pr\'ecis\'ement $$W=\big\{\mathbf{0}\big\}\cup\Big(\mathcal{M}(2;\mathbb{C})\setminus\mathcal{N}(2;\mathbb{C})\Big).$$ 

\noindent La sous-alg\`ebre $\mathfrak{sl}(2;\mathbb{C})= \big\{\mathrm{M}\in\mathcal{M}(2;\mathbb{C})\,\big\vert\,\mathrm{tr}\,\mathrm{M}=0\big\}$ est contract\'ee par $\Phi_\mathrm{Id}$ sur le centre $\mathcal{C}=\big\{\lambda\mathrm{Id}\,\big\vert\,\lambda\in~\mathbb{C}\big\}$. 
\end{eg}

\noindent Remarquons que les matrices triangulaires sup\'erieures forment un sous-espace invariant par $\Phi_\mathrm{Id}$, de m\^eme que les matrices triangulaires inf\'erieures, et bien s\^ur $\mathcal{D}$, qui est leur intersection. Leurs projectivis\'es sont invariants par $\mathbb{P}(\Phi_\mathrm{Id})$. De la m\^eme fa\c{c}on la quadrique $\mathcal{Q}=\big\{\mathrm{M}\in\mathcal{M}(2;\mathbb{C})\,\big\vert\,\det\, \mathrm{M}=0\big\}$ est compl\`etement invariante par $\Phi_\mathrm{Id}$. On peut ais\'ement v\'erifier que le d\'eterminant jacobien de $\Phi_\mathrm{Id}$ est $$\det\mathrm{jac}\,\Phi_{\mathrm{Id}_{(\mathrm{M})}}= 4(\mathrm{tr}\,\mathrm{M})^2 \det\,\mathrm{M};$$ de plus, lorsque $\det\,\mathrm{M}=0$ nous avons $\Phi_\mathrm{Id}(\mathrm{M})=(\mathrm{tr}\, \mathrm{M}) \mathrm{M}$. En r\'esulte que le projectivis\'e des matrices de trace nulle $\mathbb{P}(\big\{\mathrm{M}\in\mathcal{M}(2;\mathbb{C})\,\big\vert\,\mathrm{tr}\,\mathrm{M}=0\big\})$ est la seule surface de $\mathbb{P}^3(\mathbb{C})$ contract\'ee par $\mathbb{P}(\Phi_\mathrm{Id})$. Nous reviendrons plus tard (\S\ref{Phi2}) sur la dynamique de l'application $\Phi_\mathrm{Id}$ qui se r\'eduit peu ou prou \`a celle de $z\mapsto z^2$ dans le plan complexe (ne serait-ce que via la formule $\det \mathrm{M}^2=(\det\mathrm{M})^2$ qui induit une \og semi-conjugaison\fg) mais produit des objets int\'eressants.

\begin{defi}
Soit $\Phi\colon\mathcal{M}(2;\mathbb{C})\dashrightarrow\mathcal{M}(2;\mathbb{C})$ une application rationnelle. On dit que $\Phi$ est \emph{compatible \`a la conjugaison} si $$\mathrm{A}\Phi(\mathrm{M})\mathrm{A}^{-1}=\Phi(\mathrm{A}\mathrm{M}\mathrm{A}^{-1})$$ chaque fois que $\mathrm{M}$ et $\mathrm{A}\mathrm{M}\mathrm{A}^{-1}$ sont dans $\mathcal{M}(2;\mathbb{C})\setminus\mathrm{Ind}\,\Phi$.

\noindent Notons que si $\Phi$ est compatible \`a la conjugaison, alors $\mathrm{Ind}\,\Phi$ est invariant sous l'action adjointe $\mathrm{A}\mathrm{M}\mathrm{A}^{-1}$.
\end{defi}

\noindent Ainsi l'application $\Phi_\mathrm{Id}$ d\'efinie pr\'ec\'edemment est compatible \`a la conjugaison; il en est de m\^eme des polyn\^omes d'endomorphismes. Plus g\'en\'eralement si $r(s)=P(s)/Q(s)$ est une fonction rationnelle en une variable~$s$, alors $r(\mathrm{M})=P(\mathrm{M})Q(\mathrm{M})^{-1}$ d\'efinit une application rationnelle compatible \`a la conjugaison. Une application compatible \`a la conjugaison induit une application de l'espace quotient $\mathcal{M}(2;\mathbb{C})/\text{conjugaison}$. Tout ceci se g\'en\'eralise \'evidemment en dimension quelconque.

\begin{rem}
Une application monomiale compatible \`a la conjugaison est du type $\Phi_{\alpha,d}(\mathrm{M})= \alpha\mathrm{M}^d$, $\alpha$ dans $\mathbb{C}^*$. L'alg\`ebre $\mathfrak{sl}(2;\mathbb{C})$ n'est pas toujours contract\'ee par $\Phi_{\alpha,d};$ par exemple pour $\mathrm{M}=\left[\begin{array}{cc}x& y\\z&-x\end{array}\right]$ et $d=3$ nous avons $$\mathrm{M}^3=(x^2+yz)\left[\begin{array}{cc}x& y\\z&-x\end{array}\right].$$ Dans ce cas les matrices de trace nulle sont fixes pour l'application $\mathbb{P}(\Phi_{\alpha,3})$ induite sur $\mathbb{P}^3(\mathbb{C})=\mathbb{P}(\mathcal{M}(2;\mathbb{C}))$. On constate que l'ensemble d'ind\'etermination de $\mathbb{P}(\Phi_{\alpha,d})$ est encore le projectivis\'e de $\mathcal{N}(2;\mathbb{C})$. Plus g\'en\'eralement nous avons en restriction \`a $\mathfrak{sl}(2;\mathbb{C})$
\begin{align*}
& \mathrm{M}^{2d}=(x^2+yz)^d\mathrm{Id}, && \mathrm{M}^{2d+1}=(x^2+yz)^d\mathrm{M}.
\end{align*}
\end{rem}

\begin{rem}\label{red}
Soit $\Phi$ l'application monomiale d\'efinie par $\Phi(\mathrm{M})=\mathrm{A}_1\mathrm{M}\ldots\mathrm{A}_n\mathrm{M}\mathrm{A}_{n+1}$, $n\geq 2$. Quitte \`a conjuguer~$\Phi$ par un automorphisme de $\mathcal{M}(2;\mathbb{C})$ du type $\mathrm{M}\mapsto \mathrm{P}\mathrm{M} \mathrm{Q}$ on peut supposer que $\mathrm{A}_n=\mathrm{Id}$. On en d\'eduit apr\`es cette conjugaison que l'ensemble $\mathcal{N}(2;\mathbb{C})$ des matrices nilpotentes est aussi contract\'e sur ${\bf 0}$.  

\noindent \`A conjugaison pr\`es par $\mathrm{M}\mapsto\rho\mathrm{M}$ on peut se ramener \`a $\det\Phi(\mathrm{M})=(\det\mathrm{M})^n$. En particulier les hypersurfaces quadratiques $\big\{\mathrm{M}\in\mathcal{M}(2;\mathbb{C})\,\big\vert\,\det\mathrm{M}=\varepsilon\big\}$ sont invariantes par l'it\'er\'e $\Phi^p$ pour chaque racine $(n^p-1)-$i\`eme de l'unit\'e. \'Evidemment les quadriques $\big\{\mathrm{M}\in\mathcal{M}(2;\mathbb{C})\,\big\vert\,\det\mathrm{M}=0\big\}$ et $\mathrm{SL}(2;\mathbb{C})=\big\{\mathrm{M}\in\mathcal{M}(2;\mathbb{C})\,\big\vert\,\det\mathrm{M}=1\big\}$ sont invariantes par $\Phi$.
\end{rem}

\section{Applications compatibles \`a la conjugaison, propri\'et\'es et exemples}

\subsection{Premi\`eres propri\'et\'es}

\noindent Il est \`a peu pr\`es clair qu'une application compatible \`a la conjugaison est d\'etermin\'ee par sa restriction aux matrices diagonales. Pr\'ecisons cel\`a.

\begin{lem}\label{invdiag}
Soit $\Phi\colon\mathcal{M}(2;\mathbb{C})\dashrightarrow\mathcal{M}(2;\mathbb{C})$ une application rationnelle compatible \`a la conjugaison. Alors~$\mathcal{D}$ est invariant par $\Phi$, \emph{i.e.} $\Phi(\mathcal{D}\setminus\mathrm{Ind}\,\Phi)\subset \mathcal{D}$.
\end{lem}

\begin{proof}[D\'emonstration]
Remarquons que, par densit\'e des matrices diagonalisables, $\mathcal{D}$ n'est pas contenu dans~$\mathrm{Ind}\,\Phi$. Soit $\mathrm{M}$ un \'el\'ement de $\mathcal{D}\setminus\mathrm{Ind}\,\Phi;$ en particulier $\mathrm{M}$ commute \`a $\mathrm{diag}(1,2)$ d'o\`u $$\Phi(\mathrm{M})=\Phi(\mathrm{diag}(1,2)\,\mathrm{M}\,\mathrm{diag}(1, 1/2)).$$ Puisque $\Phi$ est compatible avec la conjugaison nous avons $\Phi(\mathrm{M})=\mathrm{diag}(1,2)\,\Phi(\mathrm{M})\,\mathrm{diag}(1,1/2);$ cette \'egalit\'e assure que $\Phi(\mathrm{M})$ appartient \`a $\mathcal{D}$. 
\end{proof}

\begin{rem} 
Si $\Phi$ est dominante, nous avons l'\'egalit\'e $\overline{\Phi(\mathcal{D}\setminus\mathrm{Ind}\, \Phi)}=\mathcal{D}$ o\`u l'adh\'erence est prise au sens ordinaire.
\end{rem}

\noindent Si $\Phi$ est compatible \`a la conjugaison nous avons $\Phi(\mathrm{diag}(\lambda_1,\lambda_2))=\mathrm{diag}(\varphi_1(\lambda_1,\lambda_2),\varphi_2(\lambda_1,\lambda_2))$ o\`u les $\varphi_i$ sont rationnels (Lemme \ref{invdiag}). Mais comme les matrices $\mathrm{diag}(\lambda_1,\lambda_2)$ et $\mathrm{diag}(\lambda_2,\lambda_1)$ sont conjugu\'ees par $\left[\begin{array}{cc}0 & 1\\1 & 0\end{array}\right]$ nous avons $$\Phi(\mathrm{diag}(\lambda_2,\lambda_1))=\mathrm{diag}(\varphi_1(\lambda_2,\lambda_1), \varphi_2 ( \lambda_2,\lambda_1))=\mathrm{diag}(\varphi_2(\lambda_1,\lambda_2),\varphi_1(\lambda_1,\lambda_2)).$$ Par suite $$\Phi(\mathrm{diag}(\lambda_1,\lambda_2))=\mathrm{diag}(\varphi_1(\lambda_1,\lambda_2),\varphi_1(\lambda_2,\lambda_1)).$$

\noindent Inversement soit $\Psi\colon\mathbb{C}^2\dashrightarrow\mathbb{C}$ une fonction rationnelle; si $\mathrm{M}$ est une matrice $2\times 2$ g\'en\'erique, $\mathrm{M}$ s'\'ecrit $\mathrm{P}\,\mathrm{diag}(\lambda_1,\lambda_2)\,\mathrm{P}^{-1}$ et on d\'efinit $$\Phi(\mathrm{M})=\mathrm{P}\,\mathrm{diag}(\Psi(\lambda_1,\lambda_2),\Psi(\lambda_2,\lambda_1))\,\mathrm{P}^{-1}.$$

\noindent Cette d\'efinition ne d\'epend ni du choix de $\mathrm{P}$, ni de l'ordre choisi pour \'enum\'erer les valeurs propres. Par suite~$\Phi$ s'\'etend en une application rationnelle de $\mathcal{M}(2;\mathbb{C})$ dans $\mathcal{M}(2;\mathbb{C})$.

\noindent Les propri\'et\'es et la dynamique de $\Phi$ sont essentiellement cod\'ees par celles de l'application correspondante $(\Psi(\lambda_1,\lambda_2),\Psi(\lambda_2, \lambda_1))$. En particulier nous avons la:

\begin{pro}\label{birationnelle}
Soit $\Phi\colon\mathcal{M}(2;\mathbb{C})\dashrightarrow\mathcal{M}(2;\mathbb{C})$ une application rationnelle compatible \`a la conjugaison;~$\Phi$ est birationnelle si et seulement si sa restriction $\Psi=\Phi_{\vert\mathcal{D}}\colon\mathcal{D}\dashrightarrow\mathcal{D}$ l'est.
\end{pro}

\begin{proof}[D\'emonstration]
Supposons que $\Psi=\Phi_{\vert\mathcal{D}}$ soit birationnelle, ou ce qui revient au m\^eme g\'en\'eriquement injective.
Soient $\mathrm{A}$, $\mathrm{B}$ dans $\mathcal{M}(2;\mathbb{C})$ diagonalisables, \emph{i.e.} $\mathrm{A}=\mathrm{P}\mathrm{diag}(\lambda_1,\lambda_2)\mathrm{P}^{-1}$ et $\mathrm{B}=\mathrm{Q}\mathrm{diag}(\mu_1,\mu_2)\mathrm{Q}^{-1}$. Si $\Phi(\mathrm{A})=\Phi(\mathrm{B})$, nous avons $$\mathrm{P}\,\mathrm{diag}( \Psi(\lambda_1, \lambda_2),\Psi(\lambda_2,\lambda_1))\, \mathrm{P}^{-1}=\mathrm{Q}\,\mathrm{diag}( \Psi(\mu_1, \mu_2),\Psi(\mu_2,\mu_1))\, \mathrm{Q}^{-1}.$$ Par suite $\Psi(\lambda_1,\lambda_2)$ vaut $\Psi(\mu_1,\mu_2)$ ou $\Psi(\mu_2,\mu_1)$. Quitte \`a modifier $\mathrm{Q}$ nous pouvons supposer que $\Psi(\lambda_1,\lambda_2)=\Psi(\mu_1,\mu_2)$ et~$\Psi(\lambda_2, \lambda_1)=\Psi(\mu_2,\mu_1)$. Puisque $\Psi=\Phi_{\vert\mathcal{D}}$ est injective nous avons $(\lambda_1,\lambda_2)=(\mu_1, \mu_2)$. Ainsi 
\begin{align*}
&\mathrm{A}=\mathrm{P}\,\mathrm{diag}(\lambda_1, \lambda_2)\,\mathrm{P}^{-1} &&\text{et} &&\mathrm{B}=\mathrm{Q}\,\mathrm{diag}(\mu_1, \mu_2)\,\mathrm{Q}^{-1}.
\end{align*}

\noindent L'\'egalit\'e $\Phi(\mathrm{A})=\Phi(\mathrm{B})$ implique $\mathrm{Q}^{-1}\mathrm{P}\,\mathrm{diag}( \Psi(\lambda_1, \lambda_2),\Psi(\lambda_2,\lambda_1))\,\mathrm{P}^{-1}\mathrm{Q}=\,\mathrm{diag}( \Psi(\lambda_1, \lambda_2),\Psi(\lambda_2,\lambda_1))\,;$ il en r\'esulte que $\mathrm{Q}^{-1}\mathrm{P}$ est une matrice diagonale $\mathrm{D}$. Par suite $$\mathrm{A}=\mathrm{P}\,\mathrm{diag}(\lambda_1, \lambda_2)\,\mathrm{P}^{-1}=\mathrm{Q}\mathrm{D}\,\mathrm{diag}(\lambda_1, \lambda_2)\,\mathrm{D}^{-1}\mathrm{Q}^{-1}=\mathrm{Q}\,\mathrm{diag}(\lambda_1, \lambda_2)\,\mathrm{Q}^{-1}=\mathrm{B}.$$ Ainsi si $\Psi$ est birationnelle, alors $\Phi$ l'est. 
\end{proof}

\noindent Dans le m\^eme ordre d'id\'ee on peut se demander si une transformation polynomiale (resp. un automorphisme polynomial) de $\mathcal{D}$ dans lui-m\^eme du type $(\Psi(\lambda_1,\lambda_2),\Psi(\lambda_2,\lambda_1))$ induit une transformation polynomiale (resp. un automorphisme polynomial) \'equivariante de $\mathcal{M}(2;\mathbb{C})$ dans lui-m\^eme. 

\smallskip

\noindent D\'esignons par $\tau$ l'involution de $\mathbb{C}^2$ d\'efinie par $\tau(\lambda_1, \lambda_2)=(\lambda_2,\lambda_1)$.

\begin{pro}\label{polynomiale}
Soit $\eta\colon\mathcal{D}\to\mathcal{D}$ une transformation polynomiale $\tau$-\'equivariante, \emph{i.e.} du type $$\eta(\lambda_1,\lambda_2)=(\Psi(\lambda_1,\lambda_2),\Psi(\lambda_2,\lambda_1)).$$ Alors $\eta$ s'\'etend en une application $\Phi\colon\mathcal{M}(2;\mathbb{C})\to\mathcal{M}(2; \mathbb{C})$ polynomiale et compatible \`a la conjugaison.
\end{pro}

\begin{proof}[D\'emonstration]
Soit $\Delta\subset\mathcal{M}(2;\mathbb{C})$ l'hypersurface discriminante: $\Delta=\big\{\mathrm{M}\in \mathcal{M}(2;\mathbb{C})\,\big\vert\,(\mathrm{tr}\,\mathrm{M})^2-4\det\mathrm{M}=~0\big\}$. Sa trace sur $\mathcal{D}$ est constitu\'ee des multiples de l'identit\'e: $\mathcal{D}\cap\Delta=\big\{ \lambda\mathrm{Id}\,\big\vert\,\lambda\in\mathbb{C}^*\big\}$. Comme nous l'avons vu on peut \'etendre $\eta$ en une application rationnelle $\Phi\colon\mathcal{M}(2;\mathbb{C})\dashrightarrow \mathcal{M}(2;\mathbb{C})$. Par construction $\Phi$ est holomorphe en restriction \`a~$\mathcal{M} (2;\mathbb{C})\setminus\Delta$. Soit $\mathrm{M}_0=\lambda_0\mathrm{Id}$ un point de $\mathcal{D}\cap\Delta;$ cette m\^eme construction assure que~$\Phi$ reste born\'ee sur un petit voisinage $\mathcal{V}(\mathrm{M}_0)$ de $\mathrm{M}_0$ priv\'e de $\Delta$. Il r\'esulte du th\'eor\`eme d'\textsc{Hartogs} que $\Phi$ s'\'etend holomorphiquement \`a $\mathcal{V}(\mathrm{M}_0)$. Ce raisonnement montre qu'en fait $\Phi$ est holomorphe sur un voisinage de~$\mathcal{D}\cap \Delta$. Mais si $\mathrm{M}$ appartient \`a $\Delta$, l'orbite adjointe de $\mathrm{M}$ coupe ce voisinage. Par cons\'equent $\Phi$ s'\'etend holomorphiquement \`a $\mathcal{M}(2;\mathbb{C})$. Comme $\Phi$ est rationnelle nous obtenons le r\'esultat annonc\'e.  
\end{proof}

\noindent Les Propositions \ref{birationnelle} et \ref{polynomiale} impliquent l'\'enonc\'e suivant.

\begin{cor}\label{corollaire}
Un automorphisme polynomial $\tau$-\'equivariant de $\mathcal{D}$ dans lui-m\^eme s'\'etend de fa\c{c}on unique en un automorphisme polynomial de $\mathcal{M}(2;\mathbb{C})$ dans $\mathcal{M}(2; \mathbb{C})$ compatible \`a la conjugaison.
\end{cor}

\begin{defis}
Soient $f\colon X\dashrightarrow X$ une application rationnelle dominante et $\mathcal{F}$ un feuilletage sur $X;$ on dit que $\mathcal{F}$ est \emph{invariant} par $f$ si $f^*\mathcal{F}=\mathcal{F}$. Plus pr\'ecis\'ement si $m$ est un point g\'en\'erique de $X$, alors $\mathcal{F}$ est r\'egulier en $m$, \emph{i.e.} donn\'e par les niveaux d'une submersion locale $g\colon\mathcal{V}(m)\to\mathbb{C}^k$ d'un voisinage de $m$ dans $\mathbb{C}^k;$ de m\^eme $m$ est g\'en\'erique pour~$f$, \emph{i.e.} $m$  est une valeur r\'eguli\`ere de $f$. Par suite $g\circ f$ est une submersion en chaque point de $f^{-1}(m)$ ce qui permet de d\'efinir les feuilles locales de $f^*\mathcal{F}$ comme les niveaux de $g\circ f$. Supposons que $\mathcal{F}$ soit d\'efini par une fibration, \emph{i.e.} la feuille g\'en\'erique de $\mathcal{F}$ est la fibre g\'en\'erique d'une application rationnelle $g\colon X\dashrightarrow Y;$ on  dit que $f$ \emph{pr\'eserve la fibration} $\mathcal{F}$ si $f^*\mathcal{F}=\mathcal{F}$. On dit que $f$ \emph{pr\'eserve la fibration~$\mathcal{F}$ fibre \`a fibre} si $g\circ f=g$.
\end{defis}

\begin{defis}
Un \emph{automorphisme \'el\'ementaire} de $\mathbb{C}^2$ est, \`a conjugaison pr\`es (dans le groupe des automorphismes), de la forme suivante
\begin{align*}
&(\alpha x+P(y),\beta y+\gamma), && \alpha,\,\beta\in\mathbb{C}^*,\,\gamma\in\mathbb{C},\, P \in\mathbb{C}[y].
\end{align*}
Une \emph{transformation de \textsc{H\'enon}} est par d\'efinition une transformation de la forme $(y,P(y)-\delta x)$ o\`u $\delta$ d\'esigne un \'el\'ement de $\mathbb{C}^*$ et $P$ un \'el\'ement de $\mathbb{C}[y]$ de degr\'e sup\'erieur ou \'egal \`a $2$.
\end{defis}

\begin{rem}
Les automorphismes \'el\'ementaires correspondent pr\'ecis\'ement aux automorphismes polynomiaux qui pr\'eservent une fibration rationnelle.
\end{rem}

\noindent \'Etant donn\'e un automorphisme polynomial $f$ de $\mathbb{C}^2$, nous avons l'alternative suivante (\cite{FM})
\begin{itemize}
\item $f$ est conjugu\'e \`a un automorphisme \'el\'ementaire;

\item $f$ est conjugu\'e \`a un produit de transformations de \textsc{H\'enon}; dans ce cas on dit que $f$ est de \emph{type \textsc{H\'enon} g\'en\'eralis\'e}.
\end{itemize}

\noindent Soit $\eta$ un automorphisme polynomial \'equivariant de $\mathcal{D}$ dans lui-m\^eme, \emph{i.e.} $\eta$ commute \`a la permuta\-tion~$\tau=(\lambda_2,\lambda_1)$. En particulier sa ligne de points fixes $y=x$ est pr\'eserv\'ee par $\eta$ qui du coup poss\`ede une courbe alg\'ebrique invariante. Il s'en suit que $\eta$ ne peut \^etre de type \textsc{H\'enon} g\'en\'eralis\'e (un tel automorphisme ne pr\'eserve pas de courbe alg\'ebrique, \cite{BS}). Il est n\'ecessairement \'el\'ementaire; nous en d\'eduisons que $\eta$ laisse invariante une fibration, d'o\`u le r\'esultat suivant.

\begin{thm}\label{fibration}
Soit $\Phi\colon\mathcal{M}(2;\mathbb{C})\to\mathcal{M}(2;\mathbb{C})$ un automorphisme polynomial compatible \`a la conjugaison. Alors~$\Phi_{\vert\mathcal{D}}$ est un automorphisme \'el\'ementaire; en particulier $\Phi$ pr\'eserve une fibration $\mathcal{L}\colon\mathcal{M}(2;\mathbb{C})\to\mathbb{C}$ polynomiale transverse \`a~$\mathcal{D}$ au sens o\`u $\mathcal{L}_{\vert\mathcal{D}}\colon\mathcal{D}\to\mathbb{C}$ est non constante.
\end{thm}

\begin{eg}\label{pasan}
L'automorphisme polynomial $(\lambda_1,\lambda_2)\mapsto (\lambda_1+(\lambda_2 -\lambda_1)^2, \lambda_2+(\lambda_1-\lambda_2)^2)$ de $\mathcal{D}$ dans $\mathcal{D}$ se rel\`eve en $$\left[ \begin{array}{cc} x & y \\ z & t \end{array}\right]\mapsto\left[\begin{array}{cc} x+(t-x)^2+4yz & y \\ z & t+(t-x)^2 +4yz \end{array}\right].$$ Ici la fibration invariante $\mathcal{L}$ est $(x-t)=$ cte.  Notons que l'hypersurface discriminante $\big\{\mathrm{M} \in\mathcal{M}(2;\mathbb{C})\,\big\vert\,(\mathrm{tr}\mathrm{M})^2-4\det\mathrm{M}=0\big\}$ est fix\'ee par $F$; on peut v\'erifier que c'est exactement l'ensemble des points fixes $\mathrm{Fix}\, F$ de $F$. \smallskip

\noindent Remarquons que la fibration $y/z=$ cte est invariante; nous verrons plus loin que ceci est un fait g\'en\'eral. 
\end{eg}

\noindent La proposition qui suit se v\'erifie par un simple calcul formel. Elle donne en particulier une version effective des \'enonc\'es \ref{birationnelle}, \ref{polynomiale} et \ref{corollaire}.

\begin{pro}\label{calculformel}
L'application rationnelle $\tau$-\'equivariante $\eta\colon(\lambda_1, \lambda_2) \dashrightarrow(\Psi(\lambda_1,\lambda_2),\Psi(\lambda_2,\lambda_1))$ se rel\`eve en l'application $\Phi\colon\mathcal{M}(2;\mathbb{C})\dashrightarrow\mathcal{M}(2;\mathbb{C})$ compatible \`a la conjugaison d\'efinie par

\begin{footnotesize}
$$\left[\begin{array}{cc}x & y \\ z & t\end{array}\right]\mapsto\left[ \begin{array}{cc} \frac{1}{2}\left(\Psi(\xi_1,\xi_2)+\Psi(\xi_2,\xi_1)+\frac{(t-x)}{\Delta}(\Psi(\xi_2,\xi_1)-\Psi(\xi_1,\xi_2))\right) & \frac{y}{\Delta}\Big(\Psi(\xi_1,\xi_2)-\Psi(\xi_2,\xi_1)\Big)\\ \frac{z}{\Delta}\Big(\Psi(\xi_1,\xi_2)-\Psi(\xi_2,\xi_1)\Big)& \frac{1}{2}\left(\Psi(\xi_1,\xi_2)+\Psi(\xi_2,\xi_1)+\frac{(t-x)}{\Delta}(\Psi(\xi_1,\xi_2)-\Psi(\xi_2,\xi_1))\right)\end{array}\right]$$
\end{footnotesize} 

\noindent o\`u $\Delta=\sqrt{(t-x)^2+4yz}$, $\xi_1=\frac{1}{2}\left(t+x+\Delta\right)$ et $\xi_2=\frac{1}{2}\left(t+x-\Delta\right)$.
\end{pro}

\noindent Dans les formules ci-dessus nous choisissons une d\'etermination de la racine de $(t-x)^2+4yz$. Un calcul de monodromie \'el\'ementaire montre que les composantes de $\Phi$ sont uniformes. Elles sont rationnelles (resp. polynomiales) si $\eta$ l'est.

\subsection{Exemples}

\begin{ap}
\begin{itemize}
\item Les applications $\eta\colon\mathcal{D}\to\mathcal{D}$ lin\'eaires $\tau$-\'equivariantes sont du type $\eta(\lambda_1,\lambda_2)=(a\lambda_1+b\lambda_2,b\lambda_1+a\lambda_2)$. Une telle application se rel\`eve en 
$$\left[\begin{array}{cc}x & y \\ z & t\end{array}\right]\mapsto\left[ \begin{array}{cc} ax+bt & y(a-b)\\z(a-b) & at+bx\end{array}\right].$$

\noindent Ces transformations forment une alg\`ebre commutative de dimension $2$. Remarquons que $\eta$ est inversible si et seulement si $a^2-b^2\not=0$; les inversibles constituent un groupe ab\'elien. 

\item Une transformation homographique $\eta$ commute \`a l'involution $\tau$ si et seulement si elle est de la forme
$$\eta(\lambda_1,\lambda_2)=\left(\frac{a\lambda_1+b\lambda_2+c}{A(\lambda_1+\lambda_2)+C},\frac{a\lambda_2+b\lambda_1+c}{A(\lambda_1+\lambda_2)+C}\right).$$ Elle se rel\`eve en $$\left[\begin{array}{cc}x & y \\ z & t\end{array}\right]\mapsto\left[ \begin{array}{cc}\frac{ax+bt+c}{A(x+t)+C} & \frac{y(a-b)}{A(x+t)+C}\\\frac{z(a-b)}{A(x+t)+C} & \frac{at+bx+c}{A(x+t)+C}\end{array}\right].$$

\noindent L'ensemble de ces relev\'es (inversibles) forme un groupe qui est le projectivis\'e du groupe $\mathrm{G}$ engendr\'e par les matrices inversibles $$\left[\begin{array}{ccccc}a & b& c & 0 & 0\\ b & a &c & 0 & 0\\ A & A & C& 0 & 0\\ 0 & 0 & 0 & a-b &0\\ 0 & 0 & 0 & 0 & a-b\end{array}\right];$$ son alg\`ebre $\mathfrak{g}$ est engendr\'ee par les matrices $\mathrm{X}$, $\mathrm{Y}$, $\mathrm{Z}$, $\mathrm{U}$ et $\mathrm{V}$ donn\'ees par
\begin{footnotesize}
\begin{align*}
& \mathrm{X}=\left[\begin{array}{ccccc}
 1& 1& 0& 0& 0\\
 1& 1& 0& 0& 0\\
 0& 0& -2& 0& 0\\
 0& 0& 0& 0& 0\\
 0& 0& 0& 0& 0
\end{array}\right], && \mathrm{Y}=\left[\begin{array}{ccccc}
0 & 0& 1& 0& 0\\
0 & 0& 1& 0& 0\\
0 & 0& 0& 0& 0\\
0 & 0& 0& 0& 0\\
0 & 0& 0& 0& 0
\end{array}\right], && \mathrm{Z}=\left[\begin{array}{ccccc}
0 & 0& 0& 0& 0\\
0 & 0& 0& 0& 0\\
1 & 1& 0& 0& 0\\
0 & 0& 0& 0& 0\\
0 & 0& 0& 0& 0
\end{array}\right], \\
& \mathrm{U}=\left[\begin{array}{ccccc}
1 &0 &0 &0 &0 \\
0 &1 &0 &0 &0 \\
0 &0 &0 &0 &0 \\
0 &0 &0 &1 &0 \\
0 &0 &0 &0 &1 
\end{array}\right], && \mathrm{V}=\left[\begin{array}{ccccc}
0 & 1& 0& 0& 0\\
1 & 0& 0& 0& 0\\
0 & 0& 0& 0& 0\\
0 & 0& 0& -1& 0\\
0 & 0& 0& 0& -1
\end{array}\right]. &&
\end{align*}
\end{footnotesize}
 
\noindent L'alg\`ebre engendr\'ee par $\mathrm{X}$, $\mathrm{Y}$ et $\mathrm{Z}$ est isomorphe \`a $\mathfrak{sl}(2;\mathbb{C})$. On peut v\'erifier que $\mathfrak{g}$ est un produit semi-direct de l'alg\`ebre commutative engendr\'ee par $\mathrm{U}$, $\mathrm{V}$ et de celle engendr\'ee par $\mathrm{X}$, $\mathrm{Y}$ et $\mathrm{Z}$.
\end{itemize}
\end{ap}

\noindent Le Th\'eor\`eme \ref{fibration}, ou le fait qu'il n'y ait pas d'automorphisme de type \textsc{H\'enon} g\'en\'eralis\'e commutant \`a $(y,x)$, ne permet pas de construire un automorphisme de $\mathcal{M}(2;\mathbb{C})$ compatible \`a la conjugaison dont la restriction \`a $\mathcal{D}$ soit de type \textsc{H\'enon} g\'en\'eralis\'e. Par contre on peut construire des transformations birationnelles compatibles dont la restriction \`a~$\mathcal{D}$ est birationnellement conjugu\'ee \`a un automorphisme de \textsc{H\'enon}. Dans l'exemple qui suit nous \og tordons\fg\, un automorphisme de \textsc{H\'enon} afin d'obtenir une transformation birationnelle qui commute \`a $\tau$, transformation qui s'\'etend \`a $\mathcal{M}(2;\mathbb{C})$.

\begin{eg}\label{henon}
Soit $\eta$ la transformation birationnelle d\'efinie par
\begin{align*} &\eta(x,y)=\left(\frac{3x^3-x^2y+5xy^2+y^3-10x^2-4xy-10y^2+12x+12y-8}{(x+y-2)(x-y)^2},\right.\\
&\hspace{1.7cm}\left.\frac{x^3+5x^2y-xy^2+3y^3-10x^2-4xy-10y^2+12x+12y-8}{(x+y-2)(x-y)^2}\right);
\end{align*}
\noindent notons que $\eta$ et $\tau$ commutent. Pour obtenir $\eta$ nous avons conjugu\'e l'automorphisme de \textsc{H\'enon} $h=(y,y^3-~x)$, automorphisme qui commute \`a $(-x,-y)$,  par la transformation $\left(\frac{x+y+1}{x},\frac{x+y-1}{x}\right)$ qui conjugue $(-x,-y)$ \`a $\tau$. Ce $\eta$ se rel\`eve en l'application $\Phi$ de $\mathcal{M}(2;\mathbb{C})$ dans lui-m\^eme donn\'ee par \begin{small}$$\left[\begin{array}{cc} x & y\\ z & t\end{array}\right]\mapsto\left[\begin{array}{cc}\frac{3x^3+t^3+5xt^2+8xyz-x^2t-4xt-10(t^2+x^2)-16yz+12(t+x)-8}{((t-x)^2+4yz)(t+x-2)}&\frac{2y}{t+x-2}\\\frac{2z}{t+x-2}& \frac{x^3+3t^3-xt^2+5x^2t+8yzt-10(x^2+t^2)-4xt-16yz+12(t+x)-8}{((t-x)^2+4yz)(t+x-2)}\end{array}\right].$$\end{small}

\noindent La transformation birationnelle $\Phi$ h\'erite des propri\'et\'es de l'automorphisme de \textsc{H\'enon}: en particulier l'ensemble de \textsc{Julia} de $h$ induit un ensemble de \og type \textsc{Julia}\fg\, invariant par $\Phi$ et l'action adjointe de $\mathrm{GL}(2;\mathbb{C})$ sur~$\mathcal{M}(2; \mathbb{C})$.
\end{eg}

\subsection{Fibrations, feuilletages, tissus, coniques invariants et application squelette}

\noindent Revenons encore \`a $\Phi_\mathrm{Id}(\mathrm{M})=\mathrm{M}^2$ pour laquelle nous allons d\'egager un ensemble de fibrations invariantes pour certaines \og universelles\fg. Comme on l'a vu, le pinceau d'hyperplans associ\'e \`a la fibration $y/z$ est invariant par~$\Phi_\mathrm{Id}$ fibre \`a fibre, ce fait se g\'en\'eralise \`a toute application compatible \`a la conjugaison (Proposition \ref{calculformel}):

\begin{pro}
Soit $\Phi\colon\mathcal{M}(2;\mathbb{C})\dashrightarrow\mathcal{M}(2;\mathbb{C})$ une application rationnelle compatible \`a la conjugaison. La fibration donn\'ee par les niveaux de $y/z$ et le feuilletage associ\'e au champ $y\frac{\partial}{\partial y}-z\frac{\partial}{\partial z}$ sont invariants par $\Phi$.
\end{pro}

\begin{rem}\label{fiber}
Sous les hypoth\`eses de la proposition pr\'ec\'edente le pinceau d'hyperplans associ\'e \`a la fibration $y/z$ est invariant par $\Phi$ fibre \`a fibre, \emph{i.e.} les orbites de $\Phi$ sont contenues dans les hyperplans $y/z=$~cte. 

\noindent Le flot du champ de vecteurs $y\frac{\partial}{\partial y}-z\frac{\partial}{\partial z}$ est donn\'e par $(x,y\mathrm{e}^s,z \mathrm{e}^{-s},t)$. Ses trajectoires sont les fibres de la fibration en coniques $\mathcal{M}(2;\mathbb{C})\to\mathbb{C}^3$, $(x,y,z,t)\mapsto(x,yz,t)$. La fibration pr\'ec\'edente est conserv\'ee globalement par $\Phi$ mais pas fibre \`a fibre.
\end{rem}

\noindent Notons que $\Phi_\mathrm{Id}$ laisse aussi la fibration $\frac{x-t}{y}=$ cte invariante fibre \`a fibre. 

\smallskip

\noindent Consid\'erons dans $\mathcal{M}(2;\mathbb{C})$ la droite $\mathcal{C}=\big\{\lambda\mathrm{Id}\,\big\vert\,\lambda\in\mathbb{C} \big\}$, \emph{i.e.} le centre de $\mathcal{M}(2;\mathbb{C})$ pour sa loi d'alg\`ebre. \`A~$\mathcal{C}$ on peut associer la fibration $\mathcal{P}$ en $2$-plans d\'efinie comme suit: si~$\mathrm{M}$ appartient \`a $\mathcal{M}(2;\mathbb{C})\setminus\mathcal{C}$, on d\'efinit~$\mathcal{P}(\mathrm{M})$ comme l'unique plan contenant $\mathrm{M}$ et $\mathcal{C}$. Si~$\mathrm{M}$ est un \'el\'ement de~$\mathcal{M}(2;\mathbb{C})$, alors~$\mathcal{P}(\mathrm{M})$ est exactement l'ensemble des matrices qui commutent \`a $\mathrm{M}$. On peut aussi remarquer que le plan g\'en\'erique $\mathcal{P}(\mathrm{M})$ est un conjugu\'e de~$\mathcal{D}:$ si $\mathrm{P}$ diagonalise $\mathrm{M}$, \emph{i.e.} $\mathrm{P}\mathrm{M}\mathrm{P}^{-1}$ appartient \`a $\mathcal{D}$, alors~$\mathcal{P}(\mathrm{M})$ est inclu dans $\mathrm{P}\mathcal{D}\mathrm{P}^{-1}$ et donc \'egal. Ceci induit un feuilletage singulier en $2$-plans qui est invariant par~$\Phi_\mathrm{Id}$ feuille \`a feuille. En effet soit $s\mathrm{M}+t \mathrm{Id}$ dans~$\mathcal{P}(\mathrm{M})$; on constate que $$\Phi_\mathrm{Id}(s\mathrm{M}+t \mathrm{Id})=s^2 \mathrm{M}^2+2st\mathrm{M}+t^2\mathrm{Id}$$ commute visiblement \`a $\mathrm{M}$. Par suite les orbites de $\Phi_\mathrm{Id}$ sont dans les $2$-plans du pinceau $\mathcal{P}$. Ainsi l'application $\mathbb{P}(\Phi_\mathrm{Id})\colon\mathbb{P}^3(\mathbb{C}) \dashrightarrow~\mathbb{P}^3(\mathbb{C})$ laisse invariante la fibration en droites $\mathbb{P} (\mathcal{P}):$ ce sont les droites de $\mathbb{P}^3(\mathbb{C})$ passant par $[\mathrm{Id}]=(1:0: 0:1)$.  Nous avons des propri\'et\'es analogues pour les transformations $\Phi_{\alpha,d}$ homog\`enes compatibles \`a la conjugaison. 

\noindent Si $\mathrm{M}$ est un point g\'en\'erique de $\mathcal{M}(2;\mathbb{C})$ l'adh\'erence de \textsc{Zariski} de son orbite par $\Phi_\mathrm{Id}$ est un $2$-plan (car c'est le cas pour un \'el\'ement g\'en\'erique de $\mathcal{D}$). Par cons\'equent il n'y a qu'une seule fibration en surfaces invariante fibre \`a fibre par $\Phi_\mathrm{Id}$, la fibration $\mathcal{P}$. Notamment la fibration $\mathcal{P}$ est donn\'ee par les niveaux de l'application rationnelle $\left( \frac{y}{z},\frac{x-t}{y}\right)$ ce que l'on peut voir de fa\c{c}on directe.

\smallskip

\noindent En fait l'invariance de $\mathcal{P}$ se g\'en\'eralise \`a toutes les applications compatibles \`a la conjugaison.

\begin{pro}
Soit $\Phi\colon\mathcal{M}(2;\mathbb{C})\dashrightarrow\mathcal{M}(2;\mathbb{C})$ une application rationnelle compatible \`a la conjugaison;~$\Phi$ pr\'eserve la fibration en $2$-plans $\mathcal{P}$ fibre \`a fibre.
\end{pro}

\begin{proof}[D\'emonstration]
D'apr\`es le Lemme \ref{invdiag} nous avons $\Phi(\mathcal{D}\setminus\mathrm{Ind}\,\Phi)\subset \mathcal{D}$. Comme un plan g\'en\'erique $\mathcal{P}(\mathrm{M})$ de $\mathcal{P}$ est un conjugu\'e de $\mathcal{D}$, $\mathcal{P}(\mathrm{M})$ est invariant par $\Phi$, d'o\`u le r\'esultat.
\end{proof}

\begin{rem}
Comme on l'a vu, par d\'efinition de la fibration $\mathcal{P}$, $\Phi$ pr\'eserve $\mathcal{P}$ fibre \`a fibre si et seulement si $\Phi(\mathrm{M})$ commute \`a $\mathrm{M}$ pour tout $\mathrm{M}$ dans $\mathcal{M}(2;\mathbb{C})$. Ceci  permet de v\'erifier qu'il existe des applications pr\'eservant $\mathcal{P}$ qui ne sont pas compatibles \`a la conjugaison.
\end{rem}

\noindent La proposition qui suit est cons\'equence directe de la d\'efinition de la compatibilit\'e.

\begin{pro}
Soit $\Phi\colon\mathcal{M}(2;\mathbb{C})\dashrightarrow\mathcal{M}(2;\mathbb{C})$ une application rationnelle compatible \`a la conjugaison. Alors $\Phi$ laisse invariant le feuilletage associ\'e \`a la fibration $\mathrm{Inv}$. Plus pr\'ecis\'ement il existe une application rationnelle $\mathrm{Sq}\,\Phi\colon\mathbb{C}^2\dashrightarrow\mathbb{C}^2$ telle que l'on ait le diagramme commutatif suivant
$$\xymatrix{\mathcal{M}(2;\mathbb{C})\ar@{-->}[rr]^\Phi\ar[d]_{\mathrm{Inv}} & & \mathcal{M}(2;\mathbb{C})\ar[d]^{\mathrm{Inv}}\\
\mathbb{C}^2\ar@{-->}[rr]_{\mathrm{Sq}\,\Phi} & & \mathbb{C}^2}$$
\end{pro}

\begin{rem}
Les fibres de $\mathrm{Inv}$ sont des surfaces quadratiques: $x+t=$ cte, $xt-yz=$ cte.
\end{rem}

\begin{rem}
Il y a des applications qui pr\'eservent la fibration $\mathrm{Inv}$ sans \^etre compatible \`a la conjugaison: par exemple soit $\mathrm{M}\mapsto P(\mathrm{M})$ une application polynomiale de $\mathcal{M}(2;\mathbb{C})$ dans lui-m\^eme telle que $\det P\not\equiv 0$. Alors l'application $\mathrm{M}\mapsto (P(\mathrm{M}))^{-1} \mathrm{M} P(\mathrm{M})$ respecte la fibration $\mathrm{Inv}$ mais n'est pas, en g\'en\'eral, compatible \`a la conjugaison.
\end{rem}

\noindent Inversement on peut se demander \`a quelle condition une application rationnelle de $\mathbb{C}^2$ se rel\`eve \`a $\mathcal{M}(2;\mathbb{C})$ via $\mathrm{Inv}$. La r\'eponse est donn\'ee par l'\'enonc\'e suivant:

\begin{pro}\label{inv}
Soit $W\colon\mathbb{C}^2\dashrightarrow\mathbb{C}^2$, $(u,v)\dashrightarrow(S(u,v),T(u,v))$ une application rationnelle. Alors $W$ s'\'ecrit $\mathrm{Sq}\,\Phi$ pour une certaine application rationnelle $\Phi\colon\mathcal{M}(2;\mathbb{C}) \dashrightarrow\mathcal{M}(2;\mathbb{C})$ compatible \`a la conjugaison si et seulement si $(S^2-4 T)\circ\mathrm{Inv}$ est un carr\'e dans le corps des fonctions rationnelles $\mathbb{C}(u,v)$.
\end{pro}

\begin{proof}[D\'emonstration]
Par restriction aux matrices diagonales il suffit de r\'esoudre en $\Psi$:
\begin{align*}
& S(\lambda_1+\lambda_2,\lambda_1\lambda_2)=\Psi(\lambda_1,\lambda_2)+\Psi(\lambda_2,\lambda_1), && T(\lambda_1+\lambda_2,\lambda_1\lambda_2)=\Psi(\lambda_1,\lambda_2)\Psi(\lambda_2,\lambda_1).
\end{align*}

\noindent Par \'elimination de $\Psi(\lambda_2,\lambda_1)$ nous obtenons $$\Psi^2(\lambda_1, \lambda_2)-S(\lambda_1+\lambda_2,\lambda_1\lambda_2)\Psi(\lambda_1,\lambda_2)+T(\lambda_1+ \lambda_2,\lambda_1\lambda_2)=0$$ d'o\`u le r\'esultat.
\end{proof}

\begin{eg}\label{firs}
L'application $\mathrm{Sq}\,\Phi_\mathrm{Id}$ est donn\'ee par $\mathrm{Sq}\,\Phi_\mathrm{Id}(u,v)=(u^2-2v,v^2)$. Alors que l'extension $\widetilde{\Phi_\mathrm{Id}}$ de $\Phi_\mathrm{Id}$ \`a~$\mathbb{P}^4(\mathbb{C})$ n'est pas un endomorphisme (puisque la quadrique $\big\{\mathrm{M}\in\mathcal{M}(2;\mathbb{C})\,\big\vert\,\mathrm{Inv}(\mathrm{M})=(0,0)\big\}$ est d'ind\'etermination), celle de $\mathrm{Sq}\,\Phi_\mathrm{Id}$ \`a $\mathbb{P}^2 (\mathbb{C})$ l'est.

\noindent Notons que la conique $v=\frac{u^2}{4}$ est invariante par $\mathrm{Sq}\,\Phi_\mathrm{Id}$: c'est l'image par $\mathrm{Inv}$ de l'hypersurface discriminante $$\big\{\mathrm{M}\in\mathcal{M}(2; \mathbb{C})\,\big\vert\,4\det\,\mathrm{M}-(\mathrm{tr}\,\mathrm{M})^2=0\big\};$$ c'est donc aussi l'image de l'orbite sous l'action adjointe du $2$-plan: $$\left\{\lambda\mathrm{Id}+\mu\left[\begin{array}{cc} 0 &1\\0 & 0\end{array}\right]\,\Big\vert\,\lambda,\,\mu\in\mathbb{C}\right\}=\mathcal{P}\left(\left[\begin{array}{cc}0 & 1 \\ 0 & 0\end{array}\right]\right).$$ Ses tangentes d\'efinissent un \og $2$-tissu\fg\, compl\`etement invariant par $\mathrm{Sq}\,\Phi_\mathrm{Id}$ (\emph{i.e.} invariant par images directe et r\'eciproque):

\begin{center}
\begin{picture}(0,0)%
\includegraphics{tissub.pstex}%
\end{picture}%
\setlength{\unitlength}{3947sp}%
\begingroup\makeatletter\ifx\SetFigFont\undefined%
\gdef\SetFigFont#1#2#3#4#5{%
  \reset@font\fontsize{#1}{#2pt}%
  \fontfamily{#3}\fontseries{#4}\fontshape{#5}%
  \selectfont}%
\fi\endgroup%
\begin{picture}(1224,1374)(1189,-1123)
\end{picture}%

\end{center}
\end{eg}

\begin{eg}\label{minv}
Consid\'erons cette fois l'application $\Phi$ compatible \`a la conjugaison d\'efinie par $\Phi(\mathrm{M})=\mathrm{M}+\mathrm{M}^{-1}$. Son extension $\mathbb{P}(\Phi)\colon\mathbb{P}^4(\mathbb{C})\dashrightarrow\mathbb{P}^4(\mathbb{C})$ est la transformation de degr\'e trois
\begin{align*}
& \mathbb{P}(\Phi)\colon (x: y:z:t:w)\mapsto(x\delta+tw^2:y(\delta-w^2):z(\delta-w^2):t\delta+xw^2: w\delta), &&\delta=xt-yz;
\end{align*}

\noindent elle pr\'eserve la fibration $y/z=$ cte. Il y a un seul point d'ind\'etermination \`a distance finie; \`a l'infini l'ensemble d'ind\'etermination est donn\'e par $\delta=w=0$. Une matrice de d\'eterminant nul va \`a l'infini sur sa matrice de cofacteurs. 

\noindent La restriction de $\Phi$ \`a $\mathcal{D}$ induit l'application $(\lambda_1,\lambda_2)\mapsto\left(\lambda_1+\lambda_1^{-1},\lambda_2+\lambda_2^{-1}\right)$. Nous en d\'eduisons l'application $\mathrm{Sq}\,\Phi\colon\mathbb{C}^2\dashrightarrow\mathbb{C}^2$ donn\'ee par
$$(u,v)\mapsto\left(\frac{u}{v}+u,v+\frac{u^2}{v}-2+ \frac{1}{v}\right).$$
Elle se prolonge en l'endomorphisme $(u,v,w)\mapsto(u(v+w):(v-w)^2+u^2:vw)$ de $\mathbb{P}^2(\mathbb{C})$ qui pr\'eserve uniquement deux droites (celle d'\'equation $u=0,$ resp. $w=0$). Ceci donne un autre exemple qui montre que l'application $\mathrm{Sq}$ peut avoir une extension holomorphe \`a $\mathbb{P}^2(\mathbb{C})$ alors que l'extension $\widetilde{\Phi}$ de $\Phi$ \`a $\mathbb{P}^4(\mathbb{C})$ est seulement rationnelle. La courbe $\mathcal{C}$ d'\'equation $v=\frac{u^2}{4}$ est encore invariante par $\mathrm{Sq}\,\Phi$. Comme dans l'Exemple \ref{firs} la famille des tangentes \`a $\mathcal{C}$ est un $2$-tissu invariant par $\mathrm{Sq}\,\Phi$. La famille des fonctions rationnelles $\lambda\left(z+\frac{1}{z}\right)$ sur~$\mathbb{P}^1( \mathbb{C})$ a \'et\'e \'etudi\'ee en d\'etail dans \cite{Mi}. Pour deux valeurs de~$\lambda$ sp\'eciales $\Big(\lambda=\pm\frac{\mathrm{i}}{2}\Big)$ il s'agit d'un exemple de \textsc{Latt\`es}. L'application induite correspondante $\zeta_\lambda\colon\mathcal{M}(2;\mathbb{C}) \dashrightarrow\mathcal{M} (2;\mathbb{C})$ d\'efinie par $\zeta_\lambda(\mathrm{M})=\lambda( \mathrm{M}+\mathrm{M}^{-1})$ h\'erite de ses propri\'et\'es; par exemple pour $\lambda=\pm\frac{\mathrm{i}}{2}$ l'ensemble $\mathrm{Per}\,\zeta_\lambda$ des points p\'eriodiques de $\zeta_\lambda$ est dense dans $\mathcal{M}(2;\mathbb{C})$.
\end{eg}

\begin{eg}\label{sssq}
Soit $\varsigma_c\colon\mathcal{M}(2;\mathbb{C})\to\mathcal{M}(2;\mathbb{C})$ d\'efinie par $\varsigma_c(\mathrm{M})=\mathrm{M}^2+c\mathrm{Id}$. La restriction de $\varsigma_c$ \`a~$\mathcal{D}$ induit l'application $(\lambda_1,\lambda_2)\mapsto(\lambda_1^2+c,\lambda_2^2+c)$ d'o\`u l'application rationnelle
\begin{align*}
&\mathrm{Sq}\,\varsigma_c\colon\mathbb{C}^2\dashrightarrow\mathbb{C}^2, &&(u,v)\mapsto(u^2-2v+2c,v^2+c(u^2-2v)+c^2);
\end{align*}

\noindent elle se prolonge en un endomorphisme de $\mathbb{P}^2(\mathbb{C})$. La conique $\mathcal{C}$ d'\'equation $v=\frac{u^2}{4}$ est invariante par $\mathrm{Sq}\,\varsigma_c$ et la dynamique de $\mathrm{Sq}\,\varsigma_{c\vert\mathcal{C}}$ est conjugu\'ee \`a $z\mapsto \frac{z^2}{2}+2c$ qui est conjugu\'ee \`a $z\mapsto z^2+c$ (via $z\mapsto\frac{z}{2}$).

\noindent De plus, comme toute application compatible, $\varsigma_c$ laisse $y/z=$ cte et $\mathcal{P}$ invariantes fibre \`a fibre.

\noindent Ici encore l'hyperplan $\mathfrak{sl}(2;\mathbb{C})$ est contract\'e sur la droite $\mathbb{C}\cdot\mathrm{Id};$ plus pr\'ecis\'ement $\varsigma_c(x,y,z,-x)=(x^2+yz)\mathrm{Id}$ et la restriction de $\varsigma_c$ \`a $\mathbb{C}\cdot\mathrm{Id}$ se traite \'evidemment comme $z\mapsto z^2+c$. 
\end{eg}

\noindent Au vu des exemples pr\'ec\'edents on peut se demander si le fait que la courbe $\mathcal{C}$ d'\'equation $v=\frac{u^2}{4}$ soit invariante par $\mathrm{Sq}\,\Phi$ est un fait g\'en\'eral. La proposition qui suit donne une r\'eponse partielle.

\begin{pro}
Soit $\Phi\colon\mathcal{M}(2;\mathbb{C})\to\mathcal{M}(2;\mathbb{C})$ un automorphisme polynomial compatible \`a la conjugaison; l'hypersurface discriminante $\Delta=\big\{\mathrm{M}\in\mathcal{M}(2;\mathbb{C})\,\big\vert\, 4\det\mathrm{M}-(\mathrm{tr}\,\mathrm{M})^2=0\big\}$ est invariante par $\Phi$ d'o\`u l'invariance de la conique $v=\frac{u^2}{4}$ par $\mathrm{Sq}\,\Phi$.
\end{pro}

\begin{proof}[D\'emonstration]
Raisonnons par l'absurde: supposons que l'hypersurface $\Delta$ ne soit pas invariante. Il existe donc $\mathrm{M}_0$ non diagonalisable dans $\Delta$ tel que $\Phi(\mathrm{M}_0)$ soit diagonalisable d'o\`u l'existence de $\mathrm{M}'_0$ dans $\Delta$ tel que $\Phi(\mathrm{M}'_0)$ soit diagonale. L'ensemble $\mathcal{D}\mathrm{M}'_0\mathcal{D}^{-1}$ est de dimension sup\'erieure ou \'egale \`a $1$ et, pour tout $\mathrm{D}$ dans $\mathcal{D}$, nous avons~$\Phi(\mathrm{D}\mathrm{M}'_0\mathrm{D}^{-1})=\mathrm{D}\Phi(\mathrm{M}'_0)\mathrm{D}^{-1} =\Phi(\mathrm{M}'_0)$: contradiction.
\end{proof}

\begin{rems}
Cette d\'emonstration s'\'etend aux applications polynomiales $\Phi$ de $\mathcal{M}(2;\mathbb{C})$ dans lui-m\^eme qui ne contractent pas de courbe sur un point. Remarquons aussi qu'une application de la forme $\frac{P(\mathrm{M})}{4\det\mathrm{M}-(\mathrm{tr}\,\mathrm{M})^2}$, o\`u $P$ d\'esigne une application polynomiale compatible \`a la conjugaison, n'est pas bien d\'efinie sur $\Delta$. 
\end{rems}

\noindent L'Exemple \ref{minv} montre que cette conique peut aussi \^etre invariante dans le cas non injectif.

\begin{pro}
Soit $\Phi\colon\mathcal{M}(2;\mathbb{C})\dashrightarrow\mathcal{M}(2;\mathbb{C})$ une transformation compatible \`a la conjugaison. Alors~$\Phi$ est birationnelle si et seulement si $\mathrm{Sq}\,\Phi$ l'est.
\end{pro}

\begin{proof}[D\'emonstration]
Supposons $\mathrm{Sq}\,\Phi$ non (g\'en\'eriquement) injective. Ceci signifie que $\Phi$ envoie (au moins) deux orbites de l'action adjointe sur une seule et ceci g\'en\'eriquement; en particulier $\Phi$ est non injective. 

\noindent R\'eciproquement on sait que $\Phi$ est birationnelle si et seulement si $\Phi_{\vert\mathcal{D}}$ l'est. Si on identifie $\Phi_{\vert\mathcal{D}}$ \`a une application rationnelle $\varphi\colon(x,t)\mapsto (a(x,t),b(x,t))$ de $\mathbb{C}^2$ dans lui-m\^eme, l'injectivit\'e de $\mathrm{Sq}\,\Phi$ se traduit par celle de~$(x,t)\mapsto\big(a(x,t)+b(x,t),a(x,t)b(x,t)\big)$ et donc celle de $\varphi$; la transformation $\Phi$ est donc birationnelle.
\end{proof}

\begin{rem}
Soit $\Phi\colon\mathcal{M}(2;\mathbb{C})\dashrightarrow\mathcal{M}(2;\mathbb{C})$ une transformation compatible \`a la conjugaison. Si l'application $\mathrm{Sq}\,\Phi$ est triviale, alors $\Phi$ est triviale ou $\Phi\colon\mathrm{M}\mapsto(\det\mathrm{M})\mathrm{M}^{-1}$. En effet si $\mathrm{Sq}\,\Phi=\mathrm{id},$ alors $\Phi_{\vert\mathcal{D}}$ co\"{i}ncide avec $(x,t)\mapsto(x,t)$ ou $(x,t)\mapsto(t,x);$ la Proposition \ref{calculformel} permet de conclure.
\end{rem}

\section{Dynamique de $\Phi_\mathrm{Id}$}\label{Phi2}

\noindent Soit $\Phi\colon\mathcal{M}(2;\mathbb{C})\dashrightarrow\mathcal{M}(2;\mathbb{C})$ une application rationnelle; notons $\mathrm{K}(\Phi)$ le corps des fonctions rationnelles invariantes par~$\Phi$. Un \'el\'ement de $\mathrm{K}(\Phi)$ est une fonction rationnelle $f\colon\mathcal{M}(2;\mathbb{C})\dashrightarrow\mathbb{P}^1(\mathbb{C})$ telle que $f\circ\Phi=f$.

\noindent Revenons \`a l'application $\Phi_\mathrm{Id}\colon\mathcal{M}(2;\mathbb{C}) \to\mathcal{M}(2;\mathbb{C}),\,\mathrm{M}\mapsto\mathrm{M}^2$: $$\left[\begin{array}{cc} x & y\\ z & t \end{array}\right]\to \left[\begin{array}{cc} x^2+yz& y(x+t)\\ z(x+t)& t^2+yz\end{array}\right].$$ Comme $\Phi_\mathrm{Id}$ laisse le feuilletage en $2$-plans $\mathcal{P}$ invariant, chacune de ses orbites est contenue dans un certain~$2$-plan de $\mathcal{P}$. Au niveau alg\'ebrique ceci se formalise de la fa\c{c}on suivante: 

\begin{thm}
Les fonctions invariantes par $\Phi_\mathrm{Id}$ sont engendr\'ees par $\frac{y}{z}$ et $\frac{x-t}{z}$, {\it i.e.} $\mathrm{K}(\Phi_\mathrm{Id})=\mathbb{C}\left(\frac{y}{z},\frac{x-t}{z}\right)$.
\end{thm}

\begin{proof}[D\'emonstration]
Soit $f$ une fonction rationnelle invariante par $\Phi_\mathrm{Id}$. Il suffit de montrer que $f$ est constante sur chaque $2$-plan de $\mathcal{P}$. Plus pr\'ecis\'ement il suffit de l'\'etablir pour un ensemble dense de tels $2$-plans, par exemple sur les $2$-plans du type suivant $$\Pi(\mathrm{P})=\left\{ \mathrm{P}\,\mathrm{diag}(x,t)\,\mathrm{P}^{-1}\,\Big\vert (x,t)\in \mathbb{C}^2\right\}$$ o\`u $\mathrm{P}$ est une matrice inversible fix\'ee. La condition d'invariance implique les \'egalit\'es suivantes 
\begin{align*}
&f\left(\mathrm{P}\,\mathrm{diag}(x^{2^n},t^{2^n})\,\mathrm{P}^{-1}\right)=f\left(\mathrm{P}\,\mathrm{diag}(x,t)\,\mathrm{P}^{-1}\right), && \text{ pour tout $n$ dans $\mathbb{Z}$.}
\end{align*}
Ainsi $f$ est constante sur l'orbite $\mathcal{O}(\mathrm{P}\,\mathrm{diag}(x_0,t_0)\,\mathrm{P}^{-1};\Phi_\mathrm{Id})$ de $\mathrm{P}\,\mathrm{diag}(x_0,t_0)\,\mathrm{P}^{-1}$ par $\Phi_\mathrm{Id}$ qui, pour un choix g\'en\'erique de $(x_0, t_0),$ est \textsc{Zariski} dense dans $\Pi(\mathrm{P})$. Par suite $f$ est constante sur chaque $\Pi(\mathrm{P})$ et donc sur chaque $2$-plan de $\mathcal{P}$.
\end{proof}

\noindent Si $\Phi\colon\mathcal{M}(2;\mathbb{C})\dashrightarrow\mathcal{M}(2;\mathbb{C})$ est compatible \`a la conjugaison, sa dynamique h\'erite de cette compatibilit\'e. Par exemple si $\mathrm{M}$ est p\'eriodique pour $\Phi$, \emph{i.e.} $\Phi^k(\mathrm{M})=\mathrm{M}$ pour un certain $k$, alors toute l'orbite $$\big\{\mathrm{A}\mathrm{M} \mathrm{A}^{-1}\,\big\vert\,\mathrm{A} \in\mathrm{GL}(2;\mathbb{C})\big\}$$ de $\mathrm{M}$ est bien s\^ur p\'eriodique. De m\^eme tout ensemble $S$ invariant par $\Phi$ a son satur\'e $\big\{\mathrm{A}S\mathrm{A}^{-1}\,\big\vert\,\mathrm{A} \in\mathrm{GL}(2;\mathbb{C})\big\}$ invariant par $\Phi$.

\subsection{Points fixes, points p\'eriodiques}

\noindent Nous nous int\'eressons aux points p\'eriodiques de l'application $\Phi_\mathrm{Id}$. Comme $\Phi_\mathrm{Id}$ est compatible avec la conjugaison il suffit de tester les matrices de \textsc{Jordan} du type $\left[\begin{array}{cc}\lambda_1& 0\\0&\lambda_2\end{array} \right]$ (c'est-\`a-dire \'etudier les points p\'eriodiques de~$\Phi_{\mathrm{Id}_{\vert\mathcal{D}}}$) et $\left[\begin{array}{cc}\lambda& 1\\0&\lambda\end{array}\right]$. On v\'erifie que les points fixes de $\Phi_\mathrm{Id}$ sont~$\mathbf{0}$,~$\mathrm{Id}$ et les matrices conjugu\'ees \`a $\left[\begin{array}{cc}1& 0\\0&0\end{array}\right]$; autrement dit les points fixes sont $$\big\{\mathbf{0},\,\mathrm{Id}\big\}\cup\big\{\mathrm{M}\in\mathcal{M}(2;\mathbb{C})\,\big\vert\,\mathrm{Inv}(\mathrm{M})=(1,0)\big\}.$$ L'\'etude de $\Phi_{\mathrm{Id}_{\vert\mathcal{D}}}$ repose sur celle de $\chi_2\colon z\mapsto z^2$ dans le plan complexe. Rappelons que l'ensemble de \textsc{Julia} de $\chi_2$ est le cercle unit\'e et que le point $0$ est super-attractant. Si $0<\vert z\vert<1$ les it\'er\'es de $z$ convergent vers~$0$; si $\vert z\vert>1$ les it\'er\'es de $z$ tendent vers l'infini (si on compactifie $\mathbb{C},$ l'infini est aussi super attracteur). Les points p\'eriodiques autres que $0$ et $\infty$ sont donc sur le cercle unit\'e. Pour d\'ecrire l'orbite de $\mathrm{e}^{2\mathrm{i}\pi \vartheta}$ sous l'action de $\chi_2$ on choisit d'\'ecrire $\vartheta$ sous forme $2$-adique: $\vartheta=\sum_{n\geq 1}\frac{\varepsilon_n}{2^n}$ avec $\varepsilon_n \in\{0,\,1\}$. On peut ainsi coder un point de $\mathbb{S}^1$ par la suite $\varepsilon=(\varepsilon_1,\, \varepsilon_2,\ldots)$. L'action de $\chi_2$ dans l'\'ecriture $2$-adique des angles est le shift $(\varepsilon_1,\varepsilon_2,\varepsilon_3,\ldots)\to(\varepsilon_2,\varepsilon_3,\varepsilon_4,\ldots)$. Un point de $\mathbb{S}^1$ associ\'e \`a une suite $\varepsilon$ p\'eriodique est lui-m\^eme p\'eriodique pour $\chi_2$; tous les points p\'eriodiques sont de cette forme. Ce sont aussi les racines $(2^n-1)$-i\`eme de l'unit\'e qui sont bien s\^ur denses dans le cercle~$\mathbb{S}^1$. Si $z$ est un point g\'en\'erique de $\mathbb{S}^1$, alors l'orbite de $z$ par $\chi_2$ est dense dans le cercle.

\bigskip

\noindent  En identifiant $\mathcal{D}$ \`a $\mathbb{C}^2$ et $\Phi_{\mathrm{Id}_{\vert\mathcal{D}}}$ \`a $(x^2,t^2)$ les points p\'eriodiques de $\Phi_{\mathrm{Id}_{\vert\mathcal{D}}}$ sont, outre $(0,0)$, de l'un des types suivants
\begin{align*}
& (x_{k_n},t_{k_n}), && (x_{k_n},0), && (0,t_{k_n})
\end{align*}

\noindent o\`u les $x_{k_n}$ et $t_{k_n}$ parcourent les racines $(2^n-1)$-i\`eme de l'unit\'e. Notons que l'ensemble $$\big\{(x,t)\in\mathbb{S}^1\times\mathbb{S}^1 \,\big\vert\, x^{2^n-1}=t^{2^n-1}=1,\, n\in\mathbb{N}\big\}$$ est dense dans le tore $\mathbb{S}^1\times\mathbb{S}^1$. Les ensembles $\mathbb{S}^1\times\big\{0\big\},$ $\big\{0\big\}\times\mathbb{S}^1$
 et $\mathbb{S}^1\times\mathbb{S}^1$ sont $\Phi_{\mathrm{Id}_{\vert\mathcal{D}}}$-invariants, de m\^eme que les ensembles $\big\{\vert x\vert\leq 1,\,\big\vert t\vert\leq 1\big\}$, $\big\{\vert x\vert\leq 1,\,\big\vert t\vert=1\big\},$ $\big\{\vert x\vert\leq 1,\,\big\vert t\vert\geq 1\big\}$ etc. Si $\mathrm{M}$ est une matrice diagonalisable, alors l'orbite $\mathcal{O}(\mathrm{M};\Phi_\mathrm{Id})$ de $\mathrm{M}$ co\"{i}ncide avec celle de $\mathrm{diag}(\lambda_1,\lambda_2)$ o\`u les $\lambda_i$ sont les valeurs propres de $\mathrm{M}$ de sorte que la description de l'orbite de $\mathrm{M}$ ne d\'epend que de ses valeurs propres. 
En particulier tous les \'el\'ements conjugu\'es aux $\mathrm{diag}(x_{k_n},t_{k_n})$ sont p\'eriodiques de p\'eriode $2^n-1$. Dans le m\^eme ordre d'id\'ee si les valeurs propres $\lambda_1$, $\lambda_2$ de $\mathrm{M}$ sont en module strictement inf\'erieures \`a $1,$ alors $\displaystyle\lim_{k\to +\infty}\Phi_\mathrm{Id}^k(\mathrm{M})=0$.

\noindent Reste \`a d\'ecrire l'orbite des matrices de type $\left[\begin{array}{cc}\lambda & \lambda \\ 0& \lambda\end{array}\right]=\lambda \left[\begin{array}{cc} 1& 1 \\ 0& 1\end{array}\right]$ dont l'it\'eration $n$i\`eme par $\Phi_\mathrm{Id}$ donne $\mathrm{M}_n(\lambda)=\lambda^{2^n}\left[\begin{array}{cc}1 &2^n\\ 0 & 1\end{array}\right]$.

\noindent Aucune valeur de $\lambda$ non triviale ne produit ici de point p\'eriodique:
\begin{itemize}
\item si $\vert\lambda\vert<1,$ alors $\mathrm{M}_n(\lambda)\to {\bf 0}$,

\item si $\vert\lambda\vert\geq 1,$ alors $\mathrm{M}_n(\lambda)$ \og s'\'echappe\fg\, dans un sens que nous pr\'ecisons maintenant. Si on compactifie $\mathcal{M}(2;\mathbb{C})\simeq\mathbb{C}^4$ par $\mathbb{P}^4(\mathbb{C})=\big\{(x:y:z:t:w)\big\},$ alors $\mathrm{M}_n (\lambda)$ correspond \`a $$(\lambda^{2^n}:\lambda^{2^n}2^n:0:\lambda^{2^n}:1)=\left(\frac{1}{2^n}:1:0:\frac{1}{2^n}:\frac{1}{\lambda^{2^n}2^n}\right)$$ qui pour $\vert\lambda\vert\geq 1$ tend vers $(0:1:0:0:0),$ point d'ind\'etermination de l'extension de~$\Phi_\mathrm{Id}$ \`a $\mathbb{P}^4(\mathbb{C})$. 
\end{itemize}

\noindent Nous d\'eduisons de ce qui pr\'ec\`ede que ${\bf 0}$ est un point super attractant. En effet si $\mathrm{M}$ est proche de~${\bf 0}$, alors $\mathrm{M}$ est conjugu\'ee \`a une matrice de \textsc{Jordan} $\left[\begin{array}{cc}\lambda_1 & 0\\ 0 & \lambda_2\end{array}\right]$ avec $\vert\lambda_i\vert<1$ ou \`a $\left[\begin{array}{cc} \lambda & \lambda \\ 0 & \lambda\end{array}\right]$ avec $\vert\lambda\vert<1$. En particulier les matrices nilpotentes sont \og pr\'e-super attractives\fg. Le bassin d'attraction $\mathrm{W}_\mathrm{Id}^s({\bf 0})$ de ${\bf 0}$ est constitu\'e des matrices ayant leurs deux valeurs propres de module strictement plus petit que $1$. Son bord est l'ensemble des matrices $\mathrm{M}$ ayant une valeur propre de module $1$ et l'autre de module plus petit ou \'egal \`a~$1$. C'est une hypersurface \textsc{Levi}-plate\footnote{\,Rappelons qu'une sous-vari\'et\'e r\'eelle $V$ de codimension $1$ dans $\mathbb{C}^n=\mathbb{R}^{2n}$ est dite \emph{\textsc{Levi}-plate} si son champ d'hyperplans tangents complexes $\mathrm{T}_m^\mathbb{C}V=\mathrm{T}_mV\cap\mathrm{i}\mathrm{T}_mV$ est int\'egrable. Ce champ induit alors un feuilletage de $V$ en sous-vari\'et\'es complexes de dimension $n-1$.} que nous allons d\'ecrire. Nous utiliserons la proposition suivante dont la preuve est \'el\'ementaire.

\begin{pro}\label{mod1}
Soit $P(z)=z^2-bz+c$ un trin\^ome du second degr\'e, $(b,c)$ d\'esignant un \'el\'ement de $\mathbb{C}^2$. Alors $P$ a une racine de module $1$ et l'autre de module plus petit que $1$ si et seulement si on a
$$\left\{\begin{array}{ll}
\vert c\vert^2-\frac{1}{2}\vert b\vert^2-\frac{1}{2}\vert b^2-4c\vert+1=0 & (i)\\
\vert c\vert \leq 1& (ii)
\end{array}
\right.$$
\end{pro}

\begin{proof}
D\'esignons par $\lambda_1$ et $\lambda_2$ les deux racines de $P$. La premi\`ere condition se traduit par $$(\vert\lambda_1\vert^2-1)(\vert\lambda_2\vert^2-1)=0;$$ en exprimant les $\lambda_i$ \`a l'aide de $b$ et $c$ nous obtenons l'\'egalit\'e $(i);$ quant \`a $(ii)$ elle est \'evidente.
\end{proof}

\noindent L'ensemble $\Sigma\subset\mathbb{C}^2$ d\'ecrit par \emph{(i)} et \emph{(ii)} est un ensemble semi-alg\'ebrique connexe dont voici une autre pr\'esentation
\begin{align*}
& \Big(\vert c\vert^2-\frac{1}{2}\vert b\vert^2+1\Big)^2
-\frac{1}{4}\Big\vert b^2-4c\Big\vert^2=0&& (i_1) \\
& \vert c\vert^2-\frac{1}{2}\vert b\vert^2+1\geq 0&& (i_2)\\
& \vert c\vert^2 \leq 1&& (ii)
\end{align*}

\noindent Cette pr\'esentation est polynomiale en les parties r\'eelles et imaginaires de $b$ et $c$. L'ensemble $\Sigma$ est invariant sous l'action de $\mathbb{S}^1$ sur $\mathbb{C}^2$ donn\'ee par $(\mathrm{e}^{\mathrm{i}\vartheta},b,c)\mapsto(b\mathrm{e}^{\mathrm{i}\vartheta},c\mathrm{e}^{2\mathrm{i}\vartheta})$. V\'erifions que $\Sigma$ est \textsc{Levi}-plat; consid\'erons l'application
\begin{align*}
& \xi\colon\mathbb{S}^1\times\mathbb{C}\to\mathbb{C}^2, && (\mathrm{e}^{\mathrm{i}\vartheta}+u,\mathrm{e}^{\mathrm{i}\vartheta}u).
\end{align*}

\noindent L'ensemble d\'ecrit par $(i_1)$ et $(i_2)$ est pr\'ecis\'ement l'image de $\xi$; il est \textsc{Levi}-plat puisque les droites $u\to(\mathrm{e}^{\mathrm{i}\vartheta}+u,\mathrm{e}^{\mathrm{i}\vartheta}u)$ sont contenues dedans. Ces droites font partie d'un $2$-tissu lin\'eaire de $\mathbb{C}^2$ que nous avons rencontr\'e dans l'Exemple \ref{firs}: le tissu des tangentes \`a la parabole $v=u^2/4$. Quant \`a l'ensemble $\Sigma$ c'est l'image de l'application $\widetilde{\xi}$, restriction de $\xi$ \`a $\mathbb{S}^1\times \overline{\mathbb{D}(0,1)}$ o\`u $\mathbb{D}(0,1)$ est le disque unit\'e de $\mathbb{C}$. Remarquons que $\widetilde{\xi}$ est injective sur $\mathbb{S}^1\times\mathbb{D}(0,1)$ et $(2:1)$ sur le bord $\mathbb{S}^1\times\mathbb{S}^1$. Par suite $\Sigma$ s'identifie topologiquement \`a l'espace compact $\mathbb{S}^1\times\overline{\mathbb{D}(0,1)}$ o\`u l'on a identifi\'e les points $(\mathrm{e}^{\mathrm{i}\vartheta},\mathrm{e}^{\mathrm{i}\varphi})$ et $(\mathrm{e}^{\mathrm{i}\varphi},\mathrm{e}^{\mathrm{i}\vartheta})$. Son int\'erieur est un tore plein $\mathbb{S}^1\times\mathbb{D}(0,1)$ et son bord $\partial\Sigma$ le quotient du tore $\mathbb{T}^2=\mathbb{S}^1\times\mathbb{S}^1$ par l'involution $(s,t)\mapsto(t,s)$; c'est une bande de \textsc{M\"{o}bius}. Nous en d\'eduisons que le bord du bassin d'attraction $\partial\mathrm{W}_\mathrm{Id}^s({\bf 0})$ de ${\bf 0}$ de $\Phi_\mathrm{Id}$  est l'ensemble semi-alg\'ebrique $\mathrm{Inv}^{-1}(\Sigma)$. Il est d\'ecrit cette fois par les in\'equations $\mathbb{R}$-polynomiales sur $\mathcal{M}(2;\mathbb{C})$:
$$\left\{\begin{array}{l}
\Big(\vert\det\vert^2-\frac{1}{2}\vert \mathrm{tr}^2\vert+1\Big)^2-\frac{1}{4}\Big\vert\mathrm{tr}^2-4\det\Big\vert^2=0\\
\vert\det\vert^2-\frac{1}{2}\vert \mathrm{tr}^2\vert+1\geq 0\\
\vert\det\vert\leq 1
\end{array}\right.$$

\noindent Il est \textsc{Levi}-plat au sens o\`u il contient un ouvert dense (remplacer la derni\`ere in\'egalit\'e par $\vert\det\vert < 1$) \textsc{Levi}-plat: les vari\'et\'es complexes contenues dans $\partial\mathrm{W}_\mathrm{Id}^s({\bf 0})$ sont les images r\'eciproques par l'application alg\'ebrique~$\mathrm{Inv}$ des disques $\mathbb{D}(0,1)\ni u\mapsto (\mathrm{e}^{\mathrm{i}\vartheta}+u,\mathrm{e}^{\mathrm{i}\vartheta}u)$, $\vartheta$ \'etant fix\'e. C'est un ensemble \'evidemment non born\'e puisque $\mathrm{W}_\mathrm{Id}^s({\bf 0})$ ne l'est pas (il contient un voisinage des nilpotentes). Ceci va \`a l'encontre du cas des applications homog\`enes g\'en\'eriques. Rappelons \`a cet effet que lorsque $f\colon\mathbb{C}^N\to \mathbb{C}^N$ est une application polynomiale homog\`ene, le bassin d'attraction de $0$ est born\'e d\`es que $f^{-1}(0)$ se r\'eduit \`a $0$ ce qui est le cas g\'en\'eriquement~(\cite{BL}). En dimension $1$ pour une transformation homog\`ene $z\mapsto z^k$ le bord du bassin d'attraction de $0$ est \'evidemment le cercle unit\'e $\mathbb{S}^1$ et tous les points p\'eriodiques sont contenus dans ce cercle (hormis le point $0$); de plus ils y sont denses. D'autre part notons que les orbites born\'ees de $z^k$ sont contenues dans la fermeture de ce bassin d'attraction. Dans le contexte de $\Phi_\mathrm{Id}$ nous avons la: 

\begin{pro}
Les points p\'eriodiques de $\Phi_\mathrm{Id}$, except\'e ${\bf 0}$, sont contenus dans le bord $\partial\mathrm{W}_\mathrm{Id}^s({\bf 0})$ du bassin d'attraction de ${\bf 0}$.
\end{pro}

\begin{proof}
C'est une application directe de la description des points p\'eriodiques de $\Phi_{\mathrm{Id}\vert\mathcal{D}}$.
\end{proof}

\noindent Nous nous int\'eressons aussi \`a l'adh\'erence des points p\'eriodiques de $\Phi_\mathrm{Id}$ (except\'e ${\bf 0}$). Cet ensemble a plusieurs composantes. Rappelons les faits \'el\'ementaires suivants:
\begin{itemize}
\item $\mathrm{M}$ a une valeur propre de module $1$ et l'autre nulle si et seulement si $\mathrm{M}$ est dans l'ensemble $$\Lambda^0=\big\{\mathrm{M}\in\mathcal{M}(2;\mathbb{C})\,\big\vert\,\det\mathrm{M}=0,\,\vert\mathrm{tr}\mathrm{M}\vert=1\big\};$$

\item $\mathrm{M}$ a une valeur propre double de module $1$ si et seulement si $\mathrm{M}$ appartient \`a $$\Lambda^1=\big\{\mathrm{M}\in\mathcal{M}(2;\mathbb{C})\,\big\vert\,\mathrm{tr}\mathrm{M}^2-4\det\mathrm{M}=0,\,\vert\det\mathrm{M}\vert=1\big\};$$

\item enfin $\mathrm{M}$ a deux valeurs propres de module $1$ distinctes si et seulement si $\mathrm{M}$ est dans $$\Lambda^2=\left\{\mathrm{M}\in\mathcal{M}(2;\mathbb{C})\,\Big\vert\,\frac{(\mathrm{tr}\mathrm{M})^2}{\det\mathrm{M}}\in [0,4[,\,\vert\det\mathrm{M}\vert=1\right\}.$$
\end{itemize}

\noindent Pour l'application $z\mapsto z^2$ l'adh\'erence des points p\'eriodiques co\"{i}ncide avec le bord du bassin d'attraction de l'origine auquel on ajoute l'origine. L'\'enonc\'e qui suit d\'ecrit l'adh\'erence des points p\'eriodiques de $\Phi_\mathrm{Id}$ et montre via des arguments de dimension que cette adh\'erence est diff\'erente du bord du bassin d'attraction.

\begin{pro}
L'adh\'erence des points p\'eriodiques de $\Phi_\mathrm{Id}$ est l'union de $\big\{{\bf 0}\big\}\cup\Lambda^0\cup\Lambda^1\cup\Lambda^2$.
\end{pro}

\begin{proof}
Le seul point \`a noter est qu'un \'el\'ement de $\Lambda^1$ est conjugu\'e \`a une matrice $\left[\begin{array}{cc}\mathrm{e}^{\mathrm{i}\vartheta}&\varepsilon\\0&\mathrm{e}^{\mathrm{i}\vartheta}\end{array}\right],$ avec $\varepsilon\in\{0,1\}$; une telle matrice est limite de matrices diagonalisables $\left[\begin{array}{cc}x_{k_n} & \varepsilon\\ 0 & t_{k_n}\end{array}\right]$ qui sont p\'eriodiques pour $\Phi_\mathrm{Id}$.
\end{proof}

\noindent Cette adh\'erence est encore un semi-alg\'ebrique d\'ecrit par l'union disjointe $$\big\{{\bf 0}\big\}\sqcup\big\{\vert\mathrm{tr}\vert=1,\,\vert\det\mathrm{M}\vert=0\big\}\sqcup\left\{\frac{(\mathrm{tr}\mathrm{M})^2}{\det\mathrm{M}}\in [0,4],\,\big\vert\det\mathrm{M}\vert=1\right\}.$$

\begin{pro}
Les orbites born\'ees de $\Phi_\mathrm{Id}$ sont contenues dans $\overline{\mathrm{W}_\mathrm{Id}^s({\bf 0})}$; tout \'el\'ement $\mathrm{M}$ de $\overline{\mathrm{W}_\mathrm{Id}^s({\bf 0})}$ a son orbite born\'ee \`a l'exception des matrices $\mathrm{M}$ ayant pour type de \textsc{Jordan} $\left[\begin{array}{cc} \lambda & 1 \\0& \lambda \end{array}\right]$ avec $\vert\lambda\vert=1$.
\end{pro}

\begin{rem}
L'image r\'eciproque par $\mathrm{Inv}$ du ruban de \textsc{M\"{o}bius} $\partial\Sigma$ est exactement l'union $$\Lambda^1 \cup\Lambda^2=\left\{\frac{(\mathrm{tr}\mathrm{M})^2}{\det\mathrm{M}}\in [0,4],\,\big\vert\det\mathrm{M}\vert=1\right\},$$ l'image r\'eciproque du bord $\partial(\partial \Sigma)$ de cette bande de \textsc{M\"{o}bius} \'etant pr\'ecis\'ement l'ensemble $\Lambda^1$.
\end{rem}

\begin{rem}
Notons que $\Sigma$ et $\partial\Sigma$ sont invariants par $\mathrm{Sq}\,\Phi_\mathrm{Id}$ et sont born\'es (ceci r\'esulte des pr\'esentations \emph{(i)} et \emph{(ii)} de la Proposition \ref{mod1}).
\end{rem}

\begin{prob}
Il serait int\'eressant de d\'ecrire les automorphismes holomorphes de l'ouvert $\mathrm{W}_\mathrm{Id}^s({\bf 0})$. Parmi ceux-ci on trouve les transformations $\sigma_\mathrm{P}\colon\mathrm{M}\mapsto\mathrm{P}\mathrm{M}\mathrm{P}^{-1}$ et les transformations lin\'eaires $\Theta\colon\mathrm{M}\mapsto\mathrm{e}^{\mathrm{i}\vartheta}\mathrm{M}$; dit autrement nous avons une action de $\mathbb{S}^1\times \mathrm{PGL}(2;\mathbb{C})$. Peut-\^etre serait-il plus facile de d\'ecrire les automorphismes holomorphes de $\mathcal{M}(2;\mathbb{C})\simeq\mathbb{C}^4$ qui pr\'eservent $\mathrm{W}_\mathrm{Id}^s({\bf 0})$ ? De tels automorphismes vont en effet pr\'eserver $\partial \mathrm{W}_\mathrm{Id}^s({\bf 0})$ et sa structure \textsc{Levi}-plate et donc le tissu image r\'eciproque par $\mathrm{Inv}$ du tissu des tangentes \`a la parabole $v=u^2/4$.
\end{prob}

\subsection{Quelques compl\'ements}

\noindent Les dynamiques de $\widetilde{\Phi_\mathrm{Id}}$ et $\mathbb{P}(\Phi_\mathrm{Id})$ se d\'eduisent plus ou moins facilement de celle de $\Phi_\mathrm{Id}$. La description de $\widetilde{\Phi_\mathrm{Id}}$ se fait d\`es que l'on connait celle de $\Phi_\mathrm{Id}\colon\mathcal{M}(2;\mathbb{C})\to\mathcal{M}(2;\mathbb{C})$ plus celle sur l'hyperplan \`a l'infini qui est un $\mathbb{P}^3(\mathbb{C})$. Mais sur l'hyperplan \`a l'infini qui est compl\`etement invariant la dynamique est exactement celle de $$\mathbb{P}(\Phi_\mathrm{Id})\colon\mathbb{P}(\mathcal{M}(2;\mathbb{C}))\simeq\mathbb{P}^3(\mathbb{C})\dashrightarrow
\mathbb{P}(\mathcal{M}(2;\mathbb{C}))\simeq\mathbb{P}^3(\mathbb{C}).$$ 

\noindent L'application $\mathbb{P}(\Phi_\mathrm{Id})$ est d\'efinie en coordonn\'ees homog\`enes par
$$\mathbb{P}(\Phi_\mathrm{Id})\colon(x:y:z:t)\mapsto(x^2+yz:y(x+t):z(x+t):yz+t^2).$$ R\'ecapitulons ses propri\'et\'es 
\begin{itemize}
\item le lieu d'ind\'etermination de $\mathbb{P}(\Phi_\mathrm{Id})$ est le projectivis\'e de l'ensemble $\mathcal{N}(2;\mathbb{C})$ des matrices nilpotentes;

\item le projectivis\'e des matrices de trace nulle (priv\'e de $\mathbb{P}(\mathcal{N}(2; \mathbb{C}))$) est contract\'e sur le point $(1:0:0:1)$ correspondant \`a la matrice identit\'e et ce point est fixe;

\item les autres points fixes sont le projectivis\'e de la quadrique $\det=0,$ priv\'e de $\mathbb{P}(\mathcal{N}(2; \mathbb{C}))$;

\item la fibre g\'en\'erique de $\mathbb{P}(\Phi_\mathrm{Id})$ a deux \'el\'ements;

\item la fibration en $2$-plans $\mathcal{P}$ qui est invariante fibre \`a fibre est maintenant la fibration en $\mathbb{P}^1(\mathbb{C})$ radiale en $(1:0:0:1)$ dans $\mathbb{P}^3(\mathbb{C})$;

\item la droite sp\'eciale $(\lambda,\lambda,0,\lambda),$ $\lambda\not=0,$ repr\'esente \`a conjugaison pr\`es les matrices non diagonalisables (inversibles) et donne un seul point bien s\^ur dans $\mathbb{P}^3(\mathbb{C}),$ la \og matrice unipotente\fg\, $(1:1:0:1)$. Nous obtenons l'orbite sp\'eciale $(1/2^n:1:0:1/2^n)$. C'est une orbite discr\`ete qui converge vers $(0:1:0:0)$ qui repr\'esente la matrice nilpotente standard. Autrement dit les it\'er\'es des matrices unipotentes convergent vers les matrices nilpotentes qui sont d'ind\'etermination;

\item la dynamique hormis les probl\`emes de type points d'ind\'etermination et ensembles contract\'es se comprend aussi essentiellement via la dynamique de $\chi_2\colon\mathbb{P}^1(\mathbb{C})\to\mathbb{P}^1( \mathbb{C})$, $z\mapsto z^2$. Par exemple dans le $\mathbb{P}^1(\mathbb{C})$ des matrices diagonales $(\lambda_1:0:0:\lambda_2)$ l'application s'\'ecrit $(\lambda_1^2:0:0:\lambda_2^2)$ c'est-\`a-dire $u\mapsto u^2$. 
\end{itemize}

\bigskip

\noindent Comme pour $\Phi_\mathrm{Id}$ on sait d\'ecrire la dynamique des polyn\^omes de matrices; par exemple celle de $\varsigma_c$ (points p\'eriodiques, ensembles invariants) s'obtient essentiellement \`a partir de $\varsigma_{c_{\vert\mathcal{D}}}$ en faisant agir l'action adjointe; celle de $\varsigma_{c_{\vert\mathcal{D}}}$ se d\'eduit quant \`a elle directement de celle de $z^2+c$. Si on note $\mathcal{J}_c\subset\mathbb{C}$ l'ensemble de \textsc{Julia} de $z^2+c$, l'ensemble $\mathcal{J}_c\times\mathcal{J}_c\subset\mathcal{D}$ est invariant par $\varsigma_c$. La fermeture de son satur\'e $\overline{\mathcal{O}( \mathcal{J}_c\times\mathcal{J}_c;\varsigma_c)}$ par l'action adjointe est un ferm\'e invariant par $\varsigma_c$ dans lequel on a densit\'e des points p\'eriodiques. Il y a d'autres ferm\'es invariants, comme par exemple le satur\'e ferm\'e de $\mathcal{J}_c\times\mathbb{C} \subset\mathcal{D},$ toujours par l'action adjointe. De m\^eme que si~$q$ est un point fixe de $z^2+c$, alors le satur\'e ferm\'e de $\mathcal{J}_c\times\big\{q\big\}$ est encore invariant. La description des orbites born\'ees est relativement raisonnable (mais zoologique), en liaison avec celles de $z^2+c$.

\subsection{Centralisateurs}

\noindent Une fa\c{c}on de mesurer la complexit\'e d'une transformation est d'examiner son centralisateur, \emph{i.e.} son groupe des commutateurs.  La transformation $\chi_2$, de m\^eme que les $z\mapsto z^k,$ a un r\^ole sp\'ecial dans la dynamique \`a une variable; comme les polyn\^omes de \textsc{Tchebychev} elle a son ensemble de \textsc{Julia} lisse et son centralisateur n'est pas r\'eduit \`a ses propres it\'er\'es. Il est donc naturel d'examiner \og certains centralisateurs\fg\, de $\Phi_\mathrm{Id}$.

\noindent Soit $f\colon X\dashrightarrow X$ une transformation rationnelle; introduisons les groupes suivants
\begin{align*}
&\mathrm{Bir}(X;f)=\big\{g\colon X\dashrightarrow X\text{ transformation birationnelle }\,\big\vert\, g\circ f=f\circ g\big\},\\
&\mathrm{Aut}(X;f)=\big\{g\colon X\to X\text{ automorphisme holomorphe }\,\big\vert\, g\circ f=f\circ g\big\}.
\end{align*}

\subsubsection{Caract\'erisation de $\mathrm{Aut}(\mathcal{M}(2;\mathbb{C});\Phi_\mathrm{Id})$}

\noindent Les transformations $\sigma_\mathrm{P}$ d\'efinies par $\sigma_\mathrm{P}(\mathrm{M})=\mathrm{P}\mathrm{M}\mathrm{P}^{-1}$ sont dans $\mathrm{Aut}(\mathcal{M}(2;\mathbb{C});\Phi_\mathrm{Id})$ et sont compatibles \`a la conjugaison. Notons que si $\Phi$ est compatible \`a la conjugaison, alors par d\'efinition les~$\sigma_{\mathrm{P}}$ sont dans $\mathrm{Aut}(\mathcal{M}(2;\mathbb{C});\Phi)$. L'application de transposition 
\begin{align*}
&\mathcal{T}\colon\mathcal{M}(2;\mathbb{C})\to\mathcal{M}(2;\mathbb{C}), &&\left[\begin{array}{cc} x & y\\ z & t\end{array}\right]\mapsto \left[\begin{array}{cc} x & z\\ y & t\end{array}\right]
\end{align*}
commute \`a $\Phi_\mathrm{Id}$; ainsi $\mathrm{Aut}(\mathcal{M}(2;\mathbb{C});\Phi_\mathrm{Id})$ contient des \'el\'ements non compatibles \`a la conjugaison bien que~$\Phi_\mathrm{Id}$ le soit. De m\^eme la transformation 
\begin{align*}
&\mathcal{I}\colon\mathcal{M}(2;\mathbb{C})\to\mathcal{M}(2;\mathbb{C}),&&\left[\begin{array}{cc} x & y\\ z & t\end{array}\right]\mapsto \left[\begin{array}{cc} t & y\\ z & x\end{array}\right]
\end{align*}

\noindent commute \`a $\Phi_\mathrm{Id}$ et n'est pas compatible \`a la conjugaison; notons que $\mathcal{I}(\mathrm{M})=\transp\, (\sigma_\mathrm{P}(\mathrm{M}))$ o\`u $\mathrm{P}=\left[ \begin{array}{cc}0 & 1\\1 &0\end{array}\right]$. 

\begin{pro}\label{centautphi2}
Le groupe $\mathrm{Aut}(\mathcal{M}(2;\mathbb{C});\Phi_\mathrm{Id})$ est engendr\'e par la transposition $\mathcal{T}$ et les $\sigma_\mathrm{P};$ plus pr\'ecis\'ement $$\mathrm{Aut}(\mathcal{M}(2;\mathbb{C});\Phi_\mathrm{Id})\simeq\mathrm{PGL}(2;\mathbb{C})\rtimes\mathbb{Z} /2\mathbb{Z}.$$
\end{pro}

\begin{proof}[D\'emonstration]
Soit $\varphi$ un automorphisme de $\mathcal{M}(2;\mathbb{C})$ qui commute \`a $\Phi_\mathrm{Id}$. Puisque ${\bf 0}$ est le seul point fixe de $\Phi_\mathrm{Id}$ o\`u la diff\'erentielle de $\Phi_\mathrm{Id}$ est identiquement nulle, ${\bf 0}$ est fixe par $\varphi$. On \'ecrit $\varphi=\varphi_1+\varphi_2+\ldots$ le d\'eveloppement de \textsc{Taylor} de $\varphi,$ chaque $\varphi_k$ \'etant holomorphe de degr\'e $k$. Notons que la partie lin\'eaire~$\varphi_1$ de~$\varphi$ est dans $\mathrm{Aut}(\mathcal{M}(2;\mathbb{C});\Phi_\mathrm{Id})$ de sorte que $\Psi=\varphi\circ\varphi_1^{-1}$ aussi. Supposons que $\Psi$ ne soit pas trivial; \'ecrivons~$\Psi$ sous la forme $\Psi=\mathrm{Id}+\Psi_k+\Psi_{k+1}+\ldots$ avec $\Psi_k$ homog\`ene de degr\'e $k>1$ et $\Psi_k$ non nul. La commutation de $\Psi$ et $\Phi_\mathrm{Id}$ implique l'\'egalit\'e suivante
$$\mathrm{M}^2+\Psi_k(\mathrm{M}^2)+\Psi_{k+1}(\mathrm{M}^2)+\ldots=\Big(\mathrm{M}+\Psi_k(\mathrm{M})+\Psi_{k+1}(\mathrm{M})+\ldots\Big)^2= \mathrm{M}^2+\mathrm{M}\Psi_k(\mathrm{M})+\Psi_k(\mathrm{M})\mathrm{M}+\ldots$$ Pour des raisons de degr\'e nous avons pour tout $\mathrm{M}$ l'\'egalit\'e $\mathrm{M}\Psi_k(\mathrm{M})+\Psi_k(\mathrm{M})\mathrm{M}=0$. \'Ecrivons $\mathrm{M}$ (resp. $\Psi_k(\mathrm{M})$) sous la forme $\left[\begin{array}{cc} x &y\\ z & t\end{array}\right]$ (resp. $\left[\begin{array}{cc} X(\mathrm{M}) & Y(\mathrm{M})\\ Z(\mathrm{M}) & T(\mathrm{M})\end{array}\right]$); alors $\mathrm{M}\Psi_k(\mathrm{M})+\Psi_k(\mathrm{M})\mathrm{M}=~0$ se r\'e\'ecrit
\begin{equation}\label{syslin}
\left\{\begin{array}{llll}
2xX+yZ+zY=0\\
xY+yT+yX+tY=0\\
zX+tZ+xZ+zT=0\\
zY+yZ+2tT=0
\end{array}
\right.
\end{equation}

\noindent On peut voir (\ref{syslin}) comme un syst\`eme lin\'eaire, les inconnues \'etant les $X,$ $Y,$ $Z,$ $T,$ dont on notera $\Delta=\Delta(\mathrm{M})$ le d\'eterminant. Nous constatons que $\Delta(\mathrm{Id})=16;$  par suite $\Delta$ est non identiquement nul. En r\'esulte que $X,$ $Y,$ $Z,$~$T,$ et donc $\Psi_k,$ sont identiquement nuls. Nous en d\'eduisons que $\varphi=\varphi_1$ est n\'ecessairement un automorphisme lin\'eaire.

\noindent 

\noindent Remarquons que $\mathrm{Id}$ est l'unique point fixe de $\Phi_\mathrm{Id}$ en lequel la diff\'erentielle de $\Phi_\mathrm{Id}$ est $2\mathrm{Id}$. Il s'en suit que $\mathrm{Id}$ est fix\'e par $\varphi_1$ et la droite $\mathcal{C}=\big\{\lambda\mathrm{Id}\,\big\vert\,\lambda\in\mathbb{C}\big\}$ est invariante point par point par $\varphi_1$. En particulier $\varphi_1$ envoie le plan $\langle\lambda\mathrm{Id}+\mu\mathrm{M}\rangle$ dans le plan $\langle\lambda\mathrm{Id}+\mu\varphi(\mathrm{M})\rangle$. Ceci implique que la fibration $\mathcal{P}$ est invariante par $\varphi_1$. Pour $\mathrm{M}$ g\'en\'erique, $\mathrm{M}$ et $\varphi_1(\mathrm{M})$ sont diagonalisables: $\mathrm{M}=\mathrm{Q}\mathrm{D}\mathrm{Q}^{-1},$ $\varphi_1(\mathrm{M})=\mathrm{Q}'\mathrm{D}' \mathrm{Q}'^{-1}$. En composant $\varphi_1$ \`a gauche et \`a droite par des $\sigma_\mathrm{P}$ ad-hoc nous nous ramenons donc au cas o\`u $\varphi_1$ respecte la fibration $\mathcal{P}$ et le groupe $\mathcal{D}$. Si on \'ecrit~$\varphi(\mathrm{M})$ sous la forme  $\left[\begin{array}{cc} X(\mathrm{M}) & Y(\mathrm{M})\\ Z(\mathrm{M}) & T(\mathrm{M})\end{array}\right]$ o\`u cette fois les transformations $X,$ $Y,$ $Z$ et $T$ sont lin\'eaires on~a: 
\begin{align*}
& Y(x,0,0,t)=Z(x,0,0,t)=0, && X(x^2,0,0,t^2)=X^2(x,0,0,t), && T(x^2,0,0,t^2)=T^2(x,0,0,t).
\end{align*}
 
\noindent Puisque $\varphi_1$ est un automorphisme, ceci entra\^ine par un calcul direct l'alternative suivante:
\begin{align*}
&\text{ ou bien $X(x,0,0,t)=x$ et $T(x,0,0,t)=t$;} && \text{ ou bien $X(x,0,0,t)=t$ et $T(x,0,0,t)=x$.}
\end{align*} 
 
\noindent Quitte \`a composer $\varphi_1$ par l'involution $\mathcal{I},$ nous pouvons supposer que nous sommes dans la premi\`ere situation. Finalement $$\varphi_1=(x+\ell_1(y,z),\ell_2(y,z),\ell_3(y,z),t+\ell_4(y,z))$$ les $\ell_i$ \'etant lin\'eaires en $y$, $z$.  En r\'e\'ecrivant que $\Phi_\mathrm{Id}$ et $\varphi_1$ commutent on constate ais\'ement que $\ell_1=\ell_4=0$ et~$\ell_2\ell_3=yz$. Ainsi apr\`es ces modifications $\varphi_1$ est de l'un des deux types suivants
\begin{align*}
&g=\left(x,\rho y,\frac{z}{\rho},t\right), && h=\left(x,\rho z,\frac{y}{\rho},t\right).
\end{align*}

\noindent Remarquons que si $\mathrm{P}=\left[\begin{array}{cc} \sqrt{\rho} & 0\\ 0 & \rho\end{array}\right],$ alors $g=\sigma_\mathrm{P}$ et que, modulo $\sigma_\mathrm{P},$ l'application $h$ est la transposition.
\end{proof}

\subsubsection{Description de $\mathrm{Bir}(\mathbb{P}\mathcal{M}(2;\mathbb{C});\mathbb{P}(\Phi_\mathrm{Id}))$}

\noindent Dans la carte $t=1$ on remarque que $\mathbb{P}(\Phi_\mathrm{Id})$ s'\'ecrit de la fa\c{c}on suivante $$\left(\frac{x^2+yz}{1+yz},\frac{y(x+1)}{1+yz},\frac{z(x+1)}{1+yz}\right).$$ La fibration en $2$-plans $\mathcal{P},$ invariante par $\Phi_\mathrm{Id}$, initialement donn\'ee par les niveaux de $\left(\frac{y}{z},\frac{x-t}{z}\right),$ produit maintenant la fibration en droites donn\'ee cette fois par les fibres de $\left(\frac{y}{z},\frac{x-1}{z}\right)$. Nous allons trivialiser cette fibration en conjuguant~$\mathbb{P}(\Phi_\mathrm{Id})$ par la transformation birationnelle $f\colon (x,y,z)\mapsto(xz+1,yz,z)$. On v\'erifie que $$G=f^{-1}\mathbb{P}(\Phi_\mathrm{Id})f=\left(x,y,\frac{z(2+xz)}{1+yz^2}\right).$$ Le fait que $\mathcal{P}$ soit l'unique fibration en surfaces invariante fibre \`a fibre par $\Phi_\mathrm{Id}$ implique que la fibration en droites $(x,y)=$ cte est l'unique fibration en courbes pr\'eserv\'ee par $G$. Il en r\'esulte que si $F$ est une application birationnelle commutant \`a~$G$, alors $F$ est n\'ecessairement du type 
\begin{align*}
&\left(X(x,y),Y(x,y),\frac{a(x,y)z+b(x,y)}{c(x,y)z+d(x,y)}\right), && a,\,b,\,c,\,d\in\mathbb{C}\{x,y\},\, ad-bc\not\equiv 0.
\end{align*}
Visiblement la transformation $(x,y)\dashrightarrow (X(x,y),Y(x,y))$ doit \^etre birationnelle. Si l'on consid\`ere $F$ et $G$ comme des applications rationnelles sur $\mathbb{C}^2\times\mathbb{P}^1$ on constate qu'elles sont fibr\'ees et holomorphes en restriction aux fibres g\'en\'eriques. Les points fixes de $G$ sont les trois surfaces irr\'eductibles $z=0,$ $z=\infty$ et $yz^2-xz-1=0$. Si $F$ commute \`a $G$, alors $F$ pr\'eserve dans leur ensemble les points fixes de $G$. Plus pr\'ecis\'ement $F$ laisse invariante la surface $yz^2-xz-1=0$ et fixe ou permute les plans $z=0$ et $z=\infty$. 

\noindent$\mathfrak{i.}$ Supposons, dans un premier temps, que $yz^2-xz-1=0$, $z=0$ et $z=\infty$ soient fix\'es par $F$. Alors $F$ est du type suivant $$F(x,y,z)=(X(x,y),Y(x,y),a(x,y)z).$$ Maintenant $(yz^2-xz-1)\circ F=Ya^2z^2-Xaz-1$ est donc un multiple (comme \'el\'ement de $\mathbb{C}(x,y)[z]$) de $yz^2-xz-1$, de sorte que nous pouvons pr\'eciser $F$: $$F(x,y,z)=\left(\frac{x}{a(x,y)},\frac{y}{a^2(x,y)}, a(x,y)z\right).$$ On peut d'ailleurs remarquer que toute transformation rationnelle comme ci dessus, \emph{i.e.} pour tout choix de~$a$, commute \`a $G$. Le fait que $F$ soit birationnelle impose des restrictions sur $a$. En effet si $F$ est birationnelle alors $H\colon(x,y)\mapsto\left( \frac{x}{a(x,y)},\frac{y}{a^2(x,y)}\right)$ l'est aussi. En conjuguant $H$ par $(x,x^2y)$ on constate que $\left(\frac{x}{a(x,x^2y)},y\right)$ est birationnelle ce qui force $\frac{x}{a(x,x^2y)}$ \`a \^etre une transformation de \textsc{M\"{o}bius} en $x$ \`a param\`etre $y$. Nous en d\'eduisons que 
\begin{equation}\label{train}
a(x,y)=x\frac{\alpha\left(\frac{y}{x^2}\right)x+\beta\left(\frac{y}{x^2}\right)}{ \gamma\left(\frac{y}{x^2}\right)x+\delta\left(\frac{y}{x^2}\right)}
\end{equation}
o\`u les $\alpha,$ $\beta,$ $\gamma,$ $\delta$ sont dans $\mathbb{C}(t)$ et $\alpha\delta -\beta\gamma\not=0$. Le groupe des transformations birationnelles $$\left\{\left(\frac{x}{a(x,y)},\frac{y}{a^2(x,y)},a(x,y)z\right)\,\big\vert\, \text{ avec $a$ du type } (\ref{train})\right\}$$ commute \`a $G$ et est isomorphe \`a $\mathrm{PGL}(2;\mathbb{C}(t))$.

\noindent$\mathfrak{ii.}$ Maintenant consid\'erons l'\'eventualit\'e o\`u $F$ fixe $yz^2-xz-1$ et permute les plans $z=0$ et $z=\infty$. Alors $F$ est du type $\left(X(x,y),Y(x,y),\frac{a(x,y)}{z}\right)$. Comme $(yz^2-xz-1)\circ F=\frac{Ya^2-Xaz-z^2}{z^2}$ nous obtenons par le m\^eme argument que pr\'ec\'edemment $$F(x,y,z)=\left(-\frac{x}{a(x,y)y},\frac{1}{a^2(x,y)y},\frac{a(x,y)}{z}\right);$$ un tel $F$ ne commute jamais \`a $G$. Nous pouvons donc \'enoncer la:

\begin{pro}
On a: $$\mathrm{Bir}(\mathbb{P}\mathcal{M}(2;\mathbb{C});\mathbb{P}(\Phi_\mathrm{Id}))=\left\{f\left(\frac{x}{a(x,y)},\frac{y}{a^2(x,y)},a(x,y)z\right)f^{-1}\,\Big\vert\, \text{ avec $a$ du type } (\ref{train}), \,f=(xz+1,yz,z)\right\};$$ en particulier $\mathrm{Bir}(\mathbb{P}\mathcal{M}(2;\mathbb{C});\mathbb{P}(\Phi_\mathrm{Id}))$ est isomorphe \`a $\mathrm{PGL}(2; \mathbb{C}(t))$.
\end{pro}

\subsubsection{Le groupe $\mathrm{Bir}(\mathcal{M}(2;\mathbb{C}); \Phi_\mathrm{Id})$}

\noindent Donnons un exemple de transformation appartenant \`a $\mathrm{Bir}(\mathcal{M}(2;\mathbb{C}); \Phi_\mathrm{Id})$. L'application 
\begin{align*}
&\iota\colon\mathcal{M}(2;\mathbb{C})\to\mathcal{M}(2;\mathbb{C}), && \mathrm{M}\mapsto \mathrm{M}^{-1}
\end{align*}

\noindent est dans $\mathrm{Bir}(\mathcal{M}(2;\mathbb{C}); \Phi_\mathrm{Id});$ elle est compatible \`a la conjugaison. Remarquons qu'elle commute aux applications $\sigma_\mathrm{P}$ ainsi qu'\`a $\mathcal{T}$. On constate que la restriction de $\iota$ aux matrices diagonales s'identifie \`a l'involution de \textsc{Cremona} standard en dimension deux $(x,t)\to\left( \frac{1}{x},\frac{1}{t}\right)$.

\noindent Soit $f$ un \'el\'ement de $\mathrm{Bir}(\mathcal{M}(2;\mathbb{C});\Phi_\mathrm{Id})$. Comme $f(\mathrm{M}^2)=f(\mathrm{M})f(\mathrm{M})$ nous avons $\det f(\mathrm{M}^2)=\Big(\det f(\mathrm{M}) \Big)^2$. Ce type de propri\'et\'e est d\'ecrit dans le lemme qui suit.

\begin{lem}\label{carre}
Soit $\psi\colon\mathcal{M}(2;\mathbb{C})\dashrightarrow\mathbb{P}^1(\mathbb{C})$ une fonction rationnelle satisfaisant l'\'equation fonctionnelle $\psi(\mathrm{M}^2)=\psi(\mathrm{M})^2$ pour tout $\mathrm{M}$ dans $\mathcal{M}(2;\mathbb{C})$. 

\noindent Alors $\psi$ est de la forme $\mathrm{M}\dashrightarrow \varepsilon(\det\mathrm{M})^k$ o\`u $k$ d\'esigne un entier relatif et $\varepsilon$ un \'el\'ement de $\{0,1\}$.
\end{lem}

\begin{proof}[D\'emonstration]
Quitte \`a changer $\psi$ en $1/\psi$ nous pouvons supposer que l'hyperplan $\mathfrak{sl}(2; \mathbb{C})$ n'est pas contenu dans les p\^oles de $\psi$. Ainsi $\psi$ est rationnelle en restriction \`a $\mathfrak{sl}(2;\mathbb{C})$ et pour $x$, $y$, $z$ g\'en\'eriques $\psi(x,y,z,-x)$ est bien d\'efini. Par suite $\psi((x^2+yz)\mathrm{Id})$ est rationnel, bien d\'efini et donc $\delta\colon t\mapsto\psi(t\mathrm{Id})$ aussi. L'hypoth\`ese $\psi(\mathrm{M}^2)=\psi(\mathrm{M})^2$ implique que $\delta(t^2)=\big(\delta(t)\big)^2$, \emph{i.e.} l'application $\delta$ commute \`a $s\mapsto s^2$; par cons\'equent $\delta(t)=~\varepsilon' t^{k'}$ pour un certain $k'$ dans $\mathbb{Z}$ et pour $\varepsilon'$ dans $\{0,1\}$. Autrement dit pour tout $\mathrm{M}$ dans $\mathfrak{sl}(2;\mathbb{C})$ nous avons $\big(\psi(\mathrm{M})\big)^2=\psi(\mathrm{M}^2)=\psi(\det(\mathrm{M})\mathrm{Id})=\varepsilon' (\det\mathrm{M})^{k'}$. La restriction de $\Psi$ \`a $\mathfrak{sl}(2;\mathbb{C})$ est donc donn\'ee par $\psi(\mathrm{M})=\varepsilon (\det\mathrm{M})^k$ pour un certain $k$ dans $\mathbb{Z}$ et pour $\varepsilon$ dans $\{0,1\}$.

\noindent Soit $\mathrm{M}$ dans $\mathcal{M}(2;\mathbb{C})$ tel que $\mathrm{M}^2$ appartienne \`a $\mathfrak{sl}(2;\mathbb{C})$. D'une part $\psi(\mathrm{M}^2)=\varepsilon(\det\mathrm{M})^{2k}$ et d'autre part $\psi(\mathrm{M}^2)=\big(\psi( \mathrm{M})\big)^2$; nous en d\'eduisons que $\psi(\mathrm{M})$ s'\'ecrit aussi $\varepsilon(\det\mathrm{M})^k$. On peut r\'ep\'eter cet argument avec tous les \'el\'ements $\mathrm{M}$ de $\mathcal{M}(2;\mathbb{C})$ pour lesquels $\mathrm{M}^{2^n}$ appartient \`a $\mathfrak{sl}(2;\mathbb{C})$ pour un certain $n$. L'ensemble de ces matrices est \textsc{Zariski} dense d'o\`u l'\'enonc\'e.
\end{proof}

\noindent Soit $f$ dans $\mathrm{Bir}(\mathcal{M}(2;\mathbb{C});\Phi_\mathrm{Id});$ on peut appliquer le Lemme \ref{carre} \`a $\psi=\det f$: il existe $k$ tel que $\det f(\mathrm{M})= \varepsilon (\det\mathrm{M})^k$, $\varepsilon\in\{0,\,1\}$. Comme $f$ est birationnelle $\varepsilon$ ne peut \^etre nul et la fibration $\det=$ cte est invariante sous l'action de $f$. N\'ecessairement ceci implique que $k=\pm 1$ et, quitte \`a composer $f$ par $\iota\colon\mathrm{M}\to\mathrm{M}^{-1}$, nous supposerons dans la suite que $k=1$, \emph{i.e.} $\det f(\mathrm{M})=\det\mathrm{M}$.

\noindent Supposons que~$f\in\mathrm{Bir}(\mathcal{M}(2;\mathbb{C});\Phi_\mathrm{Id})$ ait un p\^ole, \emph{i.e.} $f$ s'\'ecrit $\frac{g}{h}$ avec $g$, $h$ polynomiales et $h$ non constante. Nous pouvons bien \'evidemment nous ramener au cas o\`u $g$ et $h$ n'ont pas de facteur commun.

\begin{lem}\label{hypersurface}
L'hypersurface $h=0$ est compl\`etement invariante.
\end{lem}

\begin{proof}[D\'emonstration]
Pla\c{c}ons-nous dans $\mathbb{P}^4(\mathbb{C});$ notons $\widetilde{\Phi_\mathrm{Id}}$ (resp. $\widetilde{f}$) le prolongement de $\Phi_\mathrm{Id}$ (resp. $f$) \`a $\mathbb{P}^4(\mathbb{C})$. D\'esignons l'hyperplan \`a l'infini par $\mathcal{H}_\infty$. Puisque $\widetilde{\Phi_\mathrm{Id}}(x:y:z:t:s)=(x^2+yz:y(x+t): z(x+t):t^2+yz:s^2)$ l'hyperplan $\mathcal{H}_\infty$ est compl\`etement invariant par $\widetilde{\Phi_\mathrm{Id}}$. 

\noindent Soit $\mathrm{M}$ un point g\'en\'erique de $h=0$. D'une part $\mathrm{M}$ n'est pas d'ind\'etermination pour $\widetilde{f}$, d'autre part $\widetilde{f}(\mathrm{M})$ appartient \`a~$\mathcal{H}_\infty$; la compl\`ete invariance de $\mathcal{H}_\infty$ assure que $\widetilde{\Phi_\mathrm{Id}}(\widetilde{f}(\mathrm{M}))$ aussi. Comme $\widetilde{f}$ et $\widetilde{\Phi_\mathrm{Id}}$ commutent,~$\widetilde{f}(\widetilde{\Phi_\mathrm{Id}}(\mathrm{M}))$ est sur $\mathcal{H}_\infty$; il s'en suit que $\widetilde{\Phi_\mathrm{Id}}(\mathrm{M})$ appartient \`a $h=0$. 

\noindent Par ailleurs soit $\mathrm{M}$ tel que $\widetilde{\Phi_\mathrm{Id}}( \mathrm{M})$ soit contenu dans $h=0,$ $\mathrm{M}$ g\'en\'erique pour cette propri\'et\'e; alors~$\widetilde{f}(\widetilde{\Phi_\mathrm{Id}}(\mathrm{M}))$ appartient \`a $\mathcal{H}_\infty$ et, par commutation de $\widetilde{f}$ et $\widetilde{\Phi_\mathrm{Id}},$ $\widetilde{\Phi_\mathrm{Id}}(\widetilde{f}(\mathrm{M}))$ aussi. La compl\`ete invariance de $\mathcal{H}_\infty$ par $\widetilde{\Phi_\mathrm{Id}}$ assure que $\widetilde{f}(\mathrm{M})$ appartient \`a $\mathcal{H}_\infty$. 
\end{proof}

\noindent Le Lemme \ref{hypersurface} assure que l'hypersurface $h=0$ est compl\`etement invariante par $\Phi_\mathrm{Id}$. D'apr\`es \cite{CLN} il existe une constante non nulle $c$ telle que $h(\mathrm{M}^2 )=c\Big(h(\mathrm{M})\Big)^2$. Mais quitte \`a changer $h$ en $h/c$ nous pouvons supposer dans la suite que $h(\mathrm{M}^2)=h(\mathrm{M})^2$. En particulier il existe un certain entier $s\in\mathbb{N}$ tel que, pour tout~$\mathrm{M}$ dans $\mathcal{M}(2;\mathbb{C})$, nous ayons $h(\mathrm{M})=(\det\mathrm{M})^s$.

\noindent En consid\'erant l'adh\'erence de \textsc{Zariski} de l'orbite d'une matrice g\'en\'erique $\mathrm{M}$ on constate que $f$ laisse invariante la fibration en $2$-plans $\mathcal{P}$, \emph{i.e.} $\mathcal{P}(f(\mathrm{M}))=f(\mathcal{P}(\mathrm{M}))$. Quitte \`a changer $f$ en $\sigma_\mathrm{P}f\sigma_\mathrm{Q}$ avec $\mathrm{P}$, $\mathrm{Q}$ bien choisis nous pouvons supposer que $f$ est rationnelle et bien d\'efinie sur $\mathcal{D}$ et que $\mathcal{D}$ est invariant par $f$. Ceci signifie que la restriction de $f$ \`a $\mathcal{D}$ induit une application birationnelle $f_{\vert\mathcal{D}}\colon\mathcal{D}\dashrightarrow\mathcal{D}$. \'Evidemment $f_{\vert \mathcal{D}}$ commute \`a la restriction de $\Phi_\mathrm{Id}$ \`a $\mathcal{D}$. Le lemme qui suit se d\'emontre facilement.

\begin{lem}\label{comcar}
Soit $\eta\colon\mathbb{C}^2\dashrightarrow\mathbb{C}^2$ une transformation birationnelle commutant \`a $\xi\colon\mathbb{C}^2\to\mathbb{C}^2$, $(x,t)\mapsto(x^2,t^2)$. Alors~$\eta$ est de l'un des types suivants
\begin{align*}
& (x,t), &&(t,x), &&\left(\frac{1}{x},\frac{1}{t}\right), &&\left(\frac{1}{t},\frac{1}{x}\right), &&\left(x,\frac{1}{t}\right), &&\left(\frac{1}{x},t\right), &&\left(\frac{1}{t},x\right), &&\left(t,\frac{1}{x}\right).
\end{align*}

\noindent Autrement dit $\mathrm{Bir}(\mathbb{C}^2;\xi)$ est engendr\'e par les involutions $\left(\frac{1}{x},t\right)$ et $(t,x)$.
\end{lem}

\noindent L'\'enonc\'e \ref{comcar} s'applique \'evidemment \`a $f_{\vert\mathcal{D}}$ puisque $\Phi_{\mathrm{Id}_{\vert\mathcal{D}}}$ s'identifie \`a $(x,t)\mapsto(x^2,t^2)$. Mais nous avons suppos\'e que $\det f(\mathrm{M})=\det\mathrm{M},$ ce qui n'autorise  que les deux premiers mod\`eles du Lemme \ref{comcar}; nous en d\'eduisons, apr\`es les adaptations pr\'ec\'edentes, que l'application $f_{\vert\mathcal{D}}$ est soit l'identit\'e, soit $\mathrm{diag}(x,t)\mapsto\mathrm{diag}(t,x)$. Mais quitte \`a modifier encore $f$ en composant par $\mathcal{I}\colon\left[\begin{array}{cc} x & y \\ z & t \end{array}\right]\to\left[\begin{array}{cc} t & y \\ z & x \end{array}\right]$ nous nous ramenons \`a la premi\`ere \'eventualit\'e, \emph{i.e.} \`a $f_{\vert\mathcal{D}}=\mathrm{id}_{\vert \mathcal{D}}.$ R\'ecapitulons, nous pouvons supposer que
\begin{itemize}
\item $\det f(\mathrm{M})=\det\mathrm{M}$;

\item $f_{\vert\mathcal{D}}=\mathrm{id}_{\vert\mathcal{D}}$;

\item si $f$ a un p\^ole alors $f$ s'\'ecrit $f(\mathrm{M})=\frac{g(\mathrm{M})}{(\det \mathrm{M})^s}$.
\end{itemize}

\noindent On constate que $f$ est alors d\'efinie holomorphe au voisinage de l'identit\'e qui est un point fixe de $\Phi_\mathrm{Id}$ et de~$f$. La partie lin\'eaire de $\Phi_\mathrm{Id}$ au point fixe $\mathrm{Id}$ est la multiplication par $2$. En fait l'application exponentielle $\exp\colon\mathcal{M}(2; \mathbb{C})\to\mathcal{M}(2;\mathbb{C})$ lin\'earise l'application $\Phi_\mathrm{Id}$ $$\Phi_\mathrm{Id}(\exp \mathrm{M})=(\exp\mathrm{M})^2=\exp(2\mathrm{M}).$$ Maintenant les seules transformations holomorphes qui commutent \`a $\mathrm{M}\mapsto 2\mathrm{M}$ sont les transformations lin\'eaires. Ceci implique que si $f$ satisfait les propri\'et\'es qui pr\'ec\`edent, $\exp$ lin\'earise $f$ au voisinage de l'identit\'e 
\begin{equation}\label{expexp}
f(\exp\mathrm{M})=\exp f_1(\mathrm{M})
\end{equation}
\noindent o\`u $f_1=Df_{(\mathrm{Id})}$ est la partie lin\'eaire de $f$ en $\mathrm{Id}$. Notons que l'application $\mathrm{M}\mapsto f(\exp\mathrm{M})$ est partout d\'efinie holomorphe puisque $$f(\exp\mathrm{M})=\frac{g(\exp\mathrm{M})}{\Big(\det\exp \mathrm{M}\Big)^s}=\frac{g(\exp\mathrm{M})}{\Big(\exp\mathrm{tr}\,\mathrm{M}\Big)^s}.$$ Il s'en suit que $f_1$ d\'etermine $f$ et que l'identit\'e (\ref{expexp}) est en fait globale.

\noindent Nous allons maintenant d\'ecrire $f_1$ en utilisant les propri\'et\'es de $f.$ Puisque $f$ est l'identit\'e sur $\mathcal{D},$ nous avons $$f_1(x,y,z,t)=(x+a_1y+b_1z,a_2y+b_2z,a_3y+b_3z,t+ a_4y+b_4z).$$ Comme $\det$ est invariant par $f$ et que infinit\'esimalement le d\'eterminant en $\mathrm{Id}$ s'identifie \`a la trace (via $\exp$), l'application~$f_1$ pr\'eserve la trace c'est-\`a-dire $a_1+a_4=b_1+b_4=0$. Soient $\mathrm{P}$ une matrice inversible et $p,$ $q$ dans $\mathbb{Z}$; on~a $$\mathrm{Id}=\mathrm{P}\cdot\exp\left(\mathrm{diag}(2\mathrm{i}\pi p,2\mathrm{i}\pi q)\right)\cdot\mathrm{P}^{-1}=\exp\left(\mathrm{P}\cdot\mathrm{diag}(2\mathrm{i}\pi p,2\mathrm{i}\pi q)\cdot\mathrm{P}^{-1}\right)$$ et en utilisant (\ref{expexp}) $$\mathrm{Id}=f(\mathrm{Id})=f\left(\exp\left(\mathrm{P}\cdot\mathrm{diag}(2\mathrm{i}\pi p,2\mathrm{i}\pi q)\cdot\mathrm{P}^{-1}\right)\right)=\exp\left(f_1\left(\mathrm{P}\cdot\mathrm{diag}(2\mathrm{i}\pi p,2\mathrm{i}\pi q)\cdot\mathrm{P}^{-1}\right)\right).$$ En particulier $f_1\left(\mathrm{P} \cdot\mathrm{diag}(2\mathrm{i}\pi p,2\mathrm{i}\pi q)\cdot\mathrm{P}^{-1}\right)$ est un conjugu\'e d'un certain $\mathrm{diag}(2\mathrm{i}\pi p',2\mathrm{i}\pi q')$. Ceci signifie qu'il y a un ensemble d\'enombrable d'orbites de l'action adjointe qui sont envoy\'ees par $f_1$ sur d'autres orbites de cette m\^eme action. Par passage \`a l'adh\'erence de \textsc{Zariski} nous en d\'eduisons que $f_1$ envoie orbites (de l'action adjointe) dans orbites. En particulier la fibration par $\mathrm{Inv}=(\mathrm{tr},\det)$ est invariante par $f_1$. Soit $Q$ la forme quadratique $Q=\det f_1$: $$Q=(x+a_1y+b_1z)(t-a_1y-b_1z)-(a_2y+b_2z)(a_3y+b_3z).$$ Puisque $Q$ doit \^etre constant sur les niveaux de $\mathrm{Inv}$ nous obtenons en utilisant un argument de \cite{Ma}
\begin{equation}\label{malg}
Q=\varphi(\mathrm{tr},\det)
\end{equation}
avec $\varphi$ un germe d'application holomorphe en $0$, l'\'egalit\'e (\ref{malg}) \'etant comprise au voisinage de $0$. Mais en d\'eveloppant~$\varphi$ en s\'erie et en utilisant l'expression de $Q$ nous constatons que $Q=\det$ et par suite $a_1=b_1=~0$ et $\left[\begin{array}{cc} a_2 & b_2 \\ a_3 & b_3\end{array}\right]$ est l'une des matrices $\mathrm{Id}$ ou $\left[ \begin{array}{cc} 0 & 1 \\ 1 & 0\end{array}\right]$. Dans le premier cas $f_1=\mathrm{Id}$ et donc~$f(\mathrm{M})=~\mathrm{M}$; le second cas montre que $f_1$ est la diff\'erentielle de l'application transposition, et par cons\'equent $f$ co\"{i}ncide elle-m\^eme avec cette transposition (nous utilisons la commutation de $\Phi_\mathrm{Id}$ avec $\mathcal{T}$ et le fait que la partie lin\'eaire $f_1$ d\'etermine $f$). Finalement nous avons la

\begin{pro}
Le groupe $\mathrm{Bir}(\mathcal{M}(2;\mathbb{C});\Phi_\mathrm{Id})$ est engendr\'e par $\mathrm{Aut} (\mathcal{M}(2;\mathbb{C});\Phi_\mathrm{Id})\simeq \mathrm{PGL}(2;\mathbb{C})\rtimes\mathbb{Z}/2\mathbb{Z}$ et par l'involution $\iota\colon\mathrm{M}\mapsto\mathrm{M}^{-1}$.
\end{pro}

\begin{rem}
D'un point de vue abstrait $\mathrm{Bir}(\mathcal{M}(2;\mathbb{C});\Phi_\mathrm{Id})$ s'identifie \`a $(\mathrm{PGL}(2;\mathbb{C})\rtimes\mathbb{Z}/2\mathbb{Z})\times\mathbb{Z}/2\mathbb{Z}$.
\end{rem}

\begin{rem}
L'ensemble des transformations rationnelles qui commutent \`a $\Phi_\mathrm{Id}$ contient $\mathrm{Bir}(\mathcal{M}(2;\mathbb{C});\Phi_\mathrm{Id})$ mais aussi les transformations suivantes
\begin{align*}
&\mathrm{M}\mapsto\mathrm{M}^k, && \mathrm{M}\mapsto(\det\mathrm{M})^k\mathrm{Id}, && \mathrm{M}\mapsto(\det\mathrm{M}^k)\left[
\begin{array}{cc} 1 & 0 \\ 0 & 0\end{array}\right], &&  k\in\mathbb{Z}.
\end{align*}
\end{rem}

\begin{rem}
On v\'erifie que si $\Phi$ appartient \`a $\mathrm{Bir}(\mathcal{M}(2;\mathbb{C});\Phi_\mathrm{Id})$, alors $\Phi$ commute \`a toutes les transformations $\Phi_k\colon\mathrm{M}\mapsto\mathrm{M}^k$, autrement dit on a les inclusions $$\mathrm{Bir}(\mathcal{M}(2;\mathbb{C});\Phi_\mathrm{Id})\subset \mathrm{Bir}(\mathcal{M}(2;\mathbb{C});\Phi_k), \,\, k\geq 2.$$ Ces inclusions ne sont pas des \'egalit\'es: $\mathrm{M}\mapsto -\mathrm{M}$ commute \`a $\Phi_3$ mais pas \`a $\Phi_\mathrm{Id}$.
\end{rem}

\begin{rem}
Soit $\Phi\colon\mathcal{M}(2;\mathbb{C})\dashrightarrow\mathcal{M}(2;\mathbb{C})$ une application compatible \`a la multiplication: $\Phi(\mathrm{M}_1\mathrm{M}_2)=\Phi(\mathrm{M}_1) \Phi(\mathrm{M}_2)$. Un tel $\Phi$ commute \`a $\Phi_\mathrm{Id}$ et est donc dans $\mathrm{Bir}(\mathcal{M}(2;\mathbb{C});\Phi_\mathrm{Id})$ d\`es qu'il est inversible. Visiblement tous les $\sigma_\mathrm{P}$ et $\mathrm{M}\mapsto \transp\,\mathrm{M}^{-1}$ ont cette propri\'et\'e; par contre les involutions $\mathcal{T}$ et $\iota$ ne l'ont pas. Tous les \'el\'ements de $\mathrm{Bir}(\mathcal{M} (2;\mathbb{C});\Phi_\mathrm{Id})$ dont la d\'ecomposition en les $\sigma_\mathrm{P},$ $\mathcal{T}$ et $\iota$ font appara\^itre autant de fois $\mathcal{T}$ et $\iota$ satisfont cette propri\'et\'e.
\end{rem}

\begin{prob}
Quelles sont les transformations rationnelles qui commutent \`a $\Phi_\mathrm{Id}$ ?
\end{prob}

\section{\og Perturbations sp\'eciales\fg\, des applications monomiales compatibles \`a la conjugaison}\label{diagonalisable}

\noindent Nous consid\'erons dans ce qui suit quelques transformations sp\'eciales susceptibles de poss\'eder un centralisateur suffisamment \og gros\fg\, tout en pr\'esentant une dynamique plus \og riche\fg\, que celle de $\Phi_\mathrm{Id}$.

\noindent Les applications monomiales de degr\'e deux $\mathrm{M}\mapsto\mathrm{A}_1\mathrm{M} \mathrm{A}_2\mathrm{M}\mathrm{A}_3$ sont conjugu\'ees \`a celles du type $\mathrm{A}\mathrm{M}^2\mathrm{B}$ (Remarque \ref{red}); nous allons nous concentrer sur le cas sp\'ecial $\mathrm{B}=\mathrm{Id}$. Soit $\mathrm{A}$ dans $\mathrm{GL}(2;\mathbb{C})$. Consid\'erons l'application monomiale $\Phi_\mathrm{A}$ de~$\mathcal{M}(2;\mathbb{C})$ dans lui-m\^eme d\'efinie par $\Phi_\mathrm{A}(\mathrm{M})=\mathrm{A}\mathrm{M}^2$. Remarquons que $$(\sigma_\mathrm{P}^{-1} \Phi_\mathrm{A} \sigma_\mathrm{P}) (\mathrm{M})= (\mathrm{P}^{-1}\mathrm{A} \mathrm{P})\mathrm{M}^2.$$ Nous pouvons donc prendre $\mathrm{A}$ sous forme de \textsc{Jordan}. Dans un premier paragraphe nous allons consid\'erer le cas o\`u $\mathrm{A}$ est diagonale, dans le suivant nous traiterons l'autre \'eventualit\'e.

\subsection{Cas diagonalisable}

\noindent \'Ecrivons $\mathrm{A}$ sous la forme $\mathrm{diag}(\lambda_1, \lambda_2);$ remarquons que si $\lambda_1=\lambda_2,$ alors $\Phi_\mathrm{A}$ est conjugu\'ee \`a~$\Phi_\mathrm{Id}$ par une homoth\'etie. Dans la suite nous supposerons donc que $\lambda_1$ et $\lambda_2$ sont distincts. Quitte \`a conjuguer $\Phi_\mathrm{A},$ toujours par une homoth\'etie, nous pouvons supposer que $\mathrm{det}\,\mathrm{A}=1,$ \emph{i.e.} $\lambda_1\lambda_2=1$. Posons $\lambda=\lambda_1$ et $\lambda_2=\frac{1}{\lambda};$ la condition $\lambda_1\not=\lambda_2$ implique $\lambda^2\not=1$ ce que nous supposerons dans la suite. La quadrique de dimension $2$ form\'ee des matrices nilpotentes est envoy\'ee sur ${\bf 0}$ par $\Phi_\mathrm{A}$ et l'alg\`ebre $\mathfrak{sl}(2;\mathbb{C})$ est encore contract\'ee, cette fois sur $\big\{\mu\mathrm{A}\,\big\vert\,\mu\in\mathbb{C}\big\}$ qui est contenu dans $\mathcal{D}$.

\noindent La restriction de $\Phi_{\mathrm{A}}$ \`a $\mathcal{D}$ est conjugu\'ee \`a $\Phi_{\mathrm{Id}_{\vert\mathcal{D}}},$ ce que l'on voit sur l'expression de $\Phi_\mathrm{A}$: $$\left[\begin{array}{cc}x & y \\ z & t\end{array}\right]\to\left[\begin{array}{cc}\lambda(x^2+yz) & \lambda y(x+t) \\ \frac{1}{\lambda}z(x+t) & \frac{1}{\lambda}(t^2+yz)\end{array}\right].$$ 

\noindent Bien que $\Phi_{\mathrm{A}_{\vert\mathcal{D}}}$ soit conjugu\'ee \`a $\Phi_{\mathrm{Id}_{\vert\mathcal{D}}}$ nous allons voir que $\Phi_{\mathrm{A}}$ n'est pas holomorphiquement conjugu\'ee \`a $\Phi_\mathrm{Id}$. En effet, les points fixes de $\Phi_\mathrm{A}$ sont puisque $\lambda^2\not=1$ 
\begin{align*}
&{\bf 0}, &&  \mathrm{diag}\left(\frac{1}{\lambda},\lambda\right), &&\left\{\left[\begin{array}{cc}0 & 0\\z & \lambda\end{array}\right]\,\vert\, z\in\mathbb{C}\right\}, &&\left\{\left[\begin{array}{cc} \frac{1}{\lambda}&y\\0 & 0\end{array}\right]\,\vert\, y\in\mathbb{C}\right\}
\end{align*}

\noindent alors que $\Phi_\mathrm{Id}$ admet un ensemble de points fixes de dimension $2$: ceci est une obstruction \`a la conjugaison entre $\Phi_\mathrm{Id}$ et $\Phi_\mathrm{A}$.

\subsubsection{Fibrations et fonctions invariantes}\label{fibdiag}

\noindent Comme pour $\Phi_\mathrm{Id}$ la fibration $H\colon\left[\begin{array}{cc} x & y \\ z & t \end{array}\right]\dashrightarrow\frac{y}{z}$ est invariante. N\'eanmoins si $H\circ\Phi_\mathrm{Id}=H$, nous avons $H\circ\Phi_\mathrm{A}=\lambda^2 H$: il y a donc seulement deux fibres fixes, celles correspondant aux matrices triangulaires sup\'erieures, resp. inf\'erieures. Le feuilletage associ\'e \`a $y\frac{\partial}{\partial y}-z\frac{\partial}{\partial z}$ est encore invariant par $\Phi_\mathrm{A}$; en fait $\Phi_\mathrm{A}$ commute \`a toutes les transformations $(x,\mathrm{e}^sy,\mathrm{e}^{-s}z,t)$, $s\in\mathbb{C}$. En particulier l'ensemble des points p\'eriodiques, les ensembles invariants maximaux (resp. minimaux) sont eux-m\^emes invariants par le flot $y\frac{ \partial}{\partial y}-z\frac{\partial}{\partial z}$. Il en r\'esulte que la seule possibilit\'e pour que les points p\'eriodiques (de m\^eme p\'eriode) soient isol\'es est qu'ils soient confin\'es dans $\mathcal{D}$. Notons que la fibration $\det=$ cte est invariante. 

\medskip

\noindent Comme nous l'avons fait pour $\Phi_\mathrm{Id}$ nous nous int\'eressons au corps des fonctions invariantes $\mathrm{K}(\Phi_\mathrm{A})$ pour $\mathrm{A}=\mathrm{diag}\left(\frac{1}{\lambda},\lambda\right)$ g\'en\'erique. Pour cel\`a nous \'etudions $\Phi_\mathrm{A}$ au voisinage du point fixe $\mathrm{diag}\left(\frac{1}{\lambda},\lambda\right)$. On v\'erifie que la matrice jacobienne de $\Phi_\mathrm{A}$ en ce point fixe est $$\left[\begin{array}{cccc} 2 & 0 & 0 & 0 \\0 & \lambda\left(\lambda+\frac{1}{\lambda}\right) & 0 & 0 \\ 0 & 0 & \frac{1}{\lambda}\left( \lambda+\frac{1}{\lambda} \right)& 0 \\ 0 & 0 & 0 & 2 \end{array}\right].$$ Ainsi, pour $\lambda$ non r\'esonant (on demande que $2^p\lambda^q\left(\lambda+\frac{1}{\lambda} \right)^r=1$, $(p,q,r)\in\mathbb{Z}$ ait une seule solution, la solution triviale $(p,q,r)=(0,0,0)$), le germe de $\Phi_\mathrm{A}$ en $\mathrm{diag}\left(\lambda,\frac{1}{\lambda}\right)$ est formellement lin\'earisable. Nous en d\'eduisons qu'\`a conjugaison formelle pr\`es les fonctions m\'eromorphes (formelles) invariantes par le germe de $\Phi_\mathrm{A}$ en $\mathrm{diag}\left( \frac{1}{\lambda},\lambda \right)$ sont les fonctions du type $h(x/t),$ $h$ rationnelle. Soit maintenant $f$ dans $\mathrm{K}(\Phi_\mathrm{A})$ non constante. Il n'est pas difficile de voir, en vertu de ce qui pr\'ec\`ede, que $f$ est non constante sur le $2$-plan $y=z=0$. Mais dans ce $2$-plan la restriction de $\Phi_\mathrm{A}$ est donn\'ee par $(x,t)\mapsto\left(\lambda x^2, \frac{1}{\lambda}t^2\right)$ qui est conjugu\'ee \`a $(x^2,t^2)$. Or les fonctions invariantes par $(x^2,t^2)$ sont constantes. Il s'en suit la:

\begin{pro}
Pour $\lambda$ non r\'esonant $\mathrm{K}(\Phi_\mathrm{A})=\mathbb{C}$.
\end{pro}

\begin{rem}
Pour certaines valeurs sp\'eciales de $\lambda$ le corps $\mathrm{K}(\Phi_\mathrm{A})$ ne se r\'eduit pas aux constantes. En effet comme on l'a vu $H\circ\Phi_\mathrm{A}=\lambda^2 H$; en particulier si $\lambda$ est une racine $2k$-i\`eme de l'unit\'e, alors la fonction $(y/z)^k$ est invariante par $\Phi_\mathrm{A}$. 
\end{rem}

\begin{prob}
D\'ecrire le corps $\mathrm{K}(\Phi_\mathrm{A})$ pour toutes les valeurs de $\lambda$.
\end{prob}

\subsubsection{\'Etude des points p\'eriodiques pour $\lambda$ g\'en\'erique} 

\noindent Suivant la nature de $\lambda$ le comportement des points p\'eriodiques s'av\`ere diff\'erent, et diff\'erent aussi du cas $\Phi_\mathrm{Id}$ \'etudi\'e pr\'ec\'edemment. Supposons que $\lambda$ ne soit pas racine de l'unit\'e. Puisque $H\circ\Phi_\mathrm{A}=\lambda^2 H$ les points p\'eriodiques de $\Phi_\mathrm{A}$ sont contenus dans les hyperplans $z=0$ et $y=0$. La restriction $\Phi_{\mathrm{A}_{\vert z=0}}$ est conjugu\'ee \`a $\Psi (x,y,t)=\Big(x^2,y(x+\lambda^2t),t^2\Big);$ on constate que $$\Psi^n(x,y,t)=\Big(x^{2^n},y\displaystyle\prod_{i=0}^{n-1}(x^{2^i}+\lambda^2t^{2^i}),t^{2^n}\Big).$$
 
\noindent Notons que si $(\underline{x},\underline{y},\underline{t})$ est p\'eriodique pour $\Psi,$ alors $(\underline{x},\underline{t})$ l'est pour $\chi_2$; il y a donc quatre possibilit\'es pour~$(\underline{x},\underline{t})$:
\begin{itemize}
\item[$\mathfrak{a)}$] $(\underline{x},\underline{t})=(0,0)$;

\item[$\mathfrak{b)}$] $\underline{x}$ et $\underline{t}$ sont des racines $(2^n-1)$-i\`eme de l'unit\'e; 

\item[$\mathfrak{c)}$] $\underline{x}=0$ et $\underline{t}$ est une racine $(2^n-1)$-i\`eme de l'unit\'e; 

\item[$\mathfrak{d)}$] $\underline{t}=0$ et $\underline{x}$ est une racine $(2^n-1)$-i\`eme de l'unit\'e.
\end{itemize} 
 
\noindent Examinons ces possibilit\'es au cas par cas.
 
\noindent $\mathfrak{a)}$ Notons que si $(\underline{x},\underline{t})=(0,0),$ alors $\underline{y}=0$. 

\noindent $\mathfrak{b)}$ Consid\'erons l'ensemble d\'enombrable $\Lambda_n$ d\'efini par $$\Lambda_n=\Big\{(x,\,t,\,\lambda)\in\mathbb{C}^3\,\big\vert\,x^{2^n -1}=1,\,t^{2^n-1}=1,\,\lambda\displaystyle\prod_{i=0}^{n-1}\Big(x^{2^i}+\lambda^2t^{2^i}\Big)=1\Big\}$$ et  l'ensemble $\Lambda=\mathrm{pr}_3(\cup_n\Lambda_n)$ o\`u $\mathrm{pr}_3$ d\'esigne la troisi\`eme projection.

\noindent Un argument de \textsc{Baire} assure que pour $\lambda$ g\'en\'erique on a $\lambda^{n-1}\displaystyle \prod_{i=0}^{n-1}\Big(\underline{x}^{2^i}+\lambda^2\underline{t}^{2^i}\Big)\not=1$ pour tout choix $\underline{x}$ et $\underline{t}$ de racines $(2^n-1)$-i\`eme de l'unit\'e. Par suite, pour $\lambda$ n'appartenant pas \`a $\Lambda$, les points p\'eriodiques $(\underline{x}, \underline{y},\underline{t})$ de $\Psi$ tels que $\underline{x}\underline{t}\not=0$ sont exactement l'ensemble $$\displaystyle \bigcup_{n\geq 0}\Big\{(\xi,0, \eta)\,\big\vert\,\xi^{2^n-1}=\eta^{2^n-1}=1\Big\}.$$

\noindent Lorsque $(\underline{x},\underline{t},\lambda)$ appartient \`a $\Lambda_n$ (en particulier $\lambda$ appartient \`a $\Lambda$), $(\underline{x},y,\underline{t})$ est p\'eriodique pour $\Psi,$ pour tout~$y$. Par exemple $(\mathrm{j},\mathrm{j}^2,2)$ appartient \`a $\Lambda_2$ et $(\mathrm{j},y,0,\mathrm{j}^2)$ est un point p\'eriodique de p\'eriode $2$ de l'application correspondante.

\noindent $\mathfrak{c)}$ Passons maintenant au cas o\`u par exemple $\underline{t}=0$ et $\underline{x}^{2^n-1}= 1$. On v\'erifie que $\Psi^n(\underline{x},y,0)=(\underline{x},y,0)$. 

\noindent $\mathfrak{d)}$ Lorsque $\underline{x}=0$ et $\underline{t}^{2^n-1}=1$ on constate que $\Psi^n(0,y,\underline{t})=(0,\lambda^2y,\underline{t})$. 

\noindent Bien s\^ur tout ce qui est dit pour la restriction de $\Phi_\mathrm{A}$ \`a $z=0$ peut \^etre r\'ep\'et\'e pour la restriction \`a $y=0$. Nous en d\'eduisons la:

\begin{pro}
Pour $\lambda$ g\'en\'erique $(${\it i.e.} $\lambda$ n'appartient pas \`a $\Lambda)$ l'adh\'erence des points p\'eriodiques de~$\Phi_\mathrm{A}$ est constitu\'ee  
\begin{itemize}
\item d'un tore $\mathbb{S}^1\times\mathbb{S}^1$ contenu dans $y=z=0;$

\item de deux $\mathbb{S}^1\times\mathbb{C},$ pr\'ecis\'ement $\frac{1}{\lambda}\mathbb{S}^1\times\mathbb{C}\times\big\{0\big\}\times\big\{0\big\}$ et $\big\{0\big\}\times\big\{0\big\}\times\mathbb{C}\times\lambda\mathbb{S}^1$ $($en identifiant $\mathcal{M}(2;\mathbb{C})=\left\{\left[\begin{array}{cc}x & y \\ z & t\end{array}\right]\right\}$ \`a $\mathbb{C}^4=\big\{(x,y,z,t)\,\big\vert\, x,\,y,\,z,\, t\in\mathbb{C}\big\})$; 

\item de la matrice nulle ${\bf 0}$.
\end{itemize}
\end{pro}

\noindent En particulier lorsque $\lambda$ n'est pas dans $\Lambda$ les points p\'eriodiques sont d'adh\'erence de \textsc{Zariski} $\big\{y=z=0\big\}\cup\big\{z=t=0\big\}\cup\big\{x=y=0\big\}$ alors que les points p\'eriodiques de $\Phi_\mathrm{Id}$ sont \textsc{Zariski} denses. La description des points p\'eriodiques de $\Phi_\mathrm{A}$ dans le cas non g\'en\'erique (par exemple $\lambda$ racine de l'unit\'e) semble d\'elicate.

\begin{rem}
La dynamique de $\Phi_\mathrm{A}\colon\mathcal{M}(2;\mathbb{C})\to\mathcal{M}(2; \mathbb{C})$ peut \^etre pr\'ecis\'ee par l'\'etude de ses points fixes ou plus g\'en\'eralement p\'eriodiques \`a l'infini. Pour cel\`a consid\'erons l'application $\widetilde{\Phi_\mathrm{A}}\colon\mathbb{P}^4(\mathbb{C})\dashrightarrow\mathbb{P}^4(\mathbb{C})$ d\'efinie par $$(x:y:z:t:s)\mapsto\Big(\lambda(x^2+yz):\lambda y(x+t):\frac{z(x+t)}{\lambda}:\frac{t^2+yz}{\lambda}:s^2\Big).$$  On v\'erifie que $\widetilde{\Phi_\mathrm{A}}$ s'exprime de la fa\c{c}on suivante dans la carte $x=1$ $$\left(\frac{y(1+t)}{1+yz},\frac{z(1+t)}{\lambda^2(1+yz)},\frac{t^2+yz}{\lambda^2(1+yz)},\frac{s^2}{\lambda(1+yz)}\right).$$ Dans l'hyperplan $y=0$ qui est invariant, le point $(z=0,t=0,s=0)$ est un attracteur si $\vert \lambda\vert>1$. Dans la car\-te~$t=1$, o\`u $\widetilde{\Phi_\mathrm{A}}$ s'exprime cette fois sous la forme $$\left(\frac{\lambda^2(x^2+yz)}{1+yz},\frac{\lambda^2y(1+x)}{1+yz},\frac{z(1+x)}{1+yz},\frac{\lambda s^2}{1+yz}\right),$$ on voit que $z=0$ est invariant et dans ce $3$-plan le point $(x=0,y=0,s=0)$ est un attracteur si $\vert\lambda\vert<1$. Ainsi si $\vert\lambda\vert\not=1$ on trouve un hyperplan invariant contenant un point fixe attractant dans cet hyperplan. De m\^eme on peut pr\'eciser la nature des points p\'eriodiques \`a l'infini, tout du moins dans ces hyperplans invariants.
\end{rem}

\begin{prob}
D\'ecrire les points p\'eriodiques des applications $\Phi_\mathrm{A}$ pour toute valeur de $\lambda$.
\end{prob}

\subsubsection{\'Etude de quelques orbites non p\'eriodiques}\label{perper}

\noindent Nous allons maintenant pr\'eciser la nature de quelques orbites non p\'eriodiques. Dans toute \'etude de syst\`emes dynamiques on s'int\'eresse \`a des probl\`emes de stabilit\'e, de \og non-explosion\fg, ce qui conduit \`a caract\'eriser autant que faire se peut les orbites born\'ees. L'id\'eal serait d'en avoir une description compl\`ete pour chaque $\Phi_\mathrm{A}$. Nous allons simplement en donner quelques constructions. Puisque nous pouvons nous ramener \`a $\lambda_1\lambda_2=1$ nous avons $\det\Phi_\mathrm{A}(\mathrm{M})= (\det\mathrm{M})^2$ et par suite $\det\Phi_\mathrm{A}^k(\mathrm{M})= (\det\mathrm{M})^{2^k}$. Il en r\'esulte que pour $\mathrm{M}$ tel que $\vert\det\mathrm{M}\vert>1,$ la suite $(\Phi_\mathrm{A}^k(\mathrm{M}))_k$ tend vers l'infini. Pour $\mathrm{A}$ fix\'e consid\'erons le bassin d'attraction $\mathrm{W}_\mathrm{A}^s({\bf 0})$ de la matrice nulle; c'est un domaine disqu\'e, \emph{i.e.} si~$\mathrm{M}$ appartient \`a $\mathrm{W}_\mathrm{A}^s({\bf 0})$, alors le disque $\big\{\mu\mathrm{M}\,\big\vert\,\vert\mu\vert\leq 1\big\}$ est lui aussi dans $\mathrm{W}_\mathrm{A}^s({\bf 0})$. Le fait que $\Phi_\mathrm{A}$ commute aux applications lin\'eaires $f_s(x,y,z,t)=(x,\mathrm{e}^sy,\mathrm{e}^{-s}z,t)$ implique que $\mathrm{W}_\mathrm{A}^s({\bf 0})$ est non born\'e. D\'esignons par $\vert\vert\, .\vert\vert$ la norme sup sur $\mathcal{M}(2;\mathbb{C})$. Soient $K=\sup\left(\lambda,\frac{1}{\lambda}\right)$ et $\mathrm{M}$ dans le polydisque $\Delta(\rho)$ de rayon $\rho$; nous avons l'in\'egalit\'e $\vert\vert\Phi_\mathrm{A}^k(\mathrm{M}) \vert\vert\leq 2K\rho^2$ et par cons\'equent $\vert\vert\Phi_\mathrm{A}^k(\mathrm{M})\vert\vert\leq  (2K)^{2^k-1}\rho^{2^k}$. Il en r\'esulte l'inclusion du polydisque $\Delta\left(\frac{1}{2K}\right)$ dans $\mathrm{W}_\mathrm{A}^s({\bf 0})$ ainsi que de son satur\'e $\cup_s f_s\left(\Delta\left(\frac{1}{2K}\right)\right)$ par le flot $f_s$. Ceci donne d'ailleurs une preuve du fait que $\mathrm{W}_\mathrm{A}^s({\bf 0})$ est ouvert. Il est assez simple de produire des orbites born\'ees dans les plans invariants $y=0$ et~$z=0$. Pla\c{c}ons-nous par exemple dans l'hyperplan invariant $z=0$; nous travaillons de nouveau avec $\Psi$ dont l'it\'er\'e $n$-i\`eme s'\'ecrit pour $x$ non nul $$\Psi^n(x,y,t)=\Big(x^{2^n},yx^{2^n-1}\prod_{i=0}^{n-1}\left(1+\lambda^2\left(\frac{t}{x}\right)^{2^i}\right),t^{2^n}\Big).$$ En particulier les orbites des points $(x,y,0)$ sont faciles \`a d\'ecrire et produisent des orbites born\'ees lorsque $\vert x\vert=1$ (ou $0$). Remarquons que si $x$ et $t$ sont fix\'es tels que $\vert t\vert<\vert x\vert,$ alors le produit infini $\displaystyle\prod_{i=0}^{\infty}\left(1+\lambda^2\left(\frac{t}{x}\right)^{2^i}\right)$ converge vers un nombre $\mu=\mu\left(\frac{t}{x}\right)$. Consid\'erons l'application $\Psi_0$ d\'efinie par $\Psi_0(x,y,t)=(x^2,xy,t^2)$ (\emph{i.e.}  qui correspond au cas $\lambda=0$). Nous avons $$\vert\vert\Psi^n(x,y,t)-\Psi_0^n(x,y,t)\vert\vert=\vert y\vert\cdot\big\vert x^{2^n-1}\big\vert\cdot\Big\vert \prod_{i=0}^{n-1}\left(1+\lambda^2\left(\frac{t}{x}\right)^{2^i} \right)-1\Big\vert.$$ En particulier pour $\vert t \vert<\vert x\vert< 1$, nous avons $\displaystyle\lim_{n\to +\infty}\Psi^n(x,y,t)=0$ ind\'ependamment de $y$. Ainsi l'ensemble $\big\{(x,y,t)\,\big\vert\, \vert t\vert <\vert x\vert< 1\big\}$ est contenu dans le bassin d'attraction de ${\bf 0}$ pour l'application $\Psi$ ce qui donne un renseignement suppl\'ementaire pour $\mathrm{W}_\mathrm{A}^s({\bf 0})$. Lorsque $\vert\underline{x}\vert=1$, $\vert\underline{t}\vert<1$ et $y$ quelconque $\mathcal{O}\big((\underline{x},\underline{y},\underline{t});\Psi_0\big)$ est contenue dans $\vert\underline{x}\vert=1,$ $\vert y\vert=\vert\underline{y}\vert$ alors que $\displaystyle\lim_{n\to+\infty}\underline{t}^{2^n}=0$. Dans cette situation on constate que $\vert\vert\Psi^n(\underline{x},\underline{y}, \underline{t})-\Psi_0^n( \underline{x},\underline{y},\underline{t})\vert\vert$ est born\'ee et donc que $\mathcal{O}\big((\underline{x},\underline{y},\underline{t});\Psi\big)$ est born\'ee.

\noindent Remarquons que le plan $x=t$ (qui est dans le bord du domaine $\vert t\vert<\vert x\vert$) est invariant par $\Psi$. La restriction de $\Psi$ \`a ce plan s'\'ecrit $\Psi(x,y)=(x^2,(1+\lambda^2)xy)$; elle est conjugu\'ee via l'application birationnelle $(x,xy)$ \`a $\varphi=(x^2,(1+\lambda^2)y)$ dont l'it\'er\'e $n$i\`eme est $(x^{2^n},(1+\lambda^2)^ny)$. Pour tout $x$ non nul, la nature des orbites de $\Psi_{\vert x=t}$ peut se d\'eduire de celles de $\varphi$; plus pr\'ecis\'ement $\varphi$ contracte $x=t=0$ (l'axe des $y$) sur la matrice nulle. En dehors de $x=0$, $(x^2,(1+\lambda^2)xy)$ est holomorphiquement conjugu\'e \`a $(x^2,(1+\lambda^2)y)$. Par exemple lorsque $\underline{x}$ est de module $1$ g\'en\'erique et $1+\lambda^2$ est aussi de module $1$ g\'en\'erique, les orbites de $(\underline{x},y)$ sont d'adh\'erence des tores r\'eels de dimension $2$.

\noindent Une fa\c{c}on plus pr\'ecise d'appr\'ehender l'\'etude de l'application $\Psi\sim\Phi_{ \mathrm{A}\vert z=0}$ est de la conjuguer par une transformation birationnelle bien choisie. Par exemple si $E\colon\mathbb{C}^3\to\mathbb{C}^3$ est l'application d'\'eclatement de l'origine d\'efinie par $E(x,u,v)=(x,xu,xv)$ nous avons le diagramme commutatif
$$\xymatrix{\mathbb{C}^3\ar[d]_{E} \ar[rr]^\Theta& &\mathbb{C}^3\ar[d]^{E}\\
\mathbb{C}^3\ar[rr]_\Psi & & \mathbb{C}^3}$$ avec $\Theta(x,u,v)=(x^2,u(1+\lambda^2v),v^2)$. L'it\'er\'e $n$-i\`eme de $\Theta$ est donn\'e par 
\begin{align*}
&\Theta^n(x,u,v)=(x^{2^n},uT_n^\lambda(v), v^{2^n}), && T_n^\lambda(v)=\prod_{i=0}^{n-1} (1+\lambda^2v^{2^i}).
\end{align*}

\noindent Contr\^oler les $v$ appartenant \`a $\mathbb{S}^1$ pour lesquels les produits $\prod_{i=0}^{n-1} (1+\lambda^2v^{2^i})$ restent born\'es permet d'exhiber des orbites born\'ees de $\Theta$ et une orbite born\'ee de $\Theta$ induit une orbite born\'ee de $\Psi$. En particulier chaque fois que l'on sait borner l'ensemble $\big\{T_n^\lambda(v)\,\big\vert\, n\in\mathbb{N}\big\}$ pour certaines valeurs de $\lambda$ et $v$ on sait construire une orbite born\'ee de $\Psi$. Voici un exemple: supposons que $\lambda$ soit un r\'eel, $0<\lambda<1,$ et que $v$ soit une racine cubique de l'unit\'e, par exemple $\mathrm{j}$. Nous avons 
\begin{align*}
& T_1^\lambda(\mathrm{j})=(1+\lambda^2\mathrm{j}), && T_2^\lambda(\mathrm{j})=(1+\lambda^2\mathrm{j})(1+\lambda^2\mathrm{j}^2), && T_3^\lambda(\mathrm{j})=(1+\lambda^2\mathrm{j})^2(1+\lambda^2\mathrm{j}^2),\\
&T_4^\lambda(\mathrm{j})=(1+\lambda^2\mathrm{j})^2(1+\lambda^2\mathrm{j}^2)^2, && \ldots
\end{align*}
Une induction \'el\'ementaire montre que pour tout $n\geq 0$ les $\vert T_n^\lambda(\mathrm{j})\vert$ sont strictement plus petits que $1$. Par suite si $\vert x\vert\leq 1$, alors $\mathcal{O}\big((x,y,\mathrm{j});\Phi_\mathrm{A}\big)$  est born\'ee (toujours sous la condition $0<\lambda<1$).

\begin{prob}
Pour $\lambda$ fix\'e donner les $v$ de module inf\'erieur ou \'egal \`a $1$ pour lesquels $\big\{T_n^\lambda(v)\,\big\vert\, n\in\mathbb{N}\big\}$ est born\'e.
\end{prob}

\noindent Soit $\mathrm{M}$ dans $\mathcal{M}(2;\mathbb{C});$ supposons que $\mathcal{O}(\mathrm{M};\Phi_\mathrm{A})$ soit born\'ee. Alors les points limites de $\mathcal{O}(\mathrm{M};\Phi_\mathrm{A})$ sont encore \`a orbites born\'ees. Si $\vert\det\mathrm{M}\vert =1$, alors $\overline{\mathcal{O}(\mathrm{M};\Phi_\mathrm{A})}$ est contenu dans $\big\{\mathrm{M}\in\mathcal{M}(2;\mathbb{C})\,\big\vert\,\det\mathrm{M}=0\big\}$, tandis que si $\vert\det\mathrm{M}\vert<1$ ces points limites sont dans l'hypersurface $\big\{\mathrm{M}\in\mathcal{M}(2;\mathbb{C})\,\big\vert\,\det\mathrm{M}=0\big\}$. Il est donc naturel de rechercher les orbites born\'ees dans les deux ensembles 
\begin{align*}
&\big\{\mathrm{M}\in\mathcal{M}(2;\mathbb{C})\,\big\vert\,\det\mathrm{M}=0\big\} &&\text{ et } &&\big\{\mathrm{M}\in\mathcal{M}(2;\mathbb{C})\,\big\vert\,\vert\det\mathrm{M}\vert=~1\big\}
\end{align*}

\noindent ce que nous aborderons dans ce qui suit.

\bigskip

\subsubsection{Orbites born\'ees dans $\big\{\mathrm{M}\in\mathcal{M}(2;\mathbb{C})\,\big\vert\,\det\mathrm{M}=0\big\}$}

\noindent Soit $\mathrm{M}=\left[\begin{array}{cccc} x & y \\ z & t\end{array}\right]$ dans $\mathcal{M}(2;\mathbb{C})$  de d\'eterminant nul; on constate que $\Phi_\mathrm{A}(\mathrm{M})=(x+t)\left(\lambda x,\lambda y,\frac{z}{\lambda},\frac{t}{\lambda}\right)$. L'application $\mathbb{P}(\Phi_\mathrm{A})\colon \mathbb{P}^3(\mathbb{C})\dashrightarrow\mathbb{P}^3(\mathbb{C})$ co\"{i}ncide en restriction \`a la quadrique $\big\{\mathrm{M}\in\mathcal{M}(2;\mathbb{C})\, \big\vert\,\det\mathrm{M}=0\big\}$ avec $(x:y:z:t)\mapsto(\lambda x:\lambda y:z/ \lambda:~t/\lambda)$ qui est lin\'eaire. Plus g\'en\'eralement nous avons 
\begin{align*}
&\Phi_\mathrm{A}^n(\mathrm{M})=P_n\cdot\left(\lambda^nx,\lambda^ny, \frac{z}{\lambda^n},\frac{t}{\lambda^n}\right), && P_n=\displaystyle\prod_{i=0}^{n-1} \left(\lambda^ix+\frac{t}{\lambda^i}\right)^{2^{n-i-1}}.
\end{align*}

\noindent En particulier on constate le ph\'enom\`ene de r\'esonance suivant: si $\mathrm{M}=\left[\begin{array}{cccc} x & y \\ z & t\end{array}\right]$ dans $\mathcal{M}(2;\mathbb{C})$ satisfait $xt-yz=0$ et $\lambda^kx+\frac{t}{ \lambda^k}=0$ pour un certain $k$, alors $\Phi_\mathrm{A}^n(\mathrm{M})={\bf 0}$ pour $n\geq k$. Nous obtenons donc, tout du moins pour $\lambda$ g\'en\'erique (pr\'ecis\'ement pour $\lambda$ non racine de l'unit\'e), une infinit\'e de surfaces quadratiques $$\mathcal{Q}_k(\mathrm{A})=\left\{\left[\begin{array}{cc} x & y \\ z & t\end{array}\right]\in\mathcal{M}(2; \mathbb{C})\,\big\vert\,\det\mathrm{M}=0,\,\lambda^kx+t/\lambda^k=0\right\}$$ qui sont envoy\'ees sur ${\bf 0}$ apr\`es un nombre fini d'it\'erations. Comme le bassin d'attraction de ${\bf 0}$ est ouvert, il existe des voisinages ouverts de ces surfaces quadratiques contenus dans $\mathrm{W}_\mathrm{A}^s({\bf 0}),$ ce qui produit \'evidemment des orbites born\'ees. Rappelons que dans $\big\{\mathrm{M}\in\mathcal{M}(2; \mathbb{C})\,\big\vert\,\det\mathrm{M}=0\big\}$ il y a des orbites non born\'ees (celles de $(x,y,0,0)$ avec $\vert x\vert$ grand) et des orbites born\'ees non contenues dans $\mathrm{W}_\mathrm{A}^s({\bf 0})$, par exemple celles des points $(x,y,0,0)$ avec $\vert x\vert=1/\lambda$ qui sont toutefois dans le bord de $\mathrm{W}_\mathrm{A}^s({\bf 0})$. 
\medskip

\noindent Sur $\text{http://math.cmaisonneuve.qc.ca/alevesque/chaos}\_\text{fract/Julia/Julia.html}$ on trouve un programme permettant de tracer l'ensemble de \textsc{Julia} d'une application holomorphe de $\mathbb{P}^1(\mathbb{C})$ dans lui-m\^eme.
Nous proposons une \og adaptation\fg\, de ce programme \`a certaines transformations polynomiales r\'eelles. Plus pr\'ecis\'ement notons $\mathbb{D}(0,r)\subset\mathbb{R}^2$ le disque de rayon $r$ centr\'e en $0$ et $\Delta(\rho)$ le polydisque de rayon $\rho$:
\begin{align*}
\Delta(\rho)=\big\{(y,t)\in\mathbb{R}^2\,\big\vert\,
\vert y\vert<\rho,\,\vert t\vert<\rho\big\}.
\end{align*}

\noindent Soit $f$ une transformation polynomiale du plan $\mathbb{R}^2$. Soient $m=(x_0,x_1)$ un
point de $\Delta(\rho)$ et $\kappa$ un entier strictement positif. On appelle temps de sortie $N(m;r,\rho,\kappa),$ relatif aux donn\'ees de contr\^ole $r$, $\rho$, $\kappa$, du point~$m$ de l'intersection $\Delta(\rho)\cap\mathbb{D}(0,r)$ le plus grand entier $n$ dans $[0,\ldots,\kappa]$ tel que
\begin{align*}
&f^k(m)\in\mathbb{D}(0,r)
&&\forall\, 0\leq k\leq n.
\end{align*}
Consid\'erons le spectre (continu) des couleurs $[\text{rouge}\ldots\text{orange}\ldots\text{jaune}\ldots\text{vert}\ldots\text{bleu}\ldots \text{indigo}\ldots\text{rouge}]$ que l'on discr\'etise en $\kappa+1$ intervalles $[I_0,\ldots, I_\kappa]=[\text{rouge}\ldots\text{jaune}\ldots\text{bleu}\ldots\text{rouge}]$. Soit $(r,\rho,\kappa)$ un triplet de contr\^ole. Si $N(m;r,\rho,\kappa)=k$,  on colore le point $m$ de la couleur $\mathrm{Col}\,(m)=I_k$. Comme $I_0=I_\kappa=$ rouge, les points color\'es en rouge sont ceux pour lesquels le temps de passage est $0$ (sortie imm\'ediate) ou $\kappa$ (pas de sortie au bout de $\kappa$ it\'erations); sur les figures le bord de l'ensemble $\rm{Col}^{-1}\,(I_0)$ est approch\'e par la couleur $I_1\sim$ orang\'e.

\noindent Nous allons appliquer cette proc\'edure \`a $\left(\lambda x(x+t),\frac{t(x+t)}{\lambda}\right)$, application qui d\'ecrit la dynamique de la restriction de $\Phi_\mathrm{A}$ \`a $\big\{\mathrm{M}\in\mathcal{M}(2;\mathbb{R})\,\big\vert\,\det\mathrm{M}=0\big\}$. Nous obtenons pour $\rho=10$, $r=30$

\smallskip

\begin{center}
\begin{tabular}{llll}
\begin{picture}(0,0)%
\includegraphics{dessin1.pstex}%
\end{picture}%
\setlength{\unitlength}{3947sp}%
\begingroup\makeatletter\ifx\SetFigFont\undefined%
\gdef\SetFigFont#1#2#3#4#5{%
  \reset@font\fontsize{#1}{#2pt}%
  \fontfamily{#3}\fontseries{#4}\fontshape{#5}%
  \selectfont}%
\fi\endgroup%
\begin{picture}(1050,1050)(1201,-2611)
\end{picture}%
\hspace{20mm}&\hspace{20mm}\begin{picture}(0,0)%
\includegraphics{dessin2.pstex}%
\end{picture}%
\setlength{\unitlength}{3947sp}%
\begingroup\makeatletter\ifx\SetFigFont\undefined%
\gdef\SetFigFont#1#2#3#4#5{%
  \reset@font\fontsize{#1}{#2pt}%
  \fontfamily{#3}\fontseries{#4}\fontshape{#5}%
  \selectfont}%
\fi\endgroup%
\begin{picture}(1050,1050)(1201,-2611)
\end{picture}%
\hspace{20mm}&\hspace{20mm}\begin{picture}(0,0)%
\includegraphics{dessin3.pstex}%
\end{picture}%
\setlength{\unitlength}{3947sp}%
\begingroup\makeatletter\ifx\SetFigFont\undefined%
\gdef\SetFigFont#1#2#3#4#5{%
  \reset@font\fontsize{#1}{#2pt}%
  \fontfamily{#3}\fontseries{#4}\fontshape{#5}%
  \selectfont}%
\fi\endgroup%
\begin{picture}(1050,1050)(1201,-2611)
\end{picture}%
\hspace{20mm}&\hspace{20mm}\begin{picture}(0,0)%
\includegraphics{dessin4.pstex}%
\end{picture}%
\setlength{\unitlength}{3947sp}%
\begingroup\makeatletter\ifx\SetFigFont\undefined%
\gdef\SetFigFont#1#2#3#4#5{%
  \reset@font\fontsize{#1}{#2pt}%
  \fontfamily{#3}\fontseries{#4}\fontshape{#5}%
  \selectfont}%
\fi\endgroup%
\begin{picture}(1050,1050)(1201,-2611)
\end{picture}%
\\
$\kappa=10$, $\lambda=1$&\hspace{20mm}$\kappa=75$, $\lambda=1$&\hspace{20mm}$\kappa=10$, $\lambda=1.5$&\hspace{20mm}$\kappa=75$, $\lambda=1.5$\\
\end{tabular}
\end{center}

\smallskip

\noindent Les deux premi\`eres figures repr\'esentent une approximation de la projection sur le plan des $(x,t)$ du bassin d'attraction de ${\bf 0}$ par~$\Phi_\mathrm{Id}$ intersect\'e avec $\big\{\mathrm{M}\in\mathcal{M}(2;\mathbb{R})\,\big\vert\,\det\mathrm{M}=1\big\}$. On v\'erifie en effet que, dans ce cas, la bande $\vert x+t\vert < 1$ est exactement le bassin d'attraction de l'origine pour l'application consid\'er\'ee. Le domaine \og\'etoil\'e\fg\, des troisi\`eme et quatri\`eme dessins repr\'esente comme pr\'ec\'edemment  la projection sur le plan des~$(x,t)$ de $\mathrm{W}_\mathrm{A}^s{\bf 0}$ intersect\'e avec $\big\{\mathrm{M}\in\mathcal{M}(2;\mathbb{R})\,\big\vert\,\det\mathrm{M}=0\big\}$. On distingue ici quelques droites $\lambda^kx+\frac{t}{\lambda^k}=0$ et leur voisinage contenus dans le bassin d'attraction. Au vu de ces figures on peut penser que le bord du bassin d'attraction $\mathrm{W}_\mathrm{A}^s({\bf 0})$ n'est plus \textsc{Levi}-plat.

\subsubsection{Orbites born\'ees dans $\big\{\mathrm{M}\in \mathcal{M}(2;\mathbb{C})\,\big\vert\,\vert\det\mathrm{M}\vert=1\big\}$}

\noindent La dynamique dans $\big\{\mathrm{M}\in \mathcal{M}(2;\mathbb{C})\,\big\vert\,\vert\det\mathrm{M}\vert=~1\big\}$ semble difficile d'abord pour $\lambda$ quelconque. Nous nous contentons de quelques remarques concernant l'hypersurface invariante $\big\{\mathrm{M}\in \mathcal{M}(2;\mathbb{C})\,\big\vert\,\det\mathrm{M}=1\big\}$. La matrice $\mathrm{M}_0=\left[ \begin{array}{cc}\frac{1}{\lambda} &0 \\ 0&\lambda \end{array}\right]$ est fixe pour la transformation~$\Phi_\mathrm{A}$ et appartient \`a la quadrique $\big\{\mathrm{M}\in \mathcal{M}(2;\mathbb{C})\,\big\vert\,\det\mathrm{M}=1\big\}$. La matrice jacobienne de~$\Phi_\mathrm{A}$ en $\mathrm{M}_0$ est la suivante $$\mathrm{Jac}(\Phi_\mathrm{A})_{(\mathrm{M}_0)}=\left[\begin{array}{cccc}2 & 0 & 0 & 0\\ 0 & \lambda\left(\lambda+\frac{1}{\lambda}\right) & 0 & 0\\
0 & 0 & \frac{1}{\lambda}\left(\lambda+\frac{1}{\lambda}\right) & 0 \\
0 & 0 & 0 & 2\end{array}\right];$$

\noindent on constate un ph\'enom\`ene de r\'esonance entre les valeurs propres. Remarquons que si l'argument de $\lambda$ appartient \`a $]-\pi/4,\pi/4[$, alors les valeurs propres de $\mathrm{Jac}(\Phi_\mathrm{A})_{(\mathrm{M}_0)}$ sont en module strictement sup\'erieures \`a~$1$. Sous cette hypoth\`ese le point fixe $\mathrm{M}_0$ est un r\'epulseur, \emph{i.e.} il existe un voisinage $\mathcal{V}( \mathrm{M}_0)$ de $\mathrm{M}_0$ tel que $\mathcal{V}(\mathrm{M}_0)\subsetneq\Phi_\mathrm{A}( \mathcal{V}(\mathrm{M}_0))$. Mais ce fait n'est pas universel puisque pour certaines valeurs de $\lambda$ on trouve des valeurs propres de module plus petit que $1$.

\noindent En appliquant la \og proc\'edure \textsc{Julia}\fg\, introduite pr\'ec\'edemment \`a $$\left(\lambda (x^2+xt-1),\frac{(t^2+xt-1)}{\lambda}\right),$$ qui d\'ecrit la dynamique de la restriction de $\Phi_\mathrm{A}$ \`a $\big\{\mathrm{M}\in\mathcal{M}(2;\mathbb{R})\,\big\vert\,\det\mathrm{M}=1\big\}$, nous obtenons par exemple les figures qui suivent pour $\rho=10$, $r=3$ et les param\`etres $\lambda=1$ et $1.5$.  

\smallskip

\begin{center}
\begin{tabular}{llll}
\begin{picture}(0,0)%
\includegraphics{dessin5.pstex}%
\end{picture}%
\setlength{\unitlength}{3947sp}%
\begingroup\makeatletter\ifx\SetFigFont\undefined%
\gdef\SetFigFont#1#2#3#4#5{%
  \reset@font\fontsize{#1}{#2pt}%
  \fontfamily{#3}\fontseries{#4}\fontshape{#5}%
  \selectfont}%
\fi\endgroup%
\begin{picture}(1050,1050)(1201,-2611)
\end{picture}%
\hspace{20mm}&\hspace{20mm}\begin{picture}(0,0)%
\includegraphics{dessin6.pstex}%
\end{picture}%
\setlength{\unitlength}{3947sp}%
\begingroup\makeatletter\ifx\SetFigFont\undefined%
\gdef\SetFigFont#1#2#3#4#5{%
  \reset@font\fontsize{#1}{#2pt}%
  \fontfamily{#3}\fontseries{#4}\fontshape{#5}%
  \selectfont}%
\fi\endgroup%
\begin{picture}(1050,1050)(1201,-5011)
\end{picture}%
\hspace{20mm}&\hspace{20mm}\begin{picture}(0,0)%
\includegraphics{dessin7.pstex}%
\end{picture}%
\setlength{\unitlength}{3947sp}%
\begingroup\makeatletter\ifx\SetFigFont\undefined%
\gdef\SetFigFont#1#2#3#4#5{%
  \reset@font\fontsize{#1}{#2pt}%
  \fontfamily{#3}\fontseries{#4}\fontshape{#5}%
  \selectfont}%
\fi\endgroup%
\begin{picture}(1050,1050)(1201,-1411)
\end{picture}%
\hspace{20mm}&\hspace{20mm}\begin{picture}(0,0)%
\includegraphics{dessin8.pstex}%
\end{picture}%
\setlength{\unitlength}{3947sp}%
\begingroup\makeatletter\ifx\SetFigFont\undefined%
\gdef\SetFigFont#1#2#3#4#5{%
  \reset@font\fontsize{#1}{#2pt}%
  \fontfamily{#3}\fontseries{#4}\fontshape{#5}%
  \selectfont}%
\fi\endgroup%
\begin{picture}(1050,1050)(1201,-2611)
\end{picture}%
\\
$\kappa=10$, $\lambda=1$&\hspace{20mm}$\kappa=75$, $\lambda=1$&\hspace{20mm}$\kappa=10$, $\lambda=1.5$&\hspace{20mm}$\kappa=75$, $\lambda=1.5$
\end{tabular}
\end{center}

\medskip

\noindent Notons qu'il y a une similitude certaine avec les figures pr\'ec\'edentes.

\subsubsection{Pr\'ecisions sur la dynamique dans le cas quaternionique}

\noindent Lorsque la matrice $\mathrm{A}$ est une matrice de quaternions nous pouvons pr\'eciser le discours pr\'ec\'edent.

\noindent L'it\'eration des transformations polynomiales du corps des quaternions a \'et\'e abord\'ee dans \cite{GV} o\`u les auteurs proposent une adaptation de la th\'eorie de \textsc{Fatou} \textsc{Julia}.

\noindent Consid\'erons une application $\Phi_\mathrm{A}\colon\mathrm{M}\mapsto\mathrm{A}\mathrm{M}^2$ avec cette fois $\mathrm{A}~=~\left[\begin{array}{cc} a & b\\ -\overline{b}& \overline{a}\end{array}\right]$ quaternionique. Quitte \`a faire agir une homoth\'etie r\'eelle $\mathrm{M}\mapsto\rho\mathrm{M}$ sur $\Phi_\mathrm{A}$ nous pouvons supposer que $\det\mathrm{A}=1$, \emph{i.e.} que $\mathrm{A}$ est un quaternion de module $1$. Mieux en conjuguant par une transformation $\mathrm{M}\mapsto\mathrm{B}\mathrm{M}\mathrm{B}^{-1}$ o\`u $\mathrm{B}$ est encore un quaternion ad-hoc nous pouvons supposer que $\mathrm{A}$ est diagonale, \emph{i.e.} $\mathrm{A}=\mathrm{diag}(\mathrm{e}^{\mathrm{i}\vartheta},\mathrm{e}^{-\mathrm{i}\vartheta})$. On constate donc, modulo ces modifications, que le corps $\Bbb H$ des quaternions $$\Bbb H=\left\{\left[\begin{array}{cc} x & y \\ -\overline{y} & \overline{x}\end{array}\right]\,\vert\, (x,y)\in\mathbb{C}^2\right\}$$ est compl\`etement invariant par $\Phi_\mathrm{A};$ il en est de m\^eme pour l'ensemble $\Bbb H_1$ des quaternions de module $1$ $$\Bbb H_1=\big\{\mathrm{M}\in\Bbb H\,\big\vert\,\det\mathrm{M}=1\big\}.$$ Notons que $\Bbb H_{<1}=\big\{\mathrm{M}\in\Bbb H\,\big\vert\,\vert\det\mathrm{M}\vert<1\big\}$ est contenu dans le bassin d'attraction $\mathrm{W}_\mathrm{A}^s({\bf 0})$ de l'origine. \'Etant donn\'e que $\Bbb H_1$ est invariant, ceci produit des orbites born\'ees non contenues dans $\partial\mathrm{W}_\mathrm{A}^s({\bf 0})$ mais dans son bord. Puisque les $\sigma_\mathrm{D}\colon\mathrm{M}\mapsto\mathrm{D}\mathrm{M}\mathrm{D}^{-1}$, avec $\mathrm{D}$ diagonale, commutent avec $\Phi_\mathrm{A}$ les ensembles $\sigma_\mathrm{D}(\Bbb H)$, $\sigma_\mathrm{D}(\Bbb H_1)$ sont invariants. Ceci donne deux autres exemples d'ensembles invariants non born\'es:
\begin{itemize}
\item l'ensemble $\mathcal{H}$ de dimension r\'eelle $5$ qui est l'union des $\sigma_\mathrm{D}(\Bbb H)$:
$$\mathcal{H}=\left\{\left[\begin{array}{cc} x& \xi y\\ -\frac{\overline{y}}{\xi}& \overline{x}\end{array}\right]\,\vert\, (x,y)\in\mathbb{C}^2,\,\xi\in\mathbb{C}^*\right\}=\left\{\left[\begin{array}{cc} x& y\\  \alpha\overline{y} &\overline{x} \end{array}\right]\,\vert\, (x,y)\in\mathbb{C}^2,\,\alpha\in\mathbb{R}_{<0}\right\};$$

\item ainsi que l'ensemble $\mathcal{H}_1$ de dimension r\'eelle $4$ qui est l'union des $\sigma_\mathrm{D}(\Bbb H_1)$:
$$\mathcal{H}_1=\left\{\left[\begin{array}{cc} x&\xi y \\ -\frac{\overline{y}}{\xi}& \overline{x}\end{array}\right]\,\vert\, (x,y)\in\mathbb{C}^2,\,\xi\in\mathbb{C}^*\right\}.$$
\end{itemize}

\noindent Chaque \'el\'ement de $\mathcal{H}_1$ a son orbite born\'ee contenue dans le bord du bassin d'attraction $\mathrm{W}_\mathrm{A}^s({\bf 0})$. Les matrices de la forme $\mathrm{diag}(\mathrm{e}^{\mathrm{i}\varphi},\mathrm{e}^{\mathrm{i}\vartheta})$ ont aussi une orbite born\'ee par $\Phi_\mathrm{A}$ et pour la plupart ne sont pas dans $\mathcal{H}_1$.

\noindent Nous avons r\'ealis\'e quelques exp\'eriences num\'eriques concernant la restriction de $\Phi_\mathrm{A}$ \`a $\Bbb H_1$. Avant de les pr\'esenter faisons quelques remarques \'el\'ementaires. L'\'etude de la restriction de $\Phi_\mathrm{A}$ \`a $\Bbb H$ se ram\`ene bien s\^ur \`a celle de l'application analytique r\'eelle de $\mathbb{C}^2\simeq\mathbb{R}^4$ dans lui-m\^eme (on garde la m\^eme notation) induite par les deux premi\`eres composantes: $$\Phi_\mathrm{A}\colon(x,y)\mapsto\big(\mathrm{e}^{\mathrm{i}\vartheta}(x^2-\vert y\vert^2),\mathrm{e}^{\mathrm{i}\vartheta}y(x+\overline{x})\big).$$ Sur $\Bbb H_1$ nous avons $\vert x\vert^2+\vert y\vert^2=1$ de sorte que $\Phi_{\mathrm{A}\vert_{\Bbb H_1}}$ est fibr\'ee $$\Phi_{\mathrm{A}\vert_{\Bbb H_1}}\colon(x,y)\mapsto\big(\mathrm{e}^{\mathrm{i}\vartheta}(x^2+\vert x\vert^2-1),\mathrm{e}^{\mathrm{i}\vartheta}y(x+\overline{x})\big)=\big(\varphi_\vartheta(x),\mathrm{e}^{\mathrm{i}\vartheta}y(x+\overline{x})\big),$$ \emph{i.e.}  la premi\`ere composante ne d\'epend que de $x$. 

\noindent Notons que la famille des cercles verticaux param\'etr\'ee par $\eta\mapsto(x, y\mathrm{e}^{\mathrm{i}\eta})$ est globalement invariante. La seconde composante de $\Phi^k_{\mathrm{A}\vert_{\Bbb H_1}}$ a pour argument $\varphi+k\vartheta$ avec $y=\rho\mathrm{e}^{\mathrm{i}\varphi}$ qui est aussi l'argument de l'it\'er\'e $k$-i\`eme de~$y$ suivant la rotation d'angle $\vartheta$. La premi\`ere composante indique comment $\Phi_\mathrm{A}$ fait passer d'un cercle \`a l'autre. Par exemple si $\underline{x}$ est un point p\'eriodique de p\'eriode $k$ de $\varphi_\vartheta$ alors l'application $\varphi_\vartheta^k$ est une rotation d'angle $k\vartheta$ sur le cercle $\mathcal{C}=\big\{(\underline{x},y)\in\Bbb H_1\big\}$ de $\Bbb H_1\simeq\mathbb{S}^3$. Nous allons nous int\'eresser \`a cette premi\`ere composante $\varphi_\vartheta$ que l'on identifie via $x=x_1+\mathrm{i}x_2$ \`a l'application toujours not\'ee $\varphi_\vartheta$: $$\varphi_\vartheta\colon(x_1,x_2)\mapsto(\cos\vartheta (2x_1^2-1)-2\sin\vartheta x_1x_2,2\cos\vartheta x_1x_2+\sin\vartheta(2x_1^2-1)).$$

\noindent Les applications $\varphi_\vartheta$ respectent toutes le disque unit\'e $\mathbb{D}(0,1)$ de $\mathbb{R}^2$ et c'est la dynamique dans ce disque qui nous int\'eresse. Dans toute la suite nous consid\'erons la restriction de $\varphi_\vartheta$ au disque ferm\'e $\overline{\mathbb{D}(0,1)}$ en gardant la m\^eme notation $\varphi_\vartheta$. On v\'erifie sans peine que tous les $\varphi_{\vartheta\vert_{\mathbb{S}^1}}\colon\mathbb{S}^1\to\mathbb{S}^1$ sont conjugu\'es \`a l'application $z\mapsto z^2$ du cercle dans lui-m\^eme. Le diam\`etre $[-1,1]$ va sur le diam\`etre $[-\mathrm{e}^{\mathrm{i}\vartheta},\mathrm{e}^{\mathrm{i}\vartheta}]$ et est en particulier invariant lorsque $\vartheta=0$. Le diam\`etre vertical $x_1=0$ est contract\'e sur le point $-\mathrm{e}^{\mathrm{i}\vartheta}$. Ainsi un point de type $(0,x_2)$ va \`a la premi\`ere it\'eration sur $-\mathrm{e}^{\mathrm{i}\vartheta}$ puis reste par it\'eration sur le cercle $\mathbb{S}^1$ o\`u, g\'en\'eriquement sur $\vartheta,$ son orbite est dense. Remarquons que les courbes \og alg\'ebriques\fg\, r\'eelles\, $(\varphi^k)^{-1}([-\mathrm{i},\mathrm{i}])$ pr\'esentent une propri\'et\'e analogue apr\`es $k+1$ it\'erations. La description de $\varphi_0\colon(x_1,x_2)\mapsto(2x_1^2-1,2x_1x_2)$  est relativement raisonnable. Sur le cercle unit\'e $\varphi_0$ co\"{i}ncide avec $x\mapsto x^2$ et est induite par $\mathrm{M}\mapsto\mathrm{M}^2$ restreinte \`a $\Bbb H_1$. On se souvient que $\mathcal{O}(\mathrm{M};\Phi_{\mathrm{Id}})$ est trac\'ee dans le $2$-plan engendr\'e par $\mathrm{M}$, $\mathrm{Id}$ et $-\mathrm{Id}$. En se restreignant \`a $\Bbb H_1$ nous en d\'eduisons que les orbites de $\Phi_{\mathrm{Id}\vert_{\Bbb H_1}}$ sont trac\'ees sur les cercles de $\mathbb{S}^3$ obtenus comme intersection du $2$-plan, r\'eel cette fois, passant par $\mathrm{M},$ $\mathrm{Id}$ et $-\mathrm{Id}$ (except\'e lorsque $\mathrm{M}$ est $\mathrm{Id}$ ou $-\mathrm{Id}$). Les orbites de $\varphi_0$ sont donc trac\'ees sur les projections (par $(x,y)\mapsto(x,0)$) de ces cercles; ce sont les ellipses passant par les points $(1,0)$ et $(-1,0)$ et tangentes au cercle unit\'e en $(1,0)$ et $(-1,0)$ :

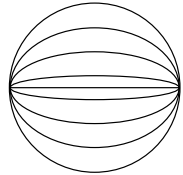
\begin{center}
\begin{picture}(0,0)%
\includegraphics{ellipses.pstex}%
\end{picture}%
\setlength{\unitlength}{3947sp}%
\begingroup\makeatletter\ifx\SetFigFont\undefined%
\gdef\SetFigFont#1#2#3#4#5{%
  \reset@font\fontsize{#1}{#2pt}%
  \fontfamily{#3}\fontseries{#4}\fontshape{#5}%
  \selectfont}%
\fi\endgroup%
\begin{picture}(1076,1074)(6063,-5098)
\end{picture}%

\end{center}

\noindent On v\'erifie facilement que la restriction de $\varphi_0$ \`a chaque ellipse est conjugu\'ee \`a $x\mapsto x^2,$ sauf \'evidemment dans le cas sp\'ecial o\`u l'ellipse d\'eg\'en\`ere sur l'intervalle $[-1,1]$. Dans ce cas on constate que $\varphi_{0\vert_{[-1,1]}}$ est l'application $x_1\mapsto 2x_1^2-1$ qui n'est rien d'autre que la c\'el\`ebre application \og logistique\fg\, \`a conjugaison pr\`es (\cite{CoEc}). En fait comme nous l'avons dit l'application $\Phi_{\mathrm{Id}}$ est compatible \`a la conjugaison. Sur $\mathbb{S}^3\simeq\Bbb H^1$ les classes de conjugaison sont d\'etermin\'ees par la trace $x+\overline{x}=2x_1$ et sont donc des $2$-sph\`eres $(x_1=$ cte $)\cap\,\mathbb{S}^3$. Par exemple si~$\mathrm{M}\in\Bbb H_1$ est p\'eriodique pour $\varphi_{\mathrm{Id}}$ alors toute la sph\`ere $\mathbb{S}^2=$ classe de $\mathrm{M}$ est form\'ee d'\'el\'ements p\'eriodiques ce que l'on peut voir directement sur l'application $\varphi_0$; celle-ci respecte les droites verticales dans leur ensemble: 

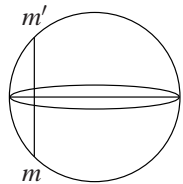
\begin{center}
\begin{picture}(0,0)%
\includegraphics{dtevert.pstex}%
\end{picture}%
\setlength{\unitlength}{3947sp}%
\begingroup\makeatletter\ifx\SetFigFont\undefined%
\gdef\SetFigFont#1#2#3#4#5{%
  \reset@font\fontsize{#1}{#2pt}%
  \fontfamily{#3}\fontseries{#4}\fontshape{#5}%
  \selectfont}%
\fi\endgroup%
\begin{picture}(1076,1156)(4863,-2285)
\put(4951,-2236){\makebox(0,0)[lb]{\smash{{\SetFigFont{12}{14.4}{\familydefault}{\mddefault}{\updefault}\begin{small}$m$\end{small}}}}}
\put(4951,-1261){\makebox(0,0)[lb]{\smash{{\SetFigFont{12}{14.4}{\familydefault}{\mddefault}{\updefault}\begin{small}$m'$\end{small}}}}}
\end{picture}%

\end{center}

 \noindent On peut d\'eduire d'ailleurs de fa\c{c}on directe la dynamique de l'application logistique (qui se fait usuellement de fa\c{c}on combinatoire via un codage RLC, \emph{voir} \cite{CoEc}) de la dynamique de l'application du cercle $x\mapsto x^2$ (qui elle se fait en utilisant l'\'ecriture diadique des angles). Si on s'int\'eresse aux points p\'eriodiques de p\'eriode~$2$ de $\varphi_0$ on trouve (outre le point $1\sim (1,0)$) le segment $[\mathrm{j},\mathrm{j}^2]$. La restriction de $\varphi_0$ \`a $[\mathrm{j},\mathrm{j}^2]$ est donn\'ee par $\left(-\frac{1}{2},x_2\right)\mapsto\left(-\frac{1}{2},-x_2\right)$ qui est \'evidemment une involution. Son point fixe $\left(-\frac{1}{2},0\right)$ correspond au point fixe~$-\frac{1}{2}$ de l'application logistique. Remarquons que l'ensemble des points p\'eriodiques de $\varphi_0$ est une union dense de segments verticaux $[m,m']$, les segments passant par les points p\'eriodiques de $x\mapsto x^2$. Comme $\varphi_\pi$ est conjugu\'e \`a $\varphi_0$ par $(x_1,x_2)\mapsto(-x_1,-x_2)$ tout ce que nous avons dit pour $\varphi_0$ s'adapte via la conjugaison \`a $\varphi_\pi.$ Nous conjecturons que pour les applications $\varphi_\vartheta$ avec $\vartheta$ g\'en\'erique les points p\'eriodiques de p\'eriode donn\'ee sont en nombre fini \`a l'inverse de $\varphi_0$. Par exemple l'application $$\varphi_{\pi/2}\colon(x_1,x_2)\mapsto(-2x_1x_2,2x_1^2-1)$$ a les points fixes suivants $(0,-1)$, $\left(\frac{1}{2},-\frac{1}{2}\right)$, $\left(-\frac{1}{2},-\frac{1}{2}\right)$ et ses points $2$-p\'eriodiques sont en nombre fini contrairement \`a ceux de $\varphi_0$. En effet $\varphi_{\pi/2}^2\colon(x_1,x_2)\mapsto(4x_1x_2(2x_1^2-1),8x_1^2x_2^2-1)$ a pour points fixes les trois points qui pr\'ec\`edent auxquels s'ajoutent $\left(\frac{\sqrt{3}}{2},-\frac{1}{2}\right)$ $\left(-\frac{\sqrt{3}}{2},\frac{1}{2}\right)$ qui sont les deux racines $3$i\`eme de l'unit\'e $\mathrm{j}$ et $\mathrm{j}^2$.

\noindent Nous avons essay\'e de pr\'eciser les points $2$-p\'eriodiques et leurs bifurcations en utilisant la technique des bases de Grobner qui fonctionne bien lorsque $$\Delta(\vartheta)=(4\cos^2\vartheta-3)(2\cos^3\vartheta - 3\cos^2\vartheta+2) (2\cos^3\vartheta + 3 \cos^2\vartheta-2) (\cos\vartheta - 1) (\cos\vartheta + 1)$$ est non nulle. Sous cette hypoth\`ese ($\Delta(\vartheta)\not=0$) l'application $\varphi_\vartheta$ a sept points p\'eriodiques de p\'eriode $2$, tous r\'eels, \emph{i.e.} dans $\mathbb{R}^2,$ qui sont les suivants: 
\begin{itemize}
\item[$\mathfrak{a}.$] le point fixe $\mathrm{e}^{-\mathrm{i}\vartheta}\simeq(\cos\vartheta,-\sin\vartheta)$;

\item[$\mathfrak{b}.$] les deux points fixes 
\begin{align*}
&\mathfrak{b.1}.\,\,\,\left( -\frac{1}{2},\frac{\cos\vartheta-1}{2\sin\vartheta}\right) && \text{ et }
&&\mathfrak{b.2}.\,\,\,\left(\frac{1}{2}, -\frac{\cos\vartheta+1}{2\sin\vartheta}\right);
\end{align*}

\item[$\mathfrak{c}.$] les deux points $2$-p\'eriodiques de la restriction $\varphi_{\vartheta_{\vert_{\mathbb{S}^1}}}\colon z \mapsto \mathrm{e}^{\mathrm{i}\vartheta}z^2$ qui sont $\mathrm{j}\mathrm{e}^{-\mathrm{i}\vartheta}$ et $\mathrm{j}^2\mathrm{e}^{-\mathrm{i}\vartheta}$;

\item[$\mathfrak{d}.$] deux autres points $2$-p\'eriodiques donn\'es par
\begin{align*}
&\mathfrak{d.1}.\,\,\,\left(\frac{\sqrt{1+4\cos^2\vartheta}-1}{4\cos\vartheta}, \frac{\sin\vartheta(\sqrt{1+4\cos^2\vartheta}+1)}{4\cos^2\vartheta}\right), 
&&\mathfrak{d.2}.\,\,\,\left(-\frac{\sqrt{1+4\cos^2\vartheta}+1}{4\cos\vartheta}, \frac{\sin\vartheta(1-\sqrt{1+4\cos^2\vartheta})}{4\cos^2\vartheta}\right).
\end{align*}
\end{itemize}

\noindent Remarquons que lorsque $\vartheta$ tend vers $0$, les points fixes de type $\mathfrak{b}.$ tendent respectivement vers  $\left(-\frac{1}{2},0\right)$ qui est point fixe de $\varphi_0$ et vers \og $\left(\frac{1}{2},\infty\right)$ \fg. Lorsque $\vartheta$ tend vers $\pi$ le ph\'enom\`ene inverse se produit: un point fixe tend vers $\left(\frac{1}{2},0\right)$ (point fixe de $\varphi_\pi$) et l'autre s'\'echappe.

\noindent Un point de type $\mathfrak{b.1}.$, resp. $\mathfrak{b.2}.$ est dans le disque unit\'e ferm\'e si et seulement si $-\frac{2\pi}{3}\leq\vartheta\leq\frac{2\pi}{3}$, resp. si et seulement si $\frac{\pi}{3}\leq\vartheta\leq \frac{5\pi}{3}$. Notons que les points fixes de $\varphi_\vartheta$ ne sont jamais contractants. Les points $2$-p\'eriodiques de type $\mathfrak{d}.$ sont dans le disque unit\'e ferm\'e si et seulement si $-\frac{\pi}{4}\leq\vartheta\leq \frac{\pi}{4}$ ou $\frac{3\pi}{4}\leq\vartheta\leq\frac{5\pi}{4}$ et ils le sont simultan\'ement.

\bigskip

\noindent Nous avons appliqu\'e la proc\'edure type \textsc{Julia} \'evoqu\'ee pr\'ec\'edemment aux applications $\varphi_\vartheta$ dans le disque unit\'e en choisissant une donn\'ee de contr\^ole $\rho$ plus petite que $1$ mais  tr\`es voisine de $1$. Les figures suivantes mesurent donc la fa\c{c}on dont les orbites s'approchent du bord du disque et \`a quelle vitesse. Par exemple dans la premi\`ere figure o\`u nous appliquons la proc\'edure \`a $\varphi_0$, les points p\'eriodiques non situ\'es sur le bord ne s'en approchent pas ce qui produit les lignes verticales \og indigo $\sim$ rouge\fg. Remarquons dans ce cas sp\'ecial la liaison avec le probl\`eme suivant: pour la transformation du cercle $z\mapsto z^2$ \'etudier comment les orbites approchent le point fixe $1$ et \`a quelle vitesse. Ceci explique la structure cantorique que l'on observe sur cette premi\`ere figure. 

\noindent Nous avons choisi de pr\'esenter des figures pour $\kappa$ petit pour la qualit\'e visuelle; on constate exp\'erimentalement qu'elles ne varient pas qualitativement quand $\kappa$ augmente. Le param\`etre $r$, resp. $\rho$ prend la valeur $1$, resp. $\sqrt{0.99}$. 

\medskip

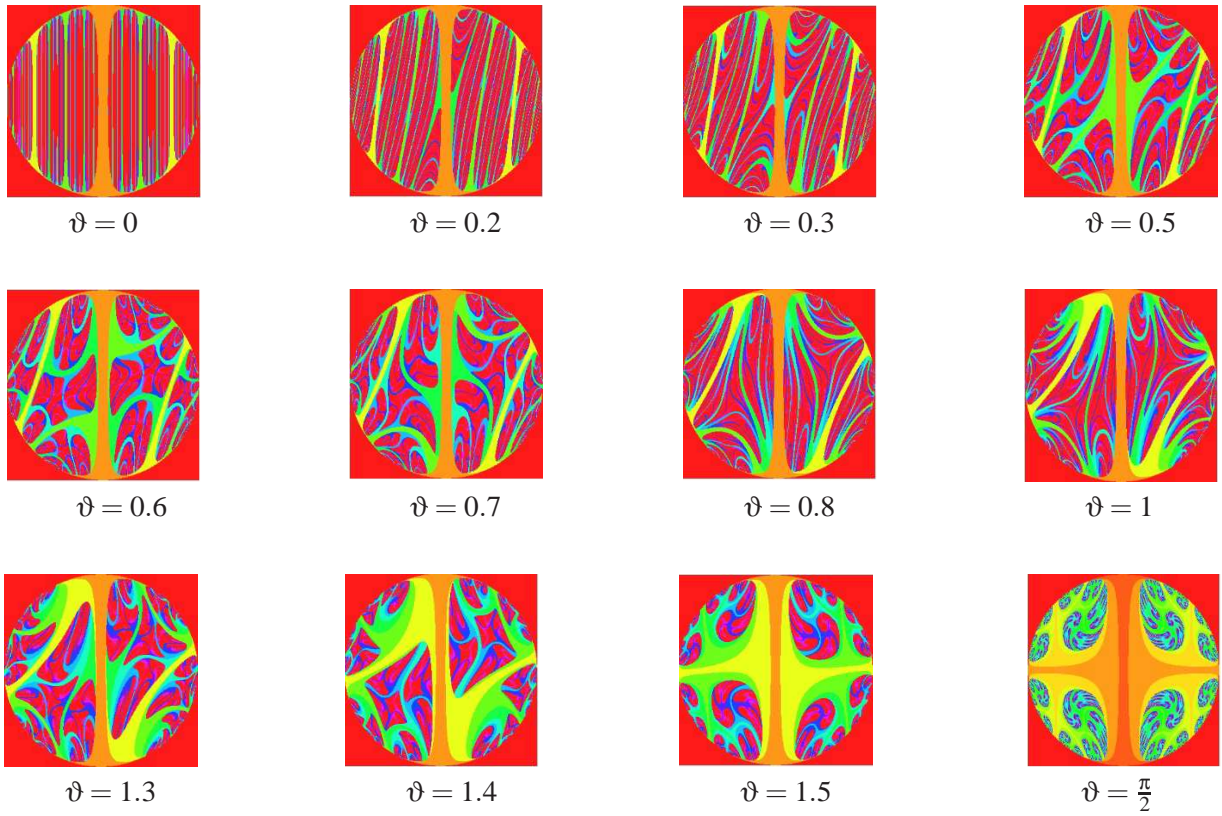
\begin{center}
\begin{tabular}{llll}
\begin{picture}(0,0)%
\includegraphics{theta0b.pstex}%
\end{picture}%
\setlength{\unitlength}{3947sp}%
\begingroup\makeatletter\ifx\SetFigFont\undefined%
\gdef\SetFigFont#1#2#3#4#5{%
  \reset@font\fontsize{#1}{#2pt}%
  \fontfamily{#3}\fontseries{#4}\fontshape{#5}%
  \selectfont}%
\fi\endgroup%
\begin{picture}(1200,1200)(1201,-5161)
\end{picture}%
 \hspace{15mm}& \hspace{15mm}\begin{picture}(0,0)%
\includegraphics{theta0.2b.pstex}%
\end{picture}%
\setlength{\unitlength}{3947sp}%
\begingroup\makeatletter\ifx\SetFigFont\undefined%
\gdef\SetFigFont#1#2#3#4#5{%
  \reset@font\fontsize{#1}{#2pt}%
  \fontfamily{#3}\fontseries{#4}\fontshape{#5}%
  \selectfont}%
\fi\endgroup%
\begin{picture}(1200,1200)(1201,-5161)
\end{picture}%
\hspace{15mm}& \hspace{15mm}\begin{picture}(0,0)%
\includegraphics{theta0.3b.pstex}%
\end{picture}%
\setlength{\unitlength}{3947sp}%
\begingroup\makeatletter\ifx\SetFigFont\undefined%
\gdef\SetFigFont#1#2#3#4#5{%
  \reset@font\fontsize{#1}{#2pt}%
  \fontfamily{#3}\fontseries{#4}\fontshape{#5}%
  \selectfont}%
\fi\endgroup%
\begin{picture}(1200,1200)(1201,-5161)
\end{picture}%
 \hspace{15mm}& \hspace{15mm}\begin{picture}(0,0)%
\includegraphics{theta0.5b.pstex}%
\end{picture}%
\setlength{\unitlength}{3947sp}%
\begingroup\makeatletter\ifx\SetFigFont\undefined%
\gdef\SetFigFont#1#2#3#4#5{%
  \reset@font\fontsize{#1}{#2pt}%
  \fontfamily{#3}\fontseries{#4}\fontshape{#5}%
  \selectfont}%
\fi\endgroup%
\begin{picture}(1200,1200)(1201,-5161)
\end{picture}%
\\
\hspace{8mm}$\vartheta=0$ &\hspace{23mm}$\vartheta=0.2$&\hspace{23mm}$\vartheta=0.3$&\hspace{23mm}$\vartheta=0.5$
\end{tabular}
\end{center}

\medskip

\begin{center}
\begin{tabular}{llll}
\begin{picture}(0,0)%
\includegraphics{theta0.6b.pstex}%
\end{picture}%
\setlength{\unitlength}{3947sp}%
\begingroup\makeatletter\ifx\SetFigFont\undefined%
\gdef\SetFigFont#1#2#3#4#5{%
  \reset@font\fontsize{#1}{#2pt}%
  \fontfamily{#3}\fontseries{#4}\fontshape{#5}%
  \selectfont}%
\fi\endgroup%
\begin{picture}(1200,1200)(1201,-5161)
\end{picture}%
 \hspace{15mm}& \hspace{15mm}\begin{picture}(0,0)%
\includegraphics{theta0.7b.pstex}%
\end{picture}%
\setlength{\unitlength}{3947sp}%
\begingroup\makeatletter\ifx\SetFigFont\undefined%
\gdef\SetFigFont#1#2#3#4#5{%
  \reset@font\fontsize{#1}{#2pt}%
  \fontfamily{#3}\fontseries{#4}\fontshape{#5}%
  \selectfont}%
\fi\endgroup%
\begin{picture}(1200,1200)(1201,-5161)
\end{picture}%
\hspace{15mm}& \hspace{15mm}\begin{picture}(0,0)%
\includegraphics{theta0.8b.pstex}%
\end{picture}%
\setlength{\unitlength}{3947sp}%
\begingroup\makeatletter\ifx\SetFigFont\undefined%
\gdef\SetFigFont#1#2#3#4#5{%
  \reset@font\fontsize{#1}{#2pt}%
  \fontfamily{#3}\fontseries{#4}\fontshape{#5}%
  \selectfont}%
\fi\endgroup%
\begin{picture}(1200,1200)(1201,-5161)
\end{picture}%
 \hspace{15mm}& \hspace{15mm}\begin{picture}(0,0)%
\includegraphics{theta1b.pstex}%
\end{picture}%
\setlength{\unitlength}{3947sp}%
\begingroup\makeatletter\ifx\SetFigFont\undefined%
\gdef\SetFigFont#1#2#3#4#5{%
  \reset@font\fontsize{#1}{#2pt}%
  \fontfamily{#3}\fontseries{#4}\fontshape{#5}%
  \selectfont}%
\fi\endgroup%
\begin{picture}(1200,1200)(1201,-5161)
\end{picture}%
\\
\hspace{8mm} $\vartheta=0.6$&\hspace{23mm}$\vartheta=0.7$&\hspace{23mm}$\vartheta=0.8$&\hspace{23mm}$\vartheta=1$
\end{tabular}
\end{center}

\medskip

\begin{center}
\begin{tabular}{llll}
\begin{picture}(0,0)%
\includegraphics{theta1.3b.pstex}%
\end{picture}%
\setlength{\unitlength}{3947sp}%
\begingroup\makeatletter\ifx\SetFigFont\undefined%
\gdef\SetFigFont#1#2#3#4#5{%
  \reset@font\fontsize{#1}{#2pt}%
  \fontfamily{#3}\fontseries{#4}\fontshape{#5}%
  \selectfont}%
\fi\endgroup%
\begin{picture}(1200,1200)(1201,-5161)
\end{picture}%
 \hspace{15mm}& \hspace{15mm}\begin{picture}(0,0)%
\includegraphics{theta1.4b.pstex}%
\end{picture}%
\setlength{\unitlength}{3947sp}%
\begingroup\makeatletter\ifx\SetFigFont\undefined%
\gdef\SetFigFont#1#2#3#4#5{%
  \reset@font\fontsize{#1}{#2pt}%
  \fontfamily{#3}\fontseries{#4}\fontshape{#5}%
  \selectfont}%
\fi\endgroup%
\begin{picture}(1200,1200)(1201,-5161)
\end{picture}%
\hspace{15mm}& \hspace{15mm}\begin{picture}(0,0)%
\includegraphics{theta1.5b.pstex}%
\end{picture}%
\setlength{\unitlength}{3947sp}%
\begingroup\makeatletter\ifx\SetFigFont\undefined%
\gdef\SetFigFont#1#2#3#4#5{%
  \reset@font\fontsize{#1}{#2pt}%
  \fontfamily{#3}\fontseries{#4}\fontshape{#5}%
  \selectfont}%
\fi\endgroup%
\begin{picture}(1200,1200)(1201,-5161)
\end{picture}%
 \hspace{15mm}& \hspace{15mm}  \begin{picture}(0,0)%
\includegraphics{thetapi2b.pstex}%
\end{picture}%
\setlength{\unitlength}{3947sp}%
\begingroup\makeatletter\ifx\SetFigFont\undefined%
\gdef\SetFigFont#1#2#3#4#5{%
  \reset@font\fontsize{#1}{#2pt}%
  \fontfamily{#3}\fontseries{#4}\fontshape{#5}%
  \selectfont}%
\fi\endgroup%
\begin{picture}(1200,1200)(1201,-3961)
\end{picture}%
\\
\hspace{8mm}$\vartheta=1.3$&\hspace{23mm}$\vartheta=1.4$&\hspace{23mm}$\vartheta=1.5$&\hspace{23mm}$\vartheta=\frac{\pi}{2}$
\end{tabular}
\end{center}

\medskip

\noindent Le trac\'e ponctuel des orbites ne s'av\`ere pas tr\`es probant. Pour pallier ce d\'efaut nous avons it\'er\'e des segments verticaux par les applications $\varphi_\vartheta$. Nous pr\'esentons ci-dessous les onze premiers it\'er\'es de la droite $\mathcal{D}$ d'\'equation $x_1=0.6$ intersect\'ee avec le disque de rayon $1$ par l'application~$\varphi_{\pi/2}$.

\medskip

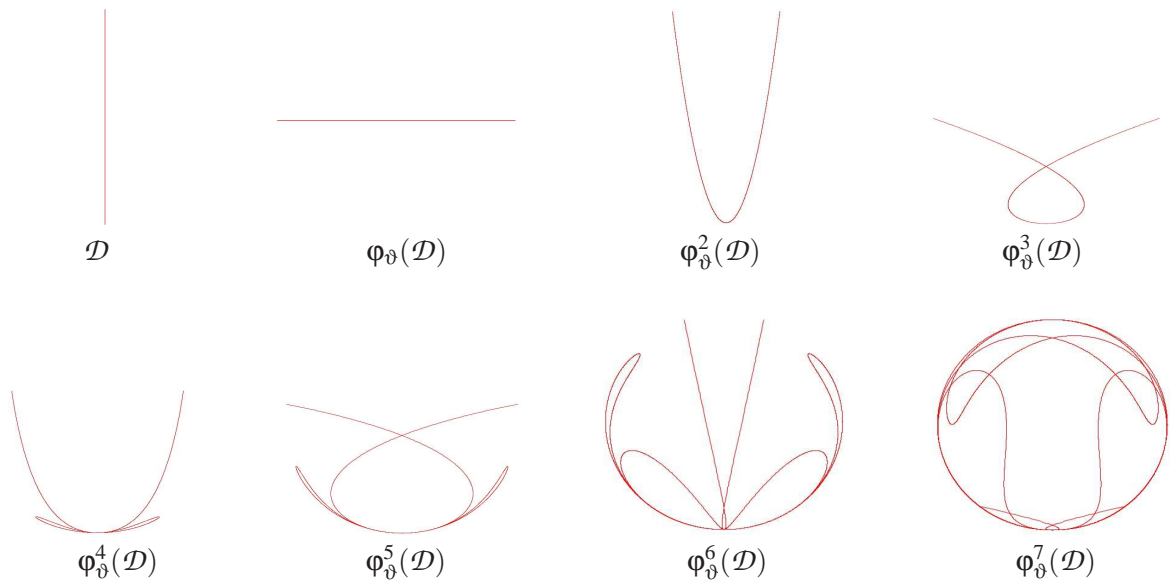
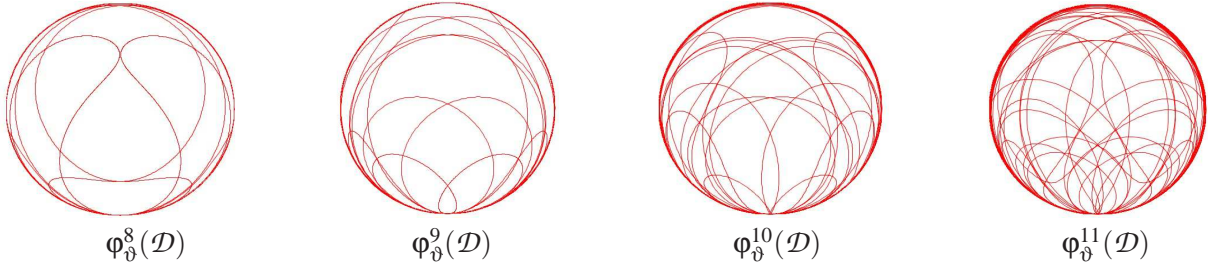
\begin{center}
\begin{tabular}{llll}
\hspace{-3mm}\begin{picture}(0,0)%
\includegraphics{itere.pstex}%
\end{picture}%
\setlength{\unitlength}{3947sp}%
\begingroup\makeatletter\ifx\SetFigFont\undefined%
\gdef\SetFigFont#1#2#3#4#5{%
  \reset@font\fontsize{#1}{#2pt}%
  \fontfamily{#3}\fontseries{#4}\fontshape{#5}%
  \selectfont}%
\fi\endgroup%
\begin{picture}(1500,1350)(1201,-1711)
\end{picture}%
 \hspace{1mm}& \hspace{2mm}\begin{picture}(0,0)%
\includegraphics{itere1.pstex}%
\end{picture}%
\setlength{\unitlength}{3947sp}%
\begingroup\makeatletter\ifx\SetFigFont\undefined%
\gdef\SetFigFont#1#2#3#4#5{%
  \reset@font\fontsize{#1}{#2pt}%
  \fontfamily{#3}\fontseries{#4}\fontshape{#5}%
  \selectfont}%
\fi\endgroup%
\begin{picture}(1500,1350)(1201,-6511)
\end{picture}%
\hspace{8mm}& \hspace{16mm}\begin{picture}(0,0)%
\includegraphics{itere2.pstex}%
\end{picture}%
\setlength{\unitlength}{3947sp}%
\begingroup\makeatletter\ifx\SetFigFont\undefined%
\gdef\SetFigFont#1#2#3#4#5{%
  \reset@font\fontsize{#1}{#2pt}%
  \fontfamily{#3}\fontseries{#4}\fontshape{#5}%
  \selectfont}%
\fi\endgroup%
\begin{picture}(750,1350)(1276,-6511)
\end{picture}%
 \hspace{1mm}& \hspace{14mm}\begin{picture}(0,0)%
\includegraphics{itere3.pstex}%
\end{picture}%
\setlength{\unitlength}{3947sp}%
\begingroup\makeatletter\ifx\SetFigFont\undefined%
\gdef\SetFigFont#1#2#3#4#5{%
  \reset@font\fontsize{#1}{#2pt}%
  \fontfamily{#3}\fontseries{#4}\fontshape{#5}%
  \selectfont}%
\fi\endgroup%
\begin{picture}(1500,675)(1201,-5236)
\end{picture}%
\\
\hspace{10mm}$\mathcal{D}$&\hspace{14mm}$\varphi_\vartheta(\mathcal{D})$&\hspace{18mm}$\varphi_\vartheta^2(\mathcal{D})$&\hspace{24mm}$\varphi_\vartheta^3(\mathcal{D})$
\end{tabular}
\end{center}

\medskip

\begin{center}
\begin{tabular}{llll}
\begin{picture}(0,0)%
\includegraphics{itere4.pstex}%
\end{picture}%
\setlength{\unitlength}{3947sp}%
\begingroup\makeatletter\ifx\SetFigFont\undefined%
\gdef\SetFigFont#1#2#3#4#5{%
  \reset@font\fontsize{#1}{#2pt}%
  \fontfamily{#3}\fontseries{#4}\fontshape{#5}%
  \selectfont}%
\fi\endgroup%
\begin{picture}(1200,900)(1201,-6061)
\end{picture}%
 \hspace{7mm}& \hspace{7mm}\begin{picture}(0,0)%
\includegraphics{itere5.pstex}%
\end{picture}%
\setlength{\unitlength}{3947sp}%
\begingroup\makeatletter\ifx\SetFigFont\undefined%
\gdef\SetFigFont#1#2#3#4#5{%
  \reset@font\fontsize{#1}{#2pt}%
  \fontfamily{#3}\fontseries{#4}\fontshape{#5}%
  \selectfont}%
\fi\endgroup%
\begin{picture}(1500,825)(1201,-5986)
\end{picture}%
\hspace{7mm}& \hspace{7mm}\begin{picture}(0,0)%
\includegraphics{itere6.pstex}%
\end{picture}%
\setlength{\unitlength}{3947sp}%
\begingroup\makeatletter\ifx\SetFigFont\undefined%
\gdef\SetFigFont#1#2#3#4#5{%
  \reset@font\fontsize{#1}{#2pt}%
  \fontfamily{#3}\fontseries{#4}\fontshape{#5}%
  \selectfont}%
\fi\endgroup%
\begin{picture}(1500,1350)(1201,-4111)
\end{picture}%
 \hspace{7mm}& \hspace{7mm}\begin{picture}(0,0)%
\includegraphics{itere7.pstex}%
\end{picture}%
\setlength{\unitlength}{3947sp}%
\begingroup\makeatletter\ifx\SetFigFont\undefined%
\gdef\SetFigFont#1#2#3#4#5{%
  \reset@font\fontsize{#1}{#2pt}%
  \fontfamily{#3}\fontseries{#4}\fontshape{#5}%
  \selectfont}%
\fi\endgroup%
\begin{picture}(1500,1350)(1201,-4111)
\end{picture}%
\\
\hspace{10mm}$\varphi_\vartheta^4(\mathcal{D})$&\hspace{18mm}$\varphi_\vartheta^5(\mathcal{D})$&\hspace{18mm}$\varphi_\vartheta^6(\mathcal{D})$&\hspace{18mm}$\varphi_\vartheta^7(\mathcal{D})$
\end{tabular}
\end{center}

\medskip

\begin{center}
\begin{tabular}{llll}
\begin{picture}(0,0)%
\includegraphics{itere8.pstex}%
\end{picture}%
\setlength{\unitlength}{3947sp}%
\begingroup\makeatletter\ifx\SetFigFont\undefined%
\gdef\SetFigFont#1#2#3#4#5{%
  \reset@font\fontsize{#1}{#2pt}%
  \fontfamily{#3}\fontseries{#4}\fontshape{#5}%
  \selectfont}%
\fi\endgroup%
\begin{picture}(1500,1350)(1201,-4111)
\end{picture}%
 \hspace{7mm}& \hspace{7mm}\begin{picture}(0,0)%
\includegraphics{itere9.pstex}%
\end{picture}%
\setlength{\unitlength}{3947sp}%
\begingroup\makeatletter\ifx\SetFigFont\undefined%
\gdef\SetFigFont#1#2#3#4#5{%
  \reset@font\fontsize{#1}{#2pt}%
  \fontfamily{#3}\fontseries{#4}\fontshape{#5}%
  \selectfont}%
\fi\endgroup%
\begin{picture}(1500,1350)(1201,-4111)
\end{picture}%
\hspace{7mm}& \hspace{7mm}\begin{picture}(0,0)%
\includegraphics{itere10.pstex}%
\end{picture}%
\setlength{\unitlength}{3947sp}%
\begingroup\makeatletter\ifx\SetFigFont\undefined%
\gdef\SetFigFont#1#2#3#4#5{%
  \reset@font\fontsize{#1}{#2pt}%
  \fontfamily{#3}\fontseries{#4}\fontshape{#5}%
  \selectfont}%
\fi\endgroup%
\begin{picture}(1500,1350)(1201,-6511)
\end{picture}%
 \hspace{7mm}& \hspace{7mm}\begin{picture}(0,0)%
\includegraphics{itere11.pstex}%
\end{picture}%
\setlength{\unitlength}{3947sp}%
\begingroup\makeatletter\ifx\SetFigFont\undefined%
\gdef\SetFigFont#1#2#3#4#5{%
  \reset@font\fontsize{#1}{#2pt}%
  \fontfamily{#3}\fontseries{#4}\fontshape{#5}%
  \selectfont}%
\fi\endgroup%
\begin{picture}(1500,1350)(1201,-4111)
\end{picture}%
 \\
\hspace{14mm}$\varphi_\vartheta^8(\mathcal{D})$&\hspace{18mm}$\varphi_\vartheta^9(\mathcal{D})$&\hspace{18mm}$\varphi_\vartheta^{10}(\mathcal{D})$&\hspace{18mm}$\varphi_\vartheta^{11}(\mathcal{D})$
\end{tabular}
\end{center}

\bigskip

\noindent On constate que rapidement les courbes images s'accumulent sur le cercle du bord tout entier.

\begin{probs}
\begin{itemize}
\item[1)] Pour $\mathrm{A}= \mathrm{diag}\left(\lambda,\frac{1}{\lambda}\right)$ caract\'eriser en fonction de $\lambda$ les orbites qui sont born\'ees en particulier celles qui ne sont pas dans $\mathrm{W}_\mathrm{A}^s({\bf 0})$.

\item[2)] Caract\'eriser le bassin d'attraction $\mathrm{W}_\mathrm{A}^s({\bf 0})$ de l'origine et si possible d\'ecrire son bord.

\item[3)] Donner la description pr\'ecise de la dynamique des $\Phi_\mathrm{A}$ dans le cas quaternionique.
\end{itemize}
\end{probs}

\subsubsection{Centralisateur}\label{centrdiag}

\noindent Comme nous l'avons fait au paragraphe pr\'ec\'edent nous allons d\'eterminer le groupe $\mathrm{Aut}(\mathcal{M}(2;\mathbb{C});\Phi_\mathrm{A})$ pour $\mathrm{A}= \mathrm{diag}\left(\lambda,\frac{1}{\lambda}\right)$, $\lambda^2\not=1$.

\begin{pro}
Soit $\mathrm{A}$ une matrice de la forme $\mathrm{diag}\left(\lambda,\frac{1}{\lambda} \right)$ avec $\lambda^2\not=1.$ Le groupe $\mathrm{Aut}(\mathcal{M}(2;\mathbb{C});\Phi_\mathrm{A})$ est engendr\'e par les $\sigma_\mathrm{P}$ avec $\mathrm{P}$ diagonale; en fait $\mathrm{Aut}(\mathcal{M}(2;\mathbb{C});\Phi_\mathrm{A})$ s'identifie \`a $\mathbb{C}^*$ agissant sur $\mathcal{M}(2;\mathbb{C})$ de la fa\c{c}on suivante: $(x,y,z,t,\alpha)\mapsto\left(x,\alpha y,\frac{z}{\alpha},t\right).$ Les orbites de cette action sont aussi celles du champ de vecteurs invariant $y\frac{\partial}{\partial y}-z\frac{\partial}{\partial z}$.
\end{pro}

\begin{proof}[D\'emonstration]
\noindent Avec des arguments analogues \`a ceux utilis\'es dans la d\'emonstration de la Proposition~\ref{centautphi2} on montre qu'un \'el\'ement $\varphi$ de $\mathrm{Aut}(\mathcal{M} (2;\mathbb{C});\Phi_\mathrm{A})$ est n\'ecessairement lin\'eaire.

\noindent \'Ecrivons $\varphi$ sous la forme $(\ell_1,\ell_2,\ell_3,\ell_4)$ les $\ell_i$ d\'esignant des formes lin\'eaires. La fibration $y/z=$ cte est invariante par~$\Phi_\mathrm{A},$ plus pr\'ecis\'ement nous avons $\frac{y}{z}\circ\Phi_\mathrm{A}=\lambda^2\frac{y}{z}$. Nous en d\'eduisons l'\'egalit\'e $\ell_3(\ell_2\circ\Phi_\mathrm{A})=\lambda^2(\ell_3\circ\Phi_\mathrm{A})\ell_2$. Ceci implique, puisque les seuls $3$-plans invariants par $\Phi_\mathrm{A}$ sont $y=0$ et $z=0$, l'alternative suivante:
\begin{align*}
& \text{ ou bien }\ell_2=\alpha y,\, \ell_3=\beta z, && \text{ ou bien }\ell_2=\alpha z,\,\ell_3=\beta y.
\end{align*}

\noindent En r\'e\'ecrivant $\frac{\ell_2}{\ell_3}\circ\Phi_\mathrm{A}=\lambda^2\frac{\ell_2} {\ell_3}$ on constate que la seconde \'eventualit\'e n'arrive que si $\lambda^4=1$. 

\noindent Dans le premier cas la commutation de $\varphi$ et $\Phi_\mathrm{A}$ entra\^ine que $\varphi=\left(x,\alpha y, \frac{z}{\alpha},t\right)$. Lorsque $\lambda^4=1$, $\lambda$ vaut $\mathrm{i}$ ou $-\mathrm{i}$ (les valeurs propres de $\mathrm{A}$ sont suppos\'ees distinctes). \'Ecrivons $\ell_1$ (resp. $\ell_4$) sous la forme $a_1x+b_1y+c_1z+d_1t$ (resp. $a_4x+b_4y+c_4z+d_4t$). L'\'egalit\'e $\Phi_\mathrm{A}\varphi=\varphi\,\Phi_\mathrm{A}$ conduit \`a $$b_1=c_1=b_4=c_4=a_1d_1=a_4d_4=0$$ et 
\begin{align*}
& a_4=-1-a_1, && d_4=-1-d_1, &&  a_1^2=a_1, && d_1^2 \lambda^2=d_1, && \lambda\alpha\beta=\frac{d_1} {\lambda}+a_1, && a_4^2=\lambda^2a_4, && a_4+\frac{d_4}{\lambda}=\frac{\alpha\beta}{\lambda}.
\end{align*}

\noindent Un calcul montre que si $a_1$ est nul, alors $\lambda=-1$, sinon $d_1=0$ et $\lambda=1$. Ces deux cas sont exclus par l'hypoth\`ese $\lambda^2\not=1$. 
\end{proof}

\subsection{Cas non diagonalisable}

\noindent Consid\'erons les applications de la forme $\mathrm{M}\mapsto\mathrm{A}\mathrm{M}^2$ avec $\mathrm{A}$ inversible non diagonalisable. Nous nous ramenons apr\`es conjugaison ad-hoc \`a $\mathrm{A}=\left[\begin{array}{cc}1 & 1 \\ 0 & 1\end{array}\right];$ par suite $$\Phi_\mathrm{A}\left(\left[\begin{array}{cc}x & y \\ z &t\end{array} \right]\right)=\left[\begin{array}{cc}(x^2+yz)+z(x+t) & y(x+t)+(t^2+yz)\\ z(x+t)& t^2+yz\end{array}\right].$$ La quadrique de dimension $2$ form\'ee des matrices nilpotentes est toujours envoy\'ee sur ${\bf 0}$ par $\Phi_\mathrm{A}$ et $\mathfrak{sl}(2;\mathbb{C})$ est encore contract\'ee, cette fois sur $\mathbb{C}\mathrm{A}$. La fibration $\frac{x-t}{z}=$ cte est invariante par~$\Phi_\mathrm{A}$; la seule fibre invariante est $z=0$. Le feuilletage $z\frac{\partial}{\partial x}+(t-x)\frac{\partial}{\partial y}-z\frac{\partial}{\partial t}$ est invariant par~$\Phi_\mathrm{A}$. On constate que l'ensemble des matrices non inversibles est invariant par $\Phi_\mathrm{A}$, les $2$-plans $x=z=0$, $y=t=0$, $z=t=0$ aussi. Toute matrice de la forme $\left[\begin{array}{cc}\lambda&\mu\\ 0 &\lambda \end{array} \right]$ commute \`a $\mathrm{A};$ par cons\'equent $\left\{\left[\begin{array}{cc}\lambda&\mu\\ 0 &\lambda \end{array} \right]\,\big\vert\,\lambda,\,\mu\in\mathbb{C}\right\}$ est invariant par multiplication par $\mathrm{A}$ et par $\Phi_\mathrm{A}$.

\smallskip

\noindent Comme au \S \ref{centrdiag} on d\'emontre l'\'enonc\'e suivant.

\begin{pro}
Soit $\mathrm{A}$ la matrice $\left[\begin{array}{cc}1&1\\0&1\end{array} \right]$. Le groupe $\mathrm{Aut}(\mathcal{M}(2;\mathbb{C});\Phi_\mathrm{A})$ est engendr\'e par les $\sigma_\mathrm{P}$ o\`u $\mathrm{P}$ commute \`a $\mathrm{A}$.
\end{pro}

\subsubsection{\'Etude des points fixes et p\'eriodiques}

\noindent On peut v\'erifier que les points fixes de $\Phi_\mathrm{A}$ sont 
\begin{align*}
& {\bf 0}, && \left[\begin{array}{cc} 1 & -1\\ 0 & 1\end{array}\right], &&
\left\{\left[\begin{array}{cc} 1 & y\\ 0 & 0\end{array}\right]\,\Big\vert\, y\in\mathbb{C}\right\}.
\end{align*}

\medskip

\noindent D'apr\`es ce qui pr\'ec\`ede les points p\'eriodiques de $\Phi_{\mathrm{A}}$ sont contenus dans l'hyperplan $z=0.$ Un calcul montre que $$\Phi_\mathrm{A}^n\left(\left[\begin{array}{cc} x &y\\ 0&t \end{array}\right] \right)=\left[\begin{array}{cc} x^{2^n}&\Big(\displaystyle\prod_{i=0}^{n-1} \Big(x^{2^i}+t^{2^i}\Big)\Big)y+t^{2^n}+\sum_{k=1}^{n-1}\left(t^{2^k}\prod_{i=k}^{n-1}\Big(x^{2^i}+t^{2^i}\Big)\right) \\ & \\ 0&t^{2^n}\end{array}\right].$$ 

\noindent Si $(\underline{x},\underline{y},\underline{z},\underline{t})$ est un point p\'eriodique de $\Phi_\mathrm{A}$ de p\'eriode $n$, alors $\underline{z}=0$ et 
\begin{itemize}
\item[$\mathfrak{a.}$] ou bien $\underline{x}=\underline{t}=0;$

\item[$\mathfrak{b.}$] ou bien $\underline{x}^{2^n-1}=1$ et $\underline{t}=0;$

\item[$\mathfrak{c.}$] ou bien $\underline{x}=0$ et $\underline{t}^{2^n-1}=1;$

\item[$\mathfrak{d.}$] ou bien $\underline{x}^{2^n-1}=\underline{t}^{2^n-1}=1.$
\end{itemize}

\noindent Examinons ces \'eventualit\'es au cas par cas.

\noindent$\mathfrak{a.}$ Si $\underline{x}= \underline{t}=0$ alors n\'ecessairement $\underline{y}=0,$ \emph{i.e.} ${\bf 0}$ est p\'eriodique. 

\noindent$\mathfrak{b.}$ Si $\underline{x}^{2^n-1}=1$ et $\underline{t}=0$, on constate que $(\underline{x},\underline{y},0,0)$, $\underline{x}^{2^n-1}=1$, est p\'eriodique de p\'eriode $n$.

\noindent$\mathfrak{c.}$ Si $\underline{x}=0$ et $\underline{t}^{2^n-1}=1$, alors $\Phi_\mathrm{A}^n(0,\underline{y},0,\underline{t})=(0,\underline{y}+n\underline{t} , 0,\underline{t})$. 

\noindent$\mathfrak{d.}$ Pla\c{c}ons-nous maintenant dans le second cas: $\underline{x}^{2^n-1}= \underline{t}^{2^n-1}=1$. Nous allons raisonner suivant que $\underline{t}= \underline{x}$ et~$\underline{t}\not=\underline{x}$. 

\noindent Dans un premier temps supposons que $\underline{t}= \underline{x}$ auquel cas $\Phi_\mathrm{A}^n(\underline{x},\underline{y},0,\underline{x})=(\underline{x}^{2^n}, 2^n\underline{y}+(2^n-1)\underline{x},0,\underline{x}^{2^n})$; ces consid\'erations produisent les points p\'eriodiques suivants $\displaystyle\bigcup_{n\geq 0}\big\{(\underline{x},-\underline{x},0,\underline{x})\,\big\vert\,\underline{x}^{2^n-1}=1\big\}$.

\noindent Reste l'\'eventualit\'e $\underline{x}\not=\underline{t}.$ Posons 
\begin{align*}
& B_n(\underline{x},\underline{t})=\prod_{i=0}^{n-1}\Big(\underline{x}^{2^i}+\underline{t}^{2^i}\Big), && C_n(\underline{x},\underline{t})= \underline{t}^{2^n}+\sum_{k=1}^{n-1}t^{2^k}\left( \prod_{i=k}^{n-1}\Big(\underline{x}^{2^i}+ \underline{t}^{2^i}\Big)\right).
\end{align*}

\noindent Sous toutes ces notations et hypoth\`eses nous avons 
\begin{align*}
&\Phi_\mathrm{A}=(\underline{x}^{2^n},B_n(\underline{x},\underline{t})y+C_n(\underline{x},\underline{t}),0,\underline{t}^{2^n}), && \underline{x}^{2^n-1}= \underline{t}^{2^n-1}=1, &&\underline{x}\not=\underline{t}.
\end{align*}

\noindent Nous sommes donc ramen\'es \`a consid\'erer l'\'equation $B_n(\underline{x},\underline{t})y+C_n(\underline{x},\underline{t})=y$ o\`u $\underline{x},$ $\underline{t}$ sont des racines $(2^n-1)$-i\`eme de l'unit\'e distinctes. On peut v\'erifier que $B_n(\underline{x},\underline{t})=1$. Si pour chaque couple $(\underline{x},\underline{t})$ de racines $(2^n-1)$-i\`eme de l'unit\'e distinctes $C_n(\underline{x},\underline{t})$ est non nul, il n'y a pas de point p\'eriodique de p\'eriode $n$ de la forme $(\underline{x},y,0,\underline{t})$ avec $\underline{x}\not=\underline{t}$, $\underline{x}^{2^n-1}=\underline{t}^{2^n-1}=1$; sinon, pour tous les $(\underline{x},\underline{t})$ tels que 
\begin{align*}
& \underline{x}^{2^n-1}=\underline{t}^{2^n-1}=1,&& \underline{x}\not=\underline{t}, && C_n(\underline{x},\underline{t})=0,
\end{align*}
les $(\underline{x},y,0,\underline{t})$ sont des points p\'eriodiques de p\'eriode $n$.

\noindent Notons qu'il arrive que $C_n(\underline{x},\underline{t})$ soit non nul, par exemple lorsque $n=2$, $\underline{t}=\mathrm{j}$ et $\underline{x}=1$. 

\bigskip

\subsubsection{\'Etude de quelques orbites non p\'eriodiques}

\noindent Comme pr\'ec\'edemment nous avons $\det\Phi_\mathrm{A}(\mathrm{M})= (\det\mathrm{M})^2$ et donc $\det\Phi_\mathrm{A}^k(\mathrm{M})= (\det\mathrm{M})^{2^k}$. Ainsi, pour tout $\mathrm{M}$ satisfaisant $\vert\det\mathrm{M}\vert>1$, nous avons $\displaystyle\lim_{k\to+\infty}\vert\vert\Phi_\mathrm{A}^k(\mathrm{M})\vert\vert=~+\infty$. 

\noindent Comme dans le cas diagonal, si $\mathrm{M}$ appartient au polydisque $\Delta(\rho),$ alors $\vert\vert\Phi_\mathrm{A}(\mathrm{M}) \vert\vert\leq 2\rho^2$ et $\vert\vert\Phi_\mathrm{A}^k(\mathrm{M})\vert\vert\leq  2^{2^k-1}\rho^{2^k}$ qui entra\^ine l'inclusion $$\Delta\left(\frac{1}{2}\right)\subset\mathrm{W}_\mathrm{A}^s ({\bf 0}).$$

\noindent Donnons quelques exemples d'orbites born\'ees. Remarquons que $\Phi_\mathrm{A}^n(x,y,0,x)=x^{2^n-1}\big(x,(2^n-1)x+2^ny, 0, x\big);$ par suite, d\`es que $\mathrm{M}$ est de la forme 
\begin{align*}
& \left[\begin{array}{cc} x & -x \\ 0 & x\end{array}\right] \text{ avec } \vert x\vert<1 &&\text{ ou }&&\left[\begin{array}{cc} x & y \\ 0 & x\end{array}\right] \text{ avec } \vert x\vert<1 \text{ et } (2^n-1)x +2^ny=0 \text{ pour un certain $n$}
\end{align*}
l'orbite de $\mathrm{M}$ est born\'ee.

\noindent Puisque $\Phi_\mathrm{A}^n(0, y, 0, t)=t^{2^n-1}\big(0, nt+y, 0, t\big),$ l'orbite de $\left[\begin{array}{cc} 0 & y \\ 0 & t\end{array}\right],$ avec $\vert t\vert<1$ et $y=-nt$ pour un certain entier $n,$ est born\'ee.

\noindent \'Etant donn\'e que $\Phi_\mathrm{A}^n(x, y, 0, 0)=x^{2^n-1}\big(x, y, 0, 0),$ toute matrice $\left[\begin{array}{cc} x & y \\ 0 & 0\end{array}\right]$ avec $\vert x\vert<1$ a une orbite born\'ee.

\begin{probs}
\begin{itemize}
\item[1)] D\'ecrire les points p\'eriodiques de l'application $\Phi_\mathrm{A}$ et leur adh\'erence.

\item[2)] D\'ecrire le bassin d'attraction $\mathrm{W}_\mathrm{A}^s({\bf 0})$ et son bord.
\end{itemize}
\end{probs}

\vspace{8mm}

\bibliographystyle{plain}
\bibliography{iteration}

\begin{thebibliography}{1}

\bibitem{BS}
E.~Bedford and J.~Smillie.
\newblock Fatou-{B}ieberbach domains arising from polynomial automorphisms.
\newblock {\em Indiana Univ. Math. J.}, 40(2):789--792, 1991.

\bibitem{BL}
F.~Berteloot and J.~J. L{\oe}b.
\newblock Holomorphic equivalence between basins of attraction in {$\bold
  C^2$}.
\newblock {\em Indiana Univ. Math. J.}, 45(2):583--589, 1996.

\bibitem{CLN}
D.~Cerveau and A.~Lins~Neto.
\newblock Hypersurfaces exceptionnelles des endomorphismes de {${\bf C}{\rm
  P}(n)$}.
\newblock {\em Bol. Soc. Brasil. Mat. (N.S.)}, 31(2):155--161, 2000.

\bibitem{CoEc}
P.~Collet and J.-P. Eckmann.
\newblock {\em Iterated maps on the interval as dynamical systems}.
\newblock Modern Birkh\"auser Classics. Birkh\"auser Boston Inc., Boston, MA,
  2009.
\newblock Reprint of the 1980 edition.

\bibitem{Fi}
G.~Fischer.
\newblock {\em Complex analytic geometry}.
\newblock Lecture Notes in Mathematics, Vol. 538. Springer-Verlag, Berlin,
  1976.

\bibitem{FM}
S.~Friedland and J.~Milnor.
\newblock Dynamical properties of plane polynomial automorphisms.
\newblock {\em Ergodic Theory Dynam. Systems}, 9(1):67--99, 1989.

\bibitem{GV}
R.~Galeeva and A.~Verjovsky.
\newblock Quaternion dynamics and fractals in {${\bf R}^4$}.
\newblock {\em Internat. J. Bifur. Chaos Appl. Sci. Engrg.}, 9(9):1771--1775,
  1999.
\newblock Discrete dynamical systems.

\bibitem{Ma}
B.~Malgrange.
\newblock Frobenius avec singularit\'es. {II}. {L}e cas g\'en\'eral.
\newblock {\em Invent. Math.}, 39(1):67--89, 1977.

\bibitem{Mi}
J.~Milnor.
\newblock Pasting together {J}ulia sets: a worked out example of mating.
\newblock {\em Experiment. Math.}, 13(1):55--92, 2004.

\end{thebibliography}
\nocite{}

\end{document}